\documentclass[11pt,a4paper]{article}

\usepackage{comment}
\usepackage{graphicx}
\usepackage{amsmath}
\usepackage{amsfonts}
\usepackage{amssymb}
\usepackage{enumerate}
\usepackage{lscape}
\usepackage{longtable}
\usepackage{rotating}
\usepackage{color}
\usepackage{url}
\usepackage{rotating}
\usepackage{algpseudocode}
\usepackage[ruled]{algorithm}
\usepackage{pgf,tikz}
\usepackage{subfig}
\usepackage[titletoc]{appendix}
\usepackage{wrapfig}
\usepackage{adjustbox}
\usepackage{mathrsfs}
\usepackage{pgfplots}
\usetikzlibrary{arrows}
\usepackage{color}

\usepackage{bbm}

\usepackage{booktabs}
\usepackage{multirow}
\usepackage{mathtools}
\usepackage{makecell}

\usepackage[sort,numbers]{natbib}
\usepackage{eqparbox}%

\numberwithin{equation}{section}%

\usepackage[symbol]{footmisc}%

\newtheorem{theorem}{Theorem}[section]
\newtheorem{proposition}{Proposition}[section]

\newtheorem{lemma}[theorem]{Lemma}

\newtheorem{remark}{Remark}[section]

\newtheorem{assumption}{Assumption}[section]

\newenvironment{proof}[1][Proof]{\noindent \textbf{#1.} }{\hfill$\Box$\par\medskip}

\DeclareMathOperator*{\argmin}{argmin}

\DeclareMathOperator{\expc}{\mathbb{E}}

\DeclareMathOperator{\train}{train}
\DeclareMathOperator{\val}{val}

\DeclareMathOperator{\tr}{tr}
\DeclareMathOperator{\va}{v}

\DeclareMathOperator{\BSG}{BSG}
\DeclareMathOperator{\DARTS}{DARTS}

\DeclareMathOperator{\sone}{S1}

\DeclareMathOperator{\all}{all}

\newcommand{\beqn}[1]{\begin{equation}\label{#1}}
\newcommand{\eeqn}{\end{equation}}

\definecolor{darkgreen}{rgb}{0,0.6,0}
\definecolor{aau2}{rgb}{0.0, 0.5, 0.69}
\definecolor{aau3}{rgb}{0.0, 0.53, 0.74}
\definecolor{aau4}{rgb}{0.0, 0.48, 0.65}
\definecolor{aau5}{rgb}{0.0, 0.45, 0.73}
\definecolor{rsap}{RGB}{130, 36, 51}
\definecolor{gsap}{RGB}{112, 164, 137}

\definecolor{tud}{rgb}{0.43,0.73,0.11}
\definecolor{verde}{rgb}{0.33,0.53,0.11}

\definecolor{ttffqq}{rgb}{0.0, 0.48, 0.65} 
\definecolor{ffqqqq}{rgb}{0.0, 0.5, 0.69} 

\usetikzlibrary{arrows}

\tikzstyle{decision} = [diamond, draw, fill=blue!20,
text width=4.5em, text badly centered, node distance=3cm, inner sep=0pt]
\tikzstyle{block} = [rectangle, draw, fill=blue!20,
text centered, rounded corners, minimum height=4em]
\tikzstyle{line} = [draw, -latex']
\tikzstyle{cloud} = [draw, ellipse,fill=red!20, node distance=3cm,
minimum height=2em]
\tikzstyle{cloud2} = [draw, ellipse,fill=green!20, node distance=3cm,
minimum height=2em]
\begin{document}
	
	\title{Inexact bilevel stochastic gradient methods for constrained and unconstrained lower-level problems}
	
	\author{
		T. Giovannelli\thanks{Department of Industrial and Systems Engineering, Lehigh University, Bethlehem, PA 18015-1582, USA ({\tt tog220@lehigh.edu}).}
		\and
		G. D. Kent\thanks{Department of Industrial and Systems Engineering, Lehigh University, Bethlehem, PA 18015-1582, USA ({\tt gdk220@lehigh.edu}).}
		\and
		L. N. Vicente\thanks{Department of Industrial and Systems Engineering, Lehigh University, Bethlehem, PA 18015-1582, USA ({\tt lnv@lehigh.edu}).}
	}
	
	\maketitle
	
	\begin{abstract}
	Two-level stochastic optimization formulations have become instrumental in a number of machine learning contexts such as continual learning, neural architecture search, adversarial learning, and hyperparameter tuning.  Practical stochastic bilevel optimization problems become challenging in optimization or learning scenarios where the number of variables is high or there are constraints.
	
    In this paper, we introduce a bilevel stochastic gradient method for bilevel problems with nonlinear and possibly nonconvex lower-level constraints. We also present a comprehensive convergence theory that addresses both the lower-level unconstrained and constrained cases and covers all inexact calculations of the adjoint gradient (also called hypergradient), such as the inexact solution of the lower-level problem, inexact computation of the adjoint formula (due to the inexact solution of the adjoint equation or use of a truncated Neumann series), and noisy estimates of the gradients, Hessians, and Jacobians involved. To promote the use of bilevel optimization in large-scale learning, we have developed new low-rank practical bilevel stochastic gradient methods (BSG-N-FD and~BSG-1) that do not require second-order derivatives and, in the lower-level unconstrained case, dismiss any matrix-vector products.  
	\end{abstract}
	
	\section{Introduction}\label{sec:introduction}

Many real-world applications are naturally formulated using hierarchical objectives, which are organized into different nested levels. In the bilevel case, the main goal is placed into an upper optimization level, while the lower optimization level aims to determine the best response to a decision made in the upper level. Bilevel optimization has a rich literature of algorithmic development and theory (see \cite{LNVicente_PHCalamai_1994,BColson_PMarcotte_GSavard_2007,JFBard_2010,ASinha_PMalo_KDeb_2018,SDempe_AZemkoho_2020,SDempe_2002} for extensive surveys and books on this topic). The main applications are found in game theory, defense industry, and optimal structural design, and one has recently seen a surge of contributions to machine learning~(ML)
(see, e.g., \cite{LFranceschi_etal_2018,HLiu_KSimonyan_YYang_2019}, and the recent review \cite{RLiu_JGao_etal_2021}). 

In this paper, we consider the following nonlinear bilevel optimization problem (BLP) formulation, where we are using a standard notation (see, e.g.,~\cite{LNVicente_PHCalamai_1994,SDempe_2002,ASinha_PMalo_KDeb_2018})
	\begin{equation}\label{prob:bilevel}
	\tag*{BLP}
	\begin{split}
	\min_{x \in \mathbb{R}^n, \, y \in \mathbb{R}^m} ~~ & f_u(x,y) \\
	\mbox{s.t.}~~ & x \in X \\
 & y \in \argmin_{y \in Y(x)} ~~ f_\ell(x,y).\\
	\end{split}
	\end{equation}
The goal of the upper-level~(UL) problem is to determine the optimal value of the UL function $f_u: \mathbb{R}^n\times\mathbb{R}^m \to \mathbb{R}$, where the UL variables~$x$ are subjected to UL constraints~($x \in X$) and the lower-level~(LL) variables~$y$ are subjected to being an optimal solution of the~LL problem. In the LL problem, the LL function $f_\ell: \mathbb{R}^n\times\mathbb{R}^m \to \mathbb{R}$ is optimized in the LL variables~$y$, subject to the~LL constraints~$y \in Y(x)$.
We will state the assumptions required for stochastic gradient descent in Subsection~\ref{subsec:general_assumptions}.
We will assume~$f_u$ to be continuously differentiable in~$(x, y)$ and~$f_{\ell}$ to be twice continuously differentiable in~$(x, y)$.
We will also assume the LL problem to be well-defined, in the sense that an~LL optimal solution~$y(x)$ exists and is unique for all $x \in X$.
Hence, problem~\ref{prob:bilevel} is equivalent to a problem posed solely in the UL variables:
\begin{equation} \label{reduced}
\min_{x \in \mathbb{R}^n} \; f(x)=f_u(x,y(x)) \quad
	\mbox{s.t.} \quad x \in X.
\end{equation}
Also, note that the UL constraints ($x \in X$) are only posed in the UL variables~$x$ as otherwise problem~\ref{prob:bilevel} could become intractable due to a disconnected feasible region in the $(x,y)$--space.
The set~$X$ will be assumed closed and convex, which will allow us to ensure~UL feasibility by applying orthogonal projections within stochastic gradient type methods.

Assuming~$\nabla^2_{yy}f_\ell(x,y(x))$ is non-singular (again, see Subsection~\ref{subsec:general_assumptions} for the statement of the assumptions), the gradient of~$f$ at~$x$, when $Y(x) = \mathbb{R}^m$, is given by the well-known so-called adjoint gradient (also called hypergradient in the~ML community)
\begin{equation} \label{adjoint}
\nabla f \; = \; \nabla_x f_u - \nabla_{xy}^2 f_\ell \nabla^2_{yy}f_\ell^{-1} \nabla_y f_u,
\end{equation}
where all gradients and Hessians on the right-hand side are evaluated at~$(x,y(x))$. We denote the steepest descent direction for~$f$ at~$x$ as~$d(x, y(x)) = - \nabla f(x)$. One arrives at the adjoint formula by first applying the chain rule to $f_u(x,y(x))$ to obtain
\begin{equation} \label{adjoint_jacobian}
\nabla f \; = \; \nabla_x f_u + \nabla y \nabla_y f_u.
\end{equation}
Then, the Jacobian~$\nabla y(x)^{\top}$ of~$y(x)$ can be calculated through the (sensitivity) \linebreak equations~$\nabla_y f_\ell(x,y(x))=0$. The implicit function theorem ensures $y(\cdot)$ to be continuously differentiable~\cite{Ruding_1953}. By taking the derivative of both sides of the equation with respect to $x$ and utilizing the chain rule, we obtain~$\nabla^2_{yx}f_\ell + \nabla^2_{yy}f_\ell\nabla y^{\top} = 0$ (all Hessians are evaluated at~$(x,y(x))$) which yields 
\begin{equation}\label{eq:jacobian_formula}
\nabla y \; = \; - \nabla^2_{xy}f_\ell \nabla^2_{yy}f_\ell^{-1}.
\end{equation}

Two approaches have been proposed in the literature to deal with~$\nabla^2_{yy}f_\ell^{-1}$ in~\eqref{adjoint}. One option is to compute the adjoint gradient by first solving the linear system given by the adjoint equation~$\nabla_{yy}^2 f_\ell \, \lambda = \nabla_y f_u$ for the adjoint variables~$\lambda=\lambda(x,y(x))$, and then calculating
$\nabla_x f_u - \nabla_{xy}^2 f_\ell \, \lambda$. The second option is to truncate the Neumann series given by~$\nabla^2_{yy}f_\ell^{-1} = \sum_{i=0}^{\infty} (I - \nabla_{yy}^2 f_\ell)^i$, which, however, requires either the strong assumption of~$\Vert \nabla_{yy}^2 f_\ell \Vert_2 < 1$ or the knowledge of a bound on the second derivatives to guarantee the convergence of the series.

When $Y(x) \neq \mathbb{R}^m$, it is still possible to use an adjoint formula to compute the gradient of~$f$ at~$x$ by using sensitivity arguments from nonlinear programming. Such an~LL constrained case will be addressed in Section~\ref{sec:constrained}.

\subsection{Bilevel machine learning}\label{sec:BML}

A variety of problems arising in machine learning can be formulated in terms of bilevel optimization. Continual learning, neural architecture search, adversarial training, and hyperparameter tuning are among the most popular examples (see~\cite{RLiu_JGao_etal_2021} for a review on this topic). 

Continual Learning (CL) aims to train ML models when the static task usually considered in learning problems (classification, regression, etc.)~is replaced by a sequence of tasks that become available one at a time~\cite{DLopezPaz_MARanzato_2017}, and for which training and validation datasets are increasingly larger.
For each task, a CL instance is formulated as a bilevel problem, where at the UL problem one minimizes the validation error on a subset of model parameters (which includes all hyperparameters),
and at the LL problem the training error is minimized on the remaining parameters. Then, a sequence of bilevel problems (one for each task) is solved. In a sense, CL is close to meta-learning~\cite{THospedales_etal_2020}, where the goal is to determine the best learning process. The increasing interest in CL is motivated by the demand for approaches that help neural networks learn new tasks without forgetting the previous ones, a phenomenon which is referred to as \textit{catastrophic forgetting}~\cite{IJGoodfellow_MMirza_DXiao_ACourville_YBengio_2013,RFrench_1999,MMcCloskey_NJCohen_1989,HSun_etal_2022}. Another relevant~ML area where bilevel optimization is used is Neural Architecture Search (NAS) for Deep Learning. The goal of this problem is to automate the task of designing Deep Neural Networks (DNNs) such that the network’s prediction error is minimized. In recent years, NAS has been proposed in a bilevel optimization formulation~\cite{HLiu_KSimonyan_YYang_2019}.

Finally, two other popular classes of ML applications that can be formulated by using bilevel optimization are adversarial training and hyperparameter tuning. Adversarial training aims to robustly address adversarial examples~\cite{CSzegedy_etal_2013} which cannot be correctly classified by ML models once a small perturbation is applied. The adversarial training problem is handled by solving a min-max problem~\cite{HJiang_etal_2018}, which can be reformulated as a bilevel one. The max/LL problem is posed on the variables that perturb the data in a worst-case fashion, where the min/UL problem attempts to minimize the training error on the ML model parameters~\cite{IGoodfellow_JShlens_CSzegedy_2015,AMadry_etall_2018}. 
Hyperparameter tuning aims to find the best values for the hyperparameters used in an ML model in order to increase its performance on unseen data~\cite{KPBennett_etall_2008,LFranceschi_etal_2018,SDhar_UKurup_MShah_2020,CChen_etal_2021}. In the bilevel formulations proposed in the literature, the UL problem optimizes the validation error over the hyperparameters, while the LL problem has the goal of finding the ML model parameters (e.g., neural network weights) that minimize the training error.     

\subsection{Bilevel stochastic descent}\label{sec:BSD}	

In bilevel stochastic optimization, $f_u$ and $f_\ell$ can be interpreted as expected values, namely, $f_u = \mathbb{E}[f_u(x,y,\vartheta^u)]$ and $f_\ell = \mathbb{E}[f_\ell(x,y,\vartheta^\ell)]$, where $\vartheta^u$ and $\vartheta^\ell$ are random variables defined in a probability space (with probability measure independent from $x$ and $y$) such that i.i.d. samples can be observed or generated.  
The same applies to the functions possibly defining $Y(x)$.
(To keep the notation simple, we are using the same $f_u$ and $f_\ell$ for deterministic and random variants.)

Having in mind~ML applications, the methods that we are considering are Stochastic Approximation (SA) techniques, of the type of the stochastic gradient (SG) method \cite{HRobbins_SMonro_1951,KLChung_1954,JSacks_1958} for single-objective optimization. In fact, the bilevel stochastic gradient (BSG) method can be seen as an~SG method applied to~(\ref{reduced}), which leads to $x_{k+1}=x_k-\alpha_k g_k^{\BSG}$, where
$\alpha_k$ is the step size or learning rate and $g_k^{\BSG}$ is a stochastic gradient of~$f$. Such a stochastic gradient is obtained by sampling the gradients and Hessians in~(\ref{adjoint}) at $(x_k,\tilde{y}_k)$, with~$\tilde{y}_k$ denoting an approximation to the~LL optimal solution. The stochastic gradient~$g_k^{\BSG}$ is inexact when $\tilde{y}_k \neq y(x_k)$, even in the full-batch (deterministic) case.

 
In general,~BSG methods have mainly been considered for the~LL unconstrained case~(i.e., $Y(x) = \mathbb{R}^m$) and are commonly classified according to the approach used to compute the~BSG direction~$g_k^{\BSG}$~\cite{RLiu_JGao_etal_2021}. In particular, a first category, referred to as \textit{implicit differentiation}, is composed of algorithms that compute the~BSG direction by applying the implicit function theorem and either solving the adjoint equation~\cite{FPedregosa_2016} or using a truncated Neumann series to approximate the inverse Hessian~$\nabla_{yy}^2 f_\ell^{-1}$~\cite{JLorraine_PVicol_DDuvenaud_2019}. Note that using a truncated Neumann series requires~$\Vert \nabla_{yy}^2 f_\ell\Vert_2 < 1$, which is a strong assumption. Therefore, a common approach is to first assume~$\nabla_y f_{\ell}$ Lipschitz continuous in~$y$ with constant~$C_0$, and then apply the truncated Neumann series to approximate~$[(1/C_0)\nabla_{yy}^2 f_\ell]^{-1}$. However, this requires the knowledge of~$C_0$, which is typically unknown in practice. A second category, referred to as \textit{iterative differentiation}, includes all the approaches based on automatic differentiation through dynamic systems~\cite{EWeinan_2017,YLu_AZhong_QLi_BDong_2017,LFranceschi_etal_2018}. 
A general convergence theory for the two classes of algorithms was proposed in~\cite{KJi_JYang_YLiang_2020}, which also shows that computing the~BSG direction by using automatic differentiation can be less computationally efficient than using an implicit differentiation method. Therefore, our paper focuses on the first category, which has been promoted in~\cite{NCouellan_WWang_2015,NCouellan_WWang_2016,SGhadimi_MWang_2018,MHong_etal_2020,TChen_YSun_WYin_2021_closingGap,DSow_KJi_YLiang_2021} (see also~\cite{RLiu_JGao_etal_2021,CChen_etal_2022} for recent reviews).
These existing approaches either focus on a specific problem structure or require the LL problem to be solved to optimality at each iteration or rely on a truncated Neumann series for the inverse Hessian approximation in the adjoint gradient, which has the issues mentioned before. Among all the algorithms proposed in the papers cited above, we emphasize StocBiO~\cite{KJi_JYang_YLiang_2020}. StocBiO employs a truncated Neumann series with automatic differentiation and a double-loop iterative scheme, which means that multiple iterations are required at the~LL problem in order for the algorithm to converge, similar to the algorithms developed in our paper.

DARTS~\cite{HLiu_KSimonyan_YYang_2019} is an optimization technique related to the BSG method and has enjoyed great popularity in NAS.
It always considers an inexact solution to the LL problem, and it starts an iteration by applying one step of SG to the LL problem,
$\tilde y_k = y_k - \eta \nabla_y f_\ell(x_k,y_k)$, where $\eta$ is a fixed step size. Then, 
it displaces the UL variables using
$x_{k+1}=x_k - \eta \, g_k^{\DARTS}$, where $g_k^{\DARTS}$ is computed by applying the chain rule to $\nabla_x f_u(x,y - \eta \nabla_y f_\ell(x,y))$, leading to (when using full-batch gradients)
\begin{equation} \label{eq:DARTS}
g_k^{\DARTS} \; = \; \nabla_x f_u(x_k, \tilde y_k) - \eta \nabla^2_{xy} f_\ell (x_k, y_k) \nabla_{y} f_u (x_k, \tilde y_k).
\end{equation}
However, in the full-batch (deterministic) case,~$g_k^{\DARTS}$ may not be a descent direction.
The matrix-vector product in~(\ref{eq:DARTS}) is approximated by finite differences, rendering DARTS free of both second-order derivatives and matrix-vector products (see Section~\ref{sec:stochasticBLP_DARTS}).

All the papers cited above focus on the~LL unconstrained case. A few approaches dealing with the~LL constrained case have been very recently proposed to tackle bilevel problems with~LL constraints. However, all these approaches are only applicable to linear constraints or constraints only depending on~$y$~\cite{ITsaknakis_PKhanduri_MHong_2022,PKhanduri_etal_2023,QXiao_etal_2023,HShen_QXiao_TChen_2023,JKwon_DKwon_SWright_etal_2023}. Among the algorithms proposed, we highlight~SIGD~\cite{PKhanduri_etal_2023}, which determines a direction by applying the implicit function theorem (similar to the approach that we use) but can only be applied to bilevel problems with linear~LL constraints in~$y$.


\subsection{Contributions of the paper} \label{sec:cont}

The first main contribution of this paper is a general framework for the~BSG method that applies to both the~LL unconstrained ($Y(x) = \mathbb{R}^m$) and constrained ($Y(x)\neq \mathbb{R}^m$) cases. In particular, our paper represents the first work that proposes a method to address the general nonlinear and possibly nonconvex~LL constrained case, which has not been covered elsewhere, neither algorithmically nor theoretically, although important~ML applications give rise to bilevel problems with~LL constraints.

The second main contribution is a comprehensive convergence theory for the~BSG method that applies to both the~LL unconstrained and constrained cases and is grounded on sensitivity principles of nonlinear optimization. The sensitivity arguments used in this paper are easily satisfied in all practical scenarios that have been proposed for~ML applications requiring bilevel optimization. Our theory also comprehensively covers all possible inexact settings such as the inexact solution of the~LL problem, inexact computation of the adjoint formula (due to the inexact solution of the adjoint equation or use of a truncated Neumann series), and noisy estimates of the gradients, Hessians, and Jacobians involved.  
Since the convergence analysis proposed is abstracted from the specifics of the approach used to handle the inverse matrix in the adjoint formula~\eqref{adjoint}, our theory unifies two different classes of~BSG methods, i.e., the ones based on the adjoint equation and the ones based on the truncated Neumann series, which have been studied separately~\cite{SGhadimi_MWang_2018,FPedregosa_2016}.

Moreover, to deal with the second-order derivatives and inverse Hessians/Jacobians in the adjoint formulas for the~LL unconstrained and constrained cases, we developed two new low-rank practical implementations of the~BSG method that can be applied to solve large-scale optimization problems arising in~ML applications. The first one, referred to as BSG-N-FD, consists of solving the adjoint equation by using an iterative method (CG or~GMRES) equipped with finite differences to dismiss any Hessian-vector products, and is grounded on theoretical principles. The second one, referred to as~BSG-1, uses rank-1 Hessian approximations to avoid explicitly solving the adjoint equation (also dismissing any matrix-vector products). The use of these rank-1 approximations is inspired by Gauss-Newton methods for nonlinear least-squares problems~\cite{JNocedal_SJWright_2006} and also from the fact that the empirical risk of misclassification in ML is often a sum of non-negative terms, matching a function to a scalar which can then be considered in a least-squares fashion~\cite{ABotev_2017,MGargiani_etal_2020}. Similar to DARTS~\cite{HLiu_KSimonyan_YYang_2019} (which is extremely popular for~NAS), \linebreak BSG-1 is also not grounded on theoretical principles. However, both methods perform well on continual learning instances in terms of training iterations and computational time.

\subsection{Organization of this paper}

	This paper is organized as follows. In Section~\ref{sec:BSGmethod}, we describe the~BSG method for the~LL unconstrained case and introduce it for the~LL constrained case. The assumptions on the problem functions and inexact calculations required for the convergence analysis of the method are reported in Section~\ref{sec:assumptions_BSGmethod}. Section~\ref{sec:convergence_rate_BSGmethod} presents the convergence rates for the nonconvex, strongly convex, and convex cases. Numerical results for synthetic quadratic bilevel problems and continual learning instances with or without~LL constraints are analyzed in Section~\ref{sec:numerical_experiments}, which also describes the practical~BSG-N-FD and~BSG-1 algorithms. Finally, in Section~\ref{sec:conclusions} we draw some concluding remarks and propose ideas for future work. By default, all norms~$\| \cdot \|$ used in this paper are the~$\ell_2$~ones.
	
\section{The bilevel stochastic gradient method}\label{sec:BSGmethod}

In this section, we introduce the bilevel stochastic gradient~(BSG) method for solving stochastic~BLPs. Let~$\{\vartheta^\ell_k\}_{k \ge 0}$ and~$\{\varsigma^\ell_k\}_{k \ge 0}$ be sequences of random variables for~LL gradient, Hessian, and Jacobian evaluations for the~LL unconstrained and constrained cases, respectively. Similarly, let~$\{\vartheta^u_k\}_{k \ge 0}$ be a sequence of random variables for~UL gradient evaluations. A realization of~$\vartheta^\ell_k$, $\varsigma^\ell_k$, and~$\vartheta^u_k$ can be interpreted as a single sample or a batch of samples for mini-batch~SG. For compactness of notation, let us set~$\xi_k=(\vartheta_k^u,\vartheta_k^\ell)$ in the~LL unconstrained case ($Y(x) = \mathbb{R}^m$), and~$\xi_k=(\vartheta_k^u,\vartheta_k^\ell,\varsigma_k^\ell)$ in the~LL constrained case ($Y(x) \neq \mathbb{R}^m$). All the assumptions included in this section will be rigorously stated in Subsection~\ref{subsec:general_assumptions}.

\subsection{The unconstrained lower-level case} \label{sec:unconstrained}

Given $(x_k,\tilde{y}_k)$, we denote by $g^u_x(x_k,\tilde{y}_k,\vartheta^u_k)$, $g^u_y(x_k,\tilde{y}_k,\vartheta^u_k)$, and $g^\ell_y(x_k,\tilde{y}_k,\vartheta^\ell_k)$ the stochastic gradient estimates that approximate $\nabla_x f_u(x_k,\tilde{y}_k)$, $\nabla_y f_u(x_k,\tilde{y}_k)$, and $\nabla_y f_\ell(x_k,\tilde{y}_k)$, respectively. The same notation applies to the stochastic Hessian estimates:
$H^\ell_{xy}(x_k,\tilde{y}_k,\vartheta^\ell_k)$ and $H^\ell_{yy}(x_k,\tilde{y}_k,\vartheta^\ell_k)$ approximate
$\nabla^2_{xy} f_\ell(x_k,\tilde{y}_k)$ and 
$\nabla^2_{yy} f_\ell(x_k,\tilde{y}_k)$, respectively. Based on Assumption~\ref{ass:cont_diff}, which will be stated in Subsection~\ref{subsubsec:LL_unc_case}, $\nabla^2_{yy} f_{\ell}^{-1}$ and~$(H^{\ell}_{yy})^{-1}$ are non-singular at all points.	
In the LL unconstrained case ($Y(x) = \mathbb{R}^m$), 
an approximate (negative)~BSG (denoted as~$-g_k^{\BSG}$ in Subsection~\ref{sec:BSD}) can be computed directly from the adjoint formula~(\ref{adjoint}), as follows:
\begin{equation} \label{stochastic-adjoint}
d(x_k, \tilde{y}_k, \xi_k) \; = \; -\left(
 g^u_x(x_k,\tilde{y}_k,\vartheta^u_k) - H^{\ell}_{xy}(x_k,\tilde{y}_k,\vartheta^\ell_k) H^{\ell}_{yy}(x_k,\tilde{y}_k,\vartheta^\ell_k)^{-1} g^u_y(x_k,\tilde{y}_k,\vartheta^u_k)\right).
\end{equation}

The data in the formula~(\ref{adjoint}) is referred to as~$D(x,y)$ or~$D(x,y(x))$, depending on the point where the gradients and Hessians are evaluated:
\begin{equation} \label{eq:D(x,y)}
D(x, y) \; = \;
\left( \nabla_x f_u(x,y), \nabla_y f_u(x,y), \nabla_{xy}^2 f_\ell(x,y), \nabla_{yy}^2 f_\ell(x,y) \right).
\end{equation}
The data in the calculation~(\ref{stochastic-adjoint}) is referred to as~$D(x_k,\tilde{y}_k,\xi_k)$:
\begin{equation} \label{eq:D(x,y,epsilon)}
D(x_k,\tilde{y}_k,\xi_k) \; = \;
\left( g^u_x(x_k,\tilde{y}_k,\vartheta^u_k), g^u_y(x_k,\tilde{y}_k,\vartheta^u_k), H^{\ell}_{xy}(x_k,\tilde{y}_k,\vartheta^\ell_k), H^{\ell}_{yy}(x_k,\tilde{y}_k,\vartheta^\ell_k) \right).
\end{equation}
Note that in~\eqref{stochastic-adjoint}, $d(x_k, \tilde{y}_k, \xi_k)$ can be seen as a function of the data~$D(x_k,\tilde{y}_k,\xi_k)$ as follows:
\begin{equation}\label{eq:dir_notation}
d(x_k, \tilde{y}_k, \xi_k) \; = \; d(D(x_k,\tilde{y}_k,\xi_k)),
\end{equation}
where we are using overlapping notation for~$d(\cdot)$ for the sake of simplicity and because we believe it is more powerful and effective than using a different letter.

As mentioned in Subsection~\ref{sec:cont}, our practical implementation~BSG-N-FD will solve the adjoint equation~$H^{\ell}_{yy}(x_k,\tilde{y}_k,\vartheta^\ell_k) \lambda = g^u_y(x_k,\tilde{y}_k,\vartheta^u_k)$ by using an iterative method equipped with finite differences (see Subsection~\ref{subsubsec:BSG-N-FD}). The~BSG-1 method will use rank-one approximations for~$H^{\ell}_{xy}(x_k,\tilde{y}_k,\vartheta^\ell_k)$ and~$H^{\ell}_{yy}(x_k,\tilde{y}_k,\vartheta^\ell_k)$, and then solve the resulting adjoint equation in the least-squares sense (see Subsection~\ref{subsubsec:BSG-1}).
    
\subsection{The constrained lower-level case} \label{sec:constrained}

Let us now handle the LL constrained case, in which we consider
\[
Y(x) \; = \; \{y \in \mathbb{R}^m ~|~ c_i(x,y) \leq 0, \; i \in I,\text{ and } c_i(x,y) = 0, \; i \in E\},
\]
where~$I$ and~$E$ are two finite sets of indices. As stated in the assumptions of Subsection~\ref{subsubsec:LL_constr_case}, each constraint function~$c_i$ is assumed twice continuously differentiable in~$(x,y)$, for all~$i \in I \cup E$.
Denoting~$c_I(x,y) = (c_i(x,y), \, i \in I)$ and~$c_E(x,y) = (c_i(x,y), \, i \in E)$, the Lagrangian function of the~LL problem is defined as~$\mathcal{L}_\ell(x,y,z) = f_\ell(x,y) + c_I(x,y)^\top z_I + c_E(x,y)^\top z_E$, where~$z_I$ and $z_E$ are Lagrange multipliers and $z=(z_I,z_E)$. We will also assume in Assumption~\ref{ass:LL_assumptions_constr} of Subsection~\ref{subsubsec:LL_constr_case} that there exists a~$y(x)$ satisfying the LL KKT conditions with associated multipliers~$(z_I(x), z_E(x))$ such that the gradients of the active constraints are linearly independent~(LICQ), the strict complementarity slackness condition~(SCS) is satisfied, and the second-order sufficient optimality conditions~(SOSC) hold. 
Under such assumptions, it is well known that the Lagrange multipliers~$z_I(x)$ and~$z_E(x)$ associated with~$y(x)$ are unique and the vector function~$v(x) = (y(x), z_I(x), z_E(x))^{\top}$, for any given $x$, is once continuously differentiable~\cite{AFiacco_GPMcCormick_1968,AFiacco_1976,AFiacco_1983,GPMcCormick_2014}.
Moreover, we can write the first-order~KKT system for the~LL problem (see~\cite{AFiacco_1976}) as
\begin{equation}\label{eq:KKT_syst}
    \begin{cases}
        \nabla_y f_\ell(x,y(x)) + \nabla_y c_I(x,y(x)) \, z_I(x) + \nabla_y c_E(x,y(x)) \, z_E(x) = 0,\\
        z_{I}(x) \circ c_{I}(x,y(x)) = 0,\\
        c_E(x,y(x)) = 0,
    \end{cases}
\end{equation}
where $\circ$ is the element-wise multiplication operation of two vectors.

Now, for any given $x$, we can rewrite the~KKT system~\eqref{eq:KKT_syst} as~$G(x, v(x)) = 0$ by introducing a corresponding vector function~$G$. 
Applying the chain rule to~$G(x, v(x)) = 0$, we obtain~$\nabla_v G^{\top} \nabla v^{\top} = - \nabla_x G^{\top}$, with 
\begin{equation}\label{eq:jacobians}
\nabla_x G^{\top} = \begin{pmatrix} \nabla_{yx}^2 \mathcal{L}_\ell \\  z_I \circ \nabla_x c_I^{\top} \\ \nabla_x c_E^{\top}  \end{pmatrix}
\; \text{ and } \;
\nabla_v G^{\top} = \begin{pmatrix} \nabla_{yy}^2 \mathcal{L}_\ell && \nabla_y c_I && \nabla_y c_E \\
z_I \circ \nabla_y c_I^{\top} && C_{I} && 0 \\
\nabla_y c_E^{\top} && 0 && 0 \end{pmatrix},
\end{equation}
where the Hessian of~$\mathcal{L}_\ell$ is evaluated at~$(x, y(x), z_I(x), z_E(x))$, the Jacobian matrices $\nabla_x c_I^{\top}$, $\nabla_x c_E^{\top}$, $\nabla_y c_I^{\top}$, and $\nabla_y c_E^{\top}$ are evaluated at~$(x,y(x))$, $C_I$ is a diagonal matrix whose elements are given by~$c_I(x,y(x))$, and~$z_I \circ \nabla_x c_I^{\top}$ is a matrix obtained by multiplying the entries of~$z_I$ by the corresponding rows of~$\nabla_x c_I^{\top}$ (a similar explanation applies to~$z_I \circ \nabla_y c_I^{\top}$). 

Since under Assumption~\ref{ass:LL_assumptions_constr} in Subsection~\ref{subsubsec:LL_constr_case}, the Jacobian~$\nabla_v G^{\top}$ is non-singular at~$(x,v(x))$ (see~\cite{GPMcCormick_2014,JNocedal_SJWright_2006}), we obtain
\begin{equation*}\label{jacobian_pi}
    \nabla v \; = \; \begin{pmatrix} \nabla y, \nabla z_I, \nabla z_E \end{pmatrix} \; = \; - \nabla_x G \, \nabla_v G^{-1}.
\end{equation*}
We can now pull out the columns of this system that correspond to the $\nabla y(x)$ term by introducing an appropriate matrix $L = \begin{pmatrix} \mathbf{I}_m & \mathbf{0}\end{pmatrix}^{\top}$, where~$\mathbf{I}_m$ is an identity matrix of size~$m$ and~$\mathbf{0}$ is a null matrix of size~$m \times \big( \vert I \vert + \vert E \vert \big)$, yielding
\begin{equation}\label{nabla_y}
    \nabla y(x) \; = \; - \nabla_x G \nabla_v G^{-1} L.
\end{equation}
Substituting~\eqref{nabla_y} into~\eqref{adjoint_jacobian}, we obtain the following adjoint gradient for the~LL constrained case:
\begin{equation}\label{constr_direction}
    \nabla f \; = \; \nabla_x f_u - \nabla_x G \nabla_v G^{-1} L \nabla_y f_u,
\end{equation}
where the gradients of~$f_u$ with respect to~$x$ and~$y$ are evaluated at~$(x,y(x))$, while the Jacobians of~$G$ with respect to~$x$ and~$v$ are evaluated at~$(x,v(x))$. The negative adjoint gradient provides the steepest descent direction for~$f$ at~$x$ in the deterministic case. Note that one can deal with the inverse Jacobian in~\eqref{constr_direction} by applying the same approaches used for the~LL unconstrained case, i.e., solving the adjoint equation or using a truncated Neumann series (see Subsection~\ref{sec:inexact_adjoint}).

Similar to the notation used for the~LL unconstrained case in Subsection~\ref{sec:unconstrained}, given $(x_k, \tilde{v}_k)$, with~$\tilde{v}_k = (\tilde{y}_k,{(\tilde{z}_I)}_k,{(\tilde{z}_E)}_k)$, the Jacobian estimates~$\mathcal{G}_x(x_k, \tilde v_k, \varsigma_k^\ell)^{\top}$ and~$\mathcal{G}_v(x_k, \tilde v_k, \varsigma_k^\ell)^{\top}$ approximate~$\nabla_x G(x_k, \tilde v_k)^{\top}$ and~$\nabla_v G(x_k, \tilde v_k)^{\top}$, respectively. Based on Assumption~\ref{ass:cont_diff-2-constr} in Subsection~\ref{subsubsec:LL_constr_case}, $\nabla_{v} G$ and~$\mathcal{G}_v$ are non-singular at all points.
An approximate (negative)~BSG can be computed directly from the adjoint formula~(\ref{constr_direction}), as follows:
\begin{equation}\label{stoch_constr_eq1}
    d(x_k,\tilde{v}_k, \xi_k) \; = \; -\left( g^u_x(x_k,\tilde{y}_k,\vartheta^u_k)  - \mathcal{G}_x(x_k, \tilde v_k, \varsigma_k^\ell) \, \mathcal{G}_v(x_k, \tilde{v}_k, \varsigma_k^\ell)^{-1} L \, g^u_y(x_k,\tilde{y}_k,\vartheta^u_k)\right).
\end{equation}

The data in \eqref{constr_direction} is referred to as~$D(x,v)$ or~$D(x,v(x))$, depending on the point where the gradients, Hessians, and Jacobians are evaluated:
\begin{equation} \label{eq:D(x,v)} 
D(x,v) \; = \;
\left( \nabla_x f_u(x,y), \,L \nabla_y f_u(x,y), \, \nabla_x G(x, v), \,\nabla_v G(x, v) \right).
\end{equation}
The data in calculation~\eqref{stoch_constr_eq1} is now referred to as $D(x_k,\tilde{v}_k,\xi_k)$:
\begin{equation} \label{eq:D(x,v,epsilon)} 
D(x_k,\tilde{v}_k,\xi_k) \; = \;
\left( g^u_x(x_k,\tilde{y}_k,\vartheta^u_k),\, L \, g^u_y(x_k,\tilde{y}_k,\vartheta^u_k), \,\mathcal{G}_x(x_k, \tilde v_k, \varsigma_k^\ell), \,\mathcal{G}_v(x_k, \tilde v_k, \varsigma_k^\ell) \right).
\end{equation}
Again, note that in~\eqref{stoch_constr_eq1}, $d(x_k, \tilde{v}_k, \xi_k)$ can be interpreted as a function of the data~$D(x_k,\tilde{v}_k,\xi_k)$ (see the explanation for~\eqref{eq:dir_notation}).

It is worth mentioning that an alternative approach to obtain a direction colinear with the negative gradient of~$f$ in the~LL constrained case was proposed in~\cite{GSavard_JGauvin_1994}. However, such an approach requires solving an auxiliary linear-quadratic bilevel problem, which is not practical in terms of solving large-scale ML application problems. Also, the approach in~\cite{GSavard_JGauvin_1994} does not yield the exact size of the gradient.  


\subsection{A unified notation for~LL unconstrained and constrained cases}\label{subsec:unified_notation}

Our goal is to propose a general algorithm that applies to both the~LL unconstrained and constrained cases. For each case, one can denote the~BSG directions used in the deterministic setting (i.e.,~\eqref{adjoint} and~\eqref{constr_direction}) and stochastic one (i.e,~\eqref{stochastic-adjoint} and~\eqref{stoch_constr_eq1}) by using the unified notation
\begin{equation}\label{eq:adjoint_general}
    d(D) = -(a - A B^{-1} b),
\end{equation}
where~$D = (a,b,A,B)$. In particular, in the LL unconstrained case, when the adjoint formula~(\ref{adjoint}) or (\ref{stochastic-adjoint}) is used, the data~$D$ is either the deterministic one~(\ref{eq:D(x,y)}) or the stochastic one~(\ref{eq:D(x,y,epsilon)}), respectively. Similarly, in the constrained case, when the adjoint formula~\eqref{constr_direction} or~\eqref{stoch_constr_eq1} is used, the data~$D$ is again either the deterministic one~\eqref{eq:D(x,v)} or the stochastic one~\eqref{eq:D(x,v,epsilon)}, respectively. 
Moreover, we use the following notation to encapsulate the~LL variables in the two cases
\begin{equation}\label{eq:notation}
w \; = \;
    \begin{cases}
        y \text{ when } Y(x) = \mathbb{R}^n,\\
        v \text{ when } Y(x) \neq \mathbb{R}^n,
    \end{cases}
\end{equation}
i.e.,~$w$ is equal to~$y$ in the unconstrained case and~$v$ in the constrained case.

\subsection{The BSG method}

The schema of the~BSG method is presented in Algorithm~\ref{alg:BSG_DS}. An initial point $(x_0,\tilde{w}_0)$ and a sequence of positive scalars~$\{\alpha_k\}$ are required as input. In Step~1, any arbitrary optimization method can be applied to approximately solve the LL problem, regardless of being unconstrained or constrained. 
In Step~2, one computes an approximate (negative)~BSG, which will be denoted by~$d(x_k, \tilde{w}_k, \xi_k)$ and computed through~\eqref{stochastic-adjoint} or~\eqref{stoch_constr_eq1}. We recall from~\eqref{eq:dir_notation} that~$d(x_k, \tilde{w}_k, \xi_k) = d(D(x_k,\tilde{w}_k,\xi_k))$, where~$D(x_k,\tilde{w}_k,\xi_k)$ is either~\eqref{eq:D(x,y,epsilon)} or~\eqref{eq:D(x,v,epsilon)}.
Finally, at Step~3, the vector~$x$ is updated by using a proper step size taken from the sequence of positive scalars. When~$X$ is a closed and convex set different from~$\mathbb{R}^n$, we need to compute the orthogonal projection of~$x_{k} + \alpha_k \, d(x_k, \tilde{w}_k, \xi_k)$ onto~$X$ (note that such a projection can be computed by solving a convex optimization problem).

     \begin{algorithm}[H]
	\caption{Bilevel Stochastic Gradient (BSG) Method}\label{alg:BSG_DS}
	\begin{algorithmic}[1]
		\medskip
		\item[] {\bf Input:} $(x_0,\tilde w_0)$, $\{\alpha_k\}_{k \geq 0} > 0$.
		\medskip
		\item[] {\bf For $k = 0, 1, 2, \ldots$ \bf do}
		\item[] \quad\quad {\bf Step 1.}
		 Obtain an approximation $\tilde w_k$ to the LL optimal solution $w(x_k)$. 
 \nonumber
		\item[] \quad\quad {\bf Step 2.} Compute a (negative) stochastic gradient approximation $d(x_k, \tilde w_k, \xi_k)$.
		\item[] \quad\quad {\bf Step 3.} Compute $x_{k+1} = P_X ( x_{k} + \alpha_k \, d(x_k, \tilde w_k, \xi_k) )$.
		\item[] {\bf End do}
		\par\bigskip\noindent
    	\end{algorithmic}
    \end{algorithm}

We point out that as is usual in the literature related to SG methods, a stopping criterion is not considered due to a lack of reasonable criteria and for the need to study the asymptotic convergence properties.

\subsection{Computing the BSG direction inexactly}\label{sec:inexact_adjoint}

The two approaches proposed in the literature to deal with the inverse matrix~$B^{-1}$ in~\eqref{eq:adjoint_general} consist of either solving the adjoint equation or using a truncated Neumann series. In particular, the first approach requires solving the adjoint equation $B \, \lambda = b$ for the adjoint variables~$\lambda$. The~BSG direction can thus be calculated by $a - A \, \lambda$. The residual error due to the inexact solution of the adjoint equation is denoted by~$\tilde{r}$, i.e., $\tilde{r} = B \, \lambda - b$. Note that the previous expression can be written as $B \, \lambda = b + \tilde{r}$, where the right-hand side can be interpreted as a perturbation of the right-hand side in the adjoint equation. Since $\lambda = B^{-1} (b + \tilde{r})$ in the inexact adjoint solve case, the~BSG direction becomes 
\begin{equation}\label{eq:inexact_adjoint_solve}
-(a - A \, B^{-1} (b + \tilde{r})).
\end{equation}
Given a positive scalar~$q > 0$ and assuming~$\Vert B \Vert < 1$, the second approach is based on the Neumann series as follows: 
\begin{equation*}
    B^{-1} \; = \; \sum_{i=0}^{\infty} (I - B)^i \; = \; \mathscr{B} + \tilde{R},
\end{equation*}
where $\mathscr{B} = \sum_{i=0}^{q} (I - B)^i$ and $\tilde{R} = \sum_{i=q+1}^{\infty} (I - B)^i$. Note that the accuracy of the approximation is an increasing function of~$q$. An approximation to $B^{-1}$ is given by~$\mathscr{B}$, i.e., $B^{-1} \simeq \mathscr{B}$. Therefore, in the inexact Neumann series case, the~BSG direction becomes
\begin{equation}\label{eq:inexact_neumann} -(a - A \, (B^{-1} - \tilde{R}) \, b).
\end{equation}
	\section{Assumptions, sensitivity, and smoothness}\label{sec:assumptions_BSGmethod}

In this section, we introduce the assumptions used in the convergence analysis of the~BSG method, which extends the convergence theory of the~SG method to the bilevel case when the stepsize is assumed to be decaying. Using the notation introduced in~\eqref{eq:adjoint_general}--\eqref{eq:notation}, the convergence theory developed in this section covers both the~LL unconstrained and constrained cases. The solution of the~LL problem is assumed to be inexact. Our theory also applies when the~BSG direction~\eqref{eq:adjoint_general} is computed inexactly by using the approaches in Subsection~\ref{sec:inexact_adjoint}, which comprehensively generalizes all existing approaches in the literature.

Given that we consider the application of the stochastic gradient method (or a similar~SA technique) to solve the~LL problem (see Step~1 of Algorithm~\ref{alg:BSG_DS}), we will denote by~$\xi_k^{\sone}$ the set of random variables for all combined iterations of the~LL solution process at iteration~$k$. Moreover, we will denote by~$\xi_k^{\all}$ the set of all random variables for both the~LL and~UL solution processes at iteration~$k$. Therefore, noticing that~$\xi_k$ denotes the set of random variables used at Step~2 of Algorithm~\ref{alg:BSG_DS}, we denote~$\xi_k^{\all}=(\xi_k^{\sone},\xi_k)$. At each iteration, the iterate~$x_k$ is completely determined by the realizations of the independent random variables~$\xi_k^{\sone}$ and~$\xi_k$.

\subsection{General assumptions}\label{subsec:general_assumptions}

\subsubsection{LL unconstrained case}\label{subsubsec:LL_unc_case}

We will start this subsection with general assumptions for the~LL unconstrained case. Assumption~\ref{ass:smoothness} below imposes the appropriate smoothness for the problem gradients and Hessians, their boundedness, and the boundedness of their stochastic counterparts.

\begin{assumption}[Smoothness and boundedness (LL unconstrained case)]\label{ass:smoothness}
    The gradient~$\nabla f_u$ and the Hessians~$\nabla_{xy}^2 f_\ell$ and~$\nabla^2_{yy} f_\ell$
    are Lipschitz continuous. Moreover, $\nabla f_u$, $\nabla_{xy}^2 f_\ell$, and~$\nabla^2_{yy} f_\ell$ and their stochastic estimates~$g^u_x$, $g^u_y$, $H_{xy}^{\ell}$, and~$H_{yy}^{\ell}$ are bounded at all points.
\end{assumption}
    
Assumption~\ref{ass:cont_diff} below ensures the existence and uniqueness of an~LL optimal solution~$y(x)$ for the original problem as well as for its corresponding stochastic approximation.

\begin{assumption}[Existence and uniqueness of solution (LL unconstrained case)]
\label{ass:cont_diff}
\phantom{.}\linebreak There exists a $y(x)$ such that $\nabla f_{\ell} (x,y(x)) = 0$, and $\nabla_{yy}^2 f_{\ell}$ is positive definite at all points.
In the stochastic case, $H_{yy}^{\ell}$ is positive definite at all points.
\end{assumption}

Assumption~\ref{ass:cont_diff-2} below requires the inverse of the Hessian of $f_\ell$ w.r.t.~$y$ to be uniformly bounded, which is equivalent to saying that~$f_{\ell}(x,\cdot)$ is strongly convex
(a very standard assumption in the literature, see, e.g.,~\cite{SGhadimi_MWang_2018}).
We also need uniform boundedness in the stochastic case. 

\begin{assumption}[Uniform convexity of problem (LL unconstrained case)]
\label{ass:cont_diff-2}
\phantom{}\linebreak
The Hessians~$\nabla^2_{yy} f_{\ell}^{-1}$ and~$(H_{xy}^{\ell})^{-1}$ are uniformly bounded at all points.
\end{assumption}

\subsubsection{LL constrained case}\label{subsubsec:LL_constr_case}

Now, we will introduce assumptions that are specific to the~LL constrained case. We will start with smoothness assumptions for the gradients and Hessians of the problem, their boundedness, and the boundedness of their stochastic counterparts.

\begin{assumption}[Smoothness and boundedness (LL constrained case)]\label{ass:smoothness-constrained}
The gradient $\nabla f_u$, the Jacobians~$\nabla c_I^{\top}$ and~$\nabla c_E^{\top}$, and the Hessians~$\nabla_{xy}^2 f_\ell$, $\nabla^2_{yy} f_\ell$, $\nabla_{yx}^2 c_i$, and~$\nabla_{yy}^2 c_i$, for all~$i \in I \cup E$, are Lipschitz continuous. In addition, $\nabla f_u$, $\nabla c_I^{\top}$, $\nabla c_E^{\top}$, $\nabla_{xy}^2 f_\ell$, $\nabla^2_{yy} f_\ell$, $\nabla_{yx}^2 c_i$, and~$\nabla_{yy}^2 c_i$, for all~$i \in I \cup E$, and their stochastic estimates are uniformly bounded at all points.
Finally, the multipliers $z$ and their stochastic counterparts are uniformly bounded at all points.
\end{assumption}

The boundedness of the multipliers is required so that all Hessians of the Lagrangian are also bounded as well as cross terms of the type $z \circ \nabla c^\top$.

We can now introduce our assumption for the existence and uniqueness of LL solutions. For that purpose, given Lagrange multipliers~$(z_I(x), z_E(x))$ associated with a solution $y(x)$ of the LL KKT conditions, let us denote the cone of critical directions~\cite{JNocedal_SJWright_2006} as follows:
\begin{equation}\label{eq:critical_cone}
Z(x) \; = \; \left\{\begin{matrix}
    \phantom{d^y \ne 0 ~:~  } \hspace{0.0cm}\nabla_y c_i(x,y(x))^\top d^y \le 0, \ \forall i \in I(x) \phantom{\text{ with } (z_I(x))_i > 0}\\
    d^y \ne 0 ~:~ \nabla_y c_i(x,y(x))^\top d^y = 0, \ \forall i \in I(x) \text{ with } (z_I(x))_i > 0\\ 
    \phantom{d^y \ne 0 : } \nabla_y c_i(x,y(x))^\top d^y = 0, \ \forall i \in E \phantom{\text{ with } (z_I(x))_i > 0}
    \end{matrix} \right\},
\end{equation}
where~$I(x)$ is the index set of the active inequality constraints at $(x,y(x))$ defined in Subsection~\ref{sec:constrained}. 

The linear independence constraint qualification~(LICQ) ensures that the \linebreak gradients $\nabla_y c_i(x,y(x))$, for all~$i \in I(x)\cup E$,
are linearly independent. The strict complementarity slackness condition~(SCS) states that for all multipliers~$(z_I(x),z_E(x))$ satisfying the~LL KKT conditions at~$(x,y(x))$, one has~$(z_I(x))_i > 0$ for all $i \in I(x)$. The second-order sufficient condition~(SOSC) states that for all multipliers~$(z_I(x), z_E(x))$ satisfying the~LL KKT conditions at~$(x,y(x))$ and for all~$d^y \in Z(x)$, where~$Z(x)$ is defined by~\eqref{eq:critical_cone}, one has~$(d^y)^\top \nabla^2_{yy} \mathcal{L}_{\ell}(x,v(x)) d^y > 0$.

\begin{assumption}[Existence and uniqueness of solution (LL constrained case)]\label{ass:LL_assumptions_constr}
\phantom{.}\linebreak There exists a $y(x)$ satisfying the LL KKT conditions with associated multipliers~$(z_I(x), z_E(x))$ such that the LICQ, SCS, and SOSC are satisfied.

This guarantees that $y(x)$ is a strict local minimizer for the~LL problem. To ensure that $y(x)$ is the unique global minimizer, one can further assume either that~$\nabla_{yy}^2 \mathcal{L}_{\ell} (x, y, z_I, z_E)$ is positive semi-definite for all~$(y, z_I, z_E)$ or that~$Y(x)$ is convex and~$\nabla_{yy}^2 \mathcal{L}_{\ell} (x, y, z_I, z_E)$ is positive semi-definite for all~$(y, z_I, z_E)$ on the tangent cone to the set~$Y(x)$ at~$y(x)$. 
\end{assumption}

As in the unconstrained case, we also need to impose some form of uniform boundedness away from singularity. In the constrained case, this is achieved through the KKT matrix of the LL problem, which is the corresponding counterpart to the Hessian of the objective function~$f_\ell$. 

\begin{assumption}[``Uniform convexity'' of problem (LL constrained case)]\label{ass:cont_diff-2-constr} The KKT matrices $\nabla_{v} G^{-1}$ and~$\mathcal{G}_v^{-1}$ are uniformly bounded at all points.
\end{assumption}

\subsubsection{Notation for constants}\label{subsubsec:constants}
In this subsection, we introduce constants that will be used to denote bounds on gradients, Hessians, and Jacobians (Remark~\ref{remark} below), and a bound on the second moment of the~BSG direction (Assumption~\ref{ass:unbiasedness} below).

\begin{remark}\label{remark}
As a consequence of Assumptions~\ref{ass:smoothness} and~\ref{ass:smoothness-constrained}, there exist positive constants~$C$ and~$\bar{C}$ such that, for any~$(x,y)$, $(x,v)$, $(x, y, \vartheta)$, $(x, y, \vartheta^{\ell})$, and~$(x, v, \varsigma^\ell)$, we have $\Vert \nabla_y f_u \Vert \le C$, $\Vert \nabla_{xy}^2 f_\ell \Vert \le C$, $\Vert \nabla_{yy}^2 f_\ell \Vert \le C$, $\Vert \nabla_x G \Vert \le C$, $\Vert \nabla_v G \Vert \le C$, $\Vert g^u_y \Vert \le \bar{C}$, $\Vert H_{xy}^{\ell} \Vert \le \bar{C}$, $\Vert H_{yy}^{\ell} \Vert \le \bar{C}$, $\Vert \mathcal{G}_x \Vert \le \bar{C}$, and~$\Vert \mathcal{G}_v \Vert \le \bar{C}$.  Assumptions~\ref{ass:cont_diff-2} and~\ref{ass:cont_diff-2-constr} imply that there exist positive constants~$C_{\ell}$ and~$\bar{C}_{\ell}$ such that, for any~$(x,y)$, $(x, y, \vartheta^{\ell})$, $(x,v)$, and~$(x, v, \varsigma^\ell)$, we have~$\Vert \nabla_{yy}^2 f_{\ell}^{-1} \Vert \le C_{\ell}$, $\Vert (H_{yy}^{\ell})^{-1} \Vert \le \bar{C}_{\ell}$, $\Vert \nabla_v G^{-1} \Vert \le C_{\ell}$, and~$\Vert \mathcal{G}_v^{-1} \Vert \le \bar{C}_{\ell}$.
\end{remark}

In Assumption~\ref{ass:unbiasedness} below, we require the BSG direction~$d(x_k,w(x_k),\xi_k)$ to have a bounded second moment, which is a classical assumption in the~SG literature~\cite{LBottou_FECurtis_JNocedal_2018}. Such an assumption will be used to derive a bound on the second moment of the approximate~BSG direction (see Lemma~\ref{lemma:Gd(r)}).	 
The expected value with respect to the probability distributions of~$\xi_k$ and~$\xi_k^{\sone}$ are denoted by~$\mathbb{E}_{\xi_k}[\cdot]$ and~$\mathbb{E}_{\xi_k^{\sone}}[\cdot]$, respectively. The expected value with respect to the joint distribution of~$\xi_k$ and~$\xi_k^{\sone}$ is denoted by~$\mathbb{E}_{\xi_k^{\all}} = \mathbb{E}_{\xi_k}[\mathbb{E}_{\xi_k^{\sone}}[\cdot]]$.

	\begin{assumption}[Bound on the second moment of the BSG direction]\label{ass:unbiasedness}
	There exists a positive scalar $V_d$ such that the vector $d(x_k,w(x_k),\xi_k)$ satisfies the following condition: 
	\begin{align*}
	    \mathbb{E}_{\xi_k}[\Vert d(x_k,w(x_k),\xi_k) \Vert^2] \; \leq \; V_d.
	    \label{eq5010}
	\end{align*}
	\end{assumption}

\subsection{Sensitivity of the approximate bilevel stochastic gradient direction}\label{sec:sensitivity}

Let us recall that~$f(x) = f_u(x,y(x))$. To bound the second moment of the approximate~BSG direction~$d(x_k,\tilde w_k,\xi_k)$ and the expectation of the error between the negative gradient~$- \nabla f(x_k)$ and~$d(x_k,\tilde w_k,\xi_k)$, we will need to apply sensitivity analysis arguments from nonlinear optimization.
We start by assuming that at each iteration, the calculation process of the BSG direction~\eqref{eq:adjoint_general} is (approximately) Lipschitz continuous with respect to changes in its data. Note that this result is presented as a general assumption that any stochastic algorithm for solving bilevel problems needs to satisfy in order for our convergence theory in Section~\ref{sec:convergence_rate_BSGmethod} to hold. However, such an assumption is not restrictive, and Proposition~\ref{prop:sensitivity_BSG_dir} shows that it can be easily enforced in all practical scenarios that have been proposed for~ML applications requiring bilevel optimization. 
	
	\begin{assumption}[Sensitivity of the BSG direction]\label{ass:lipschitz_continuity}
	Given any pair of data $(D_1)_k$ \linebreak and~$(D_2)_k$, there exists a constant $L_{BSG} > 0$ such that 
    	\begin{equation}\label{Ass:33}
    	    \Vert d((D_1)_k) - d((D_2)_k) \Vert \; \le \; L_{BSG} (\Vert (D_1)_k - (D_2)_k \Vert  +  \Vert (r_1)_k - (r_2)_k \Vert),
    	\end{equation}
     where~$(r_1)_k$ and~$(r_2)_k$ are the residual errors in the inexact computations of~$d((D_1)_k)$ \linebreak and~$d((D_2)_k)$. Moreover, when~$(D_1)_k = D(x_k, w(x_k), \xi_k)$ and~$(D_2)_k = D(x_k, \tilde{w}_k, \xi_k)$, one has
     \begin{equation}\label{Ass:33-additional}
	    \Vert (D_1)_k - (D_2)_k \Vert \; \le \; \bar{L}_{LL} \Vert w(x_k) - \tilde{w}_k \Vert,
	\end{equation} 	
 where~$\bar{L}_{LL} > 0$ is a constant only dependent on the Lipschitz constants of the stochastic gradients, Hessians, and Jacobians of Assumptions~\ref{ass:smoothness} and~\ref{ass:smoothness-constrained}.
    \end{assumption}

Proposition~\ref{prop:sensitivity_BSG_dir} below shows that the inexact ways~\eqref{eq:inexact_adjoint_solve} and~\eqref{eq:inexact_neumann} of calculating adjoint gradients or~BSG directions do ensure that inequality~\eqref{Ass:33} of Assumption~\ref{ass:lipschitz_continuity} is satisfied. Parts of the proof have been published elsewhere~\cite{NCouellan_WWang_2016}. In fact, the arguments used are the known facts that sum is Lipschitz continuous, multiplication is Lipschitz continuous if the factors are bounded, and matrix inversion is Lipschitz continuous if its singular values are bounded away from zero. Moreover, Proposition~\ref{prop:sensitivity_BSG_dir} shows that~inequality~\eqref{Ass:33-additional} can be satisfied when assuming the Lipschitz continuity of the stochastic estimates used in the~BSG directions. We point out that such an assumption is met when the gradients, Hessians, and Jacobians have a finite-sum structure (as commonly found in~ML application problems) and all the terms included in the sums are Lipschitz continuous. Note that this aligns with the typical scenario in~ML, where the same gradient, Hessian, or Jacobian function is evaluated at different data points. In this case, the randomness~$\xi$ consists of drawing a batch, which trivially renders the assumption true. 

\begin{proposition}[Sensitivity of the BSG direction]\label{prop:sensitivity_BSG_dir}
Under Assumptions~\ref{ass:smoothness}--\ref{ass:cont_diff-2-constr}, given any pair of data~$(D_1)_k$ and~$(D_2)_k$, let~$(r_1)_k$ and~$(r_2)_k$ be the residual errors incurred when~$d((D_1)_k)$ and~$d((D_2)_k)$ are computed inexactly by either~\eqref{eq:inexact_adjoint_solve} (in which case~$(r_1)_k$, $(r_2)_k$ are~$(\tilde{r}_1)_k$, $(\tilde{r}_2)_k$) or~\eqref{eq:inexact_neumann} (in which case~$(r_1)_k$, $(r_2)_k$ are~$(\tilde{R}_1)_k$, $(\tilde{R}_2)_k$). Then, there exists a constant~$L_{BSG} > 0$ such that inequality~\eqref{Ass:33} of Assumption~\ref{ass:lipschitz_continuity} is satisfied. When~$(D_1)_k = D(x_k, w(x_k), \xi_k)$ and~$(D_2)_k = D(x_k, \tilde{w}_k, \xi_k)$, assuming the stochastic estimates of the gradients, Hessians, and Jacobians in~$(D_1)_k$ and~$(D_2)_k$ to be Lipschitz continuous in~$w$ for all~$\xi_k$, there exists a positive constant~$\Bar L_{LL}$ such that inequality~\eqref{Ass:33-additional} of Assumption~\ref{ass:lipschitz_continuity} is satisfied.
 \end{proposition}

 \begin{proof}
     See Appendix~\ref{appendix:sensitivity_BSG_dir} for the proof, where we omit the dependence on~$k$ for simplicity. 
 \end{proof}

We now introduce Assumption~\ref{ass:assumption_SGD} on the absolute error of the~LL optimal solution, whose validity in practice is discussed in Subsection~\ref{sec:bound_SGD}.

	\begin{assumption}\label{ass:assumption_SGD}
    	There exists a positive scalar~$C_w$ such that
    	   \begin{equation*} \label{eq:conv_y}
        	    \mathbb{E}_{\xi_k^{\sone}}[\Vert w(x_k) - \tilde w_k\Vert^2] \; \le \; (C_w \, \alpha_k)^2.
    	    \end{equation*}
	\end{assumption}

By applying Jensen's inequality, one also has~$(\mathbb{E}_{\xi_k^{\sone}}[\Vert w(x_k) - \tilde w_k\Vert])^2 \; \le \; \mathbb{E}_{\xi_k^{\sone}}[\Vert w(x_k) - \tilde w_k\Vert^2]$, thus
\begin{equation} \label{eq:conv_y_jensen}
    \mathbb{E}_{\xi_k^{\sone}}[\Vert w(x_k) - \tilde w_k\Vert] \; \le \; C_w \, \alpha_k.
\end{equation}

We also need Assumption~\ref{ass:assumption_sampling} below to hold which essentially amounts to the sampling error in the data. To enforce this assumption in practice, we refer the reader to the discussion reported in Section~\ref{sec:bound_biasedness}.
Recall that~$D(x_k,\tilde{w}_k)$ represents the deterministic data defining the quantity~$d(x_k,\tilde{w}_k)$ (see~\eqref{eq:D(x,y)} and~\eqref{eq:D(x,v)}), and~$D(x_k,\tilde{w}_k,\xi_k)$
the stochastic data of the calculation of~$d(x_k,\tilde{w}_k,\xi_k)$ (see~\eqref{eq:D(x,y,epsilon)} and~\eqref{eq:D(x,v,epsilon)}). 

	\begin{assumption}\label{ass:assumption_sampling}
    	There exists a positive scalar~$C_D$ such that
    	   \[
        	    \expc_{\xi_k^{\all}}[ \Vert D(x_k,\tilde{w}_k) - D(x_k,\tilde{w}_k,\xi_k) \Vert] \; \le \; C_D \, \alpha_k.
    	    \]
	\end{assumption}
	
Finally, we need to bound the residual errors introduced in Assumption~\ref{ass:lipschitz_continuity}. In practice, when the~BSG direction is computed inexactly (see Subsection~\ref{sec:inexact_adjoint}), this assumption can be enforced by either increasing the accuracy in the inexact adjoint equation solve or increasing the value of~$q$ used to truncate the Neumann series.
\begin{assumption}\label{ass:assumption_residual}
	There exists a positive scalar $C_e$ such that, for all realizations of the algorithm,
	   \[
    	   \Vert r_k \Vert  \; \le \; C_e \, \alpha_k,
	    \]
	where~$r_k$ is either~$(r_1)_k$ or~$(r_2)_k$ in Assumption~\ref{ass:lipschitz_continuity}.
\end{assumption}

One is ready to establish the desired bounds (see Lemmas~\ref{lemma:Gd(r)} and~\ref{lemma:Cd(r)} below), for which the constants will depend on the above-introduced constants.

\begin{lemma}\label{lemma:Gd(r)}
Under Assumptions \ref{ass:smoothness}--\ref{ass:assumption_SGD} and~\ref{ass:assumption_residual}, and assuming the stepsize sequence~$\{\alpha_k\}$ bounded from above by a constant~$C_s > 0$, one has
\begin{equation}\label{eq501}
	    \mathbb{E}_{\xi_k^{\all}}[\Vert d(x_k,\tilde{w}_k,\xi_k) \Vert^2] \; \leq \; G_d,
\end{equation}
where~$G_d = 2 \, (U_d + V_d)$, $U_d = L_{BSG}^2 C_s^{2} (L_{LL}^2 C_w^2 + 4 \bar{L}_{LL} C_e C_w + 4 C_e^2)$, and~$\bar{L}_{LL} > 0$ is the constant introduced in Assumption~\ref{ass:lipschitz_continuity}.
\end{lemma}
\begin{proof}
	By considering the data $D_1 = D(x_k, w(x_k), \xi_k)$ and $D_2 = D(x_k, \tilde{w}_k, \xi_k)$ in Assumption~\ref{ass:lipschitz_continuity}, we obtain
	\begin{equation}\label{eq:eq_55555}
	\begin{split}
	    \Vert d(x_k, w(x_k), \xi_k) - d(x_k, \tilde{w}_k, \xi_k)\Vert  \; \le \; &L_{BSG} \Vert D(x_k, w(x_k), \xi_k) - D(x_k, \tilde{w}_k, \xi_k) \Vert \\ \; &+ \; L_{BSG} \Vert r_1 - r_2 \Vert.
	\end{split}
	\end{equation}
From~\eqref{eq:eq_55555}, inequality~\eqref{Ass:33-additional} of Assumption~\ref{ass:lipschitz_continuity}, and Assumption~\ref{ass:assumption_residual}, we obtain
	\begin{equation}
	    \Vert d(x_k,w(x_k),\xi_k) - d(x_k,\tilde{w}_k,\xi_k) \Vert \; \le \; L_{BSG} (\bar{L}_{LL} \Vert w(x_k) - \tilde{w}_k \Vert + 2 C_e \alpha_k). \label{eq:555500}
	\end{equation}
By raising both sides of~\eqref{eq:555500} to the second power, we have  
	\begin{equation*}
	    \Vert d(x_k,w(x_k),\xi_k) - d(x_k,\tilde{w}_k,\xi_k) \Vert^2 \; \le \; L_{BSG}^2 (\bar{L}_{LL}^2 \Vert w(x_k) - \tilde{w}_k \Vert^2 + 4 \bar{L}_{LL} C_e \alpha_k \Vert w(x_k) - \tilde{w}_k \Vert  + 4 C_e^2 \alpha_k^2).
	\end{equation*}
    Therefore, by taking expectations with respect to the distribution of~$\xi_k^{\all}$, considering Assumption~\ref{ass:assumption_SGD}, \eqref{eq:conv_y_jensen}, and the bound on the stepsize sequence, and denoting~$U_d = L_{BSG}^2 C_s^{2} (\bar{L}_{LL}^2 C_w^2 + 4 \bar{L}_{LL} C_e C_w + 4 C_e^2)$, we obtain
	\begin{equation}
	    \mathbb{E}_{\xi_k^{\all}}[\Vert d(x_k,w(x_k),\xi_k) - d(x_k,\tilde{w}_k,\xi_k) \Vert^2] \; \le \; U_d, \label{eq:55553}
	\end{equation}
	where we have used that~$\mathbb{E}_{\xi_k^{\all}}[\Vert w(x_k) - \tilde w_k\Vert^{i}] \; = \; \mathbb{E}_{\xi_k^{\sone}}[\Vert w(x_k) - \tilde w_k\Vert^{i}]$, with~$i \in \{1,2\}$.
	
From~Assumption~\ref{ass:unbiasedness} and~\eqref{eq:55553}, we can obtain the desired bound on the second moment of the approximate~BSG direction by adding and subtracting~$d(x_k,w(x_k),\xi_k)$ and considering that $\mathbb{E}_{\xi_k^{\all}}[\Vert d(x_k,w(x_k),\xi_k) \Vert^{2}] = \mathbb{E}_{\xi_k}[\Vert d(x_k,w(x_k),\xi_k) \Vert^{2}]$. In particular, we have
\begin{equation*}\label{eq501_1}
\begin{split}
	    \mathbb{E}_{\xi_k^{\all}}[\Vert d(x_k,\tilde{w}_k,\xi_k) \Vert^2] \; \leq \;\; &2\, \mathbb{E}_{\xi_k^{\all}}[\Vert d(x_k,\tilde{w}_k,\xi_k) - d(x_k,w(x_k),\xi_k)\Vert^2]\\ + &2\,\mathbb{E}_{\xi_k}[\Vert d(x_k,w(x_k),\xi_k) \Vert^2]\\ \; \leq \;\; &2 \, (U_d + V_d).
	\end{split}
\end{equation*}
\end{proof}

\begin{lemma}\label{lemma:Cd(r)}
Under Assumptions~\ref{ass:smoothness}--\ref{ass:cont_diff-2-constr} and~\ref{ass:lipschitz_continuity}--\ref{ass:assumption_residual},
\begin{equation}
	    \mathbb{E}_{\xi_k^{\all}}[\Vert -\nabla f(x_k) - d(x_k,\tilde{w}_k,\xi_k) \Vert] \; \le \;  C_d \, \alpha_k, \label{eq:5555}
	\end{equation}
	where~$C_d = L_{BSG} (L_{LL} C_w + C_D + 4 C_e)$, and~$L_{LL} > 0$ is a constant only dependent on the Lipschitz constants of the gradients, Hessians, and Jacobians of Assumptions~\ref{ass:smoothness} and~\ref{ass:smoothness-constrained}.
\end{lemma}
\begin{proof}
    By adding and subtracting the term $d(x_k,\tilde{w}_k)$ and using the triangle inequality, we have
	\begin{align}
	    \expc_{\xi_k^{\all}}[\Vert d(x_k,w(x_k)) - d(x_k,\tilde{w}_k,\xi_k) \Vert] \; \le \;\;  &\expc_{\xi_k^{\all}}[\Vert d(x_k,w(x_k)) - d(x_k,\tilde{w}_k) \Vert] \label{eq:3332}\\ \vspace{1ex} \quad \; + &\expc_{\xi_k^{\all}}[\Vert d(x_k,\tilde{w}_k) - d(x_k,\tilde{w}_k,\xi_k)\Vert]. \label{eq:3333}
	\end{align}
	
	Now we derive a bound for the right-hand side in~\eqref{eq:3332}. By considering the data $D_1 = D(x_k,w(x_k))$ and $D_2 = D(x_k,\tilde{w}_k)$ in Assumption~\ref{ass:lipschitz_continuity}, and taking the expectation, we obtain
	\begin{equation}\label{eq:1_1}
	    \expc_{\xi_k^{\all}}[\Vert d(x_k,w(x_k)) - d(x_k,\tilde{w}_k) \Vert] \; \le \; L_{BSG} \expc_{\xi_k^{\all}}[\Vert D(x_k,w(x_k)) - D(x_k,\tilde{w}_k) \Vert] \; + \; 2 L_{BSG} C_e \alpha_k,
	\end{equation}
    where we have applied Assumption~\ref{ass:assumption_residual} on~$\Vert (r_1)_k - (r_2)_k \Vert$.

	Note that the right-hand side of~\eqref{eq:1_1} contains exact~BLP  gradients and Hessians (or Jacobians in the~LL constrained case). Therefore, the Lipschitz continuity of those mappings  (Assumptions~\ref{ass:smoothness} and~\ref{ass:smoothness-constrained}) implies the existence of a constant~$L_{LL} > 0$ such that
	\begin{equation}
	    \Vert D(x_k,w(x_k)) - D(x_k,\tilde{w}_k) \Vert \; \le \; L_{LL} \Vert w(x_k) - \tilde{w}_k \Vert. \label{eq:4444}
	\end{equation}
	Taking expectations with respect to the distribution of~$\xi_k^{\all}$ on both sides of~(\ref{eq:4444}), we can write
	\begin{equation}
	    \mathbb{E}_{\xi_k^{\all}}[\Vert D(x_k,w(x_k)) - D(x_k,\tilde{w}_k) \Vert] \; \le \; L_{LL} \, \mathbb{E}_{\xi_k^{\sone}}[\Vert w(x_k) - \tilde{w}_k\Vert], \label{eq1000}
	\end{equation}
	where we have used that~$\mathbb{E}_{\xi_k^{\all}}[\Vert w(x_k) - \tilde w_k\Vert] \; = \; \mathbb{E}_{\xi_k^{\sone}}[\Vert w(x_k) - \tilde w_k\Vert]$. Therefore, from \eqref{eq:1_1}, \eqref{eq1000}, and~\eqref{eq:conv_y_jensen}, we obtain
	\begin{equation}
	    \expc_{\xi_k^{\all}}[\Vert d(x_k,w(x_k)) - d(x_k,\tilde{w}_k) \Vert] \; \le \; L_{BSG} \, ( L_{LL} \, C_w \, \; + \; 2 C_e) \alpha_k. \label{eq1001}
	\end{equation}
	
Now we derive a bound for \eqref{eq:3333}. By considering the data $D_1 = D(x_k,\tilde{w}_k)$ and $D_2 = D(x_k,\tilde{w}_k,\xi_k)$ in Assumption~\ref{ass:lipschitz_continuity} and applying Assumption~\ref{ass:assumption_residual}, we have
	\begin{equation} \label{eq:sensitivity_analysis_inexactness}
	  \Vert d(x_k, \tilde{w}_k) - d(x_k, \tilde{w}_k, \xi_k)\Vert \; \le \; L_{BSG} \Vert D(x_k, \tilde{w}_k) - D(x_k, \tilde{w}_k, \xi_k) \Vert \; + \; 2 L_{BSG} C_e \alpha_k.
	\end{equation}
Taking expectations with respect to the distribution of~$\xi_k^{\all}$ on both sides of (\ref{eq:sensitivity_analysis_inexactness}), we obtain
	\begin{equation}\label{eq:1_2}
	    \expc_{\xi_k^{\all}}[\Vert d(x_k,\tilde{w}_k) - d(x_k,\tilde{w}_k,\xi_k)\Vert] \; \le \; L_{BSG} \expc_{\xi_k^{\all}}[\Vert D(x_k,\tilde{w}_k) - D(x_k,\tilde{w}_k,\xi_k)\Vert] \; + \; 2 L_{BSG} C_e \alpha_k.
	\end{equation}	
Hence, from \eqref{eq:1_2} and Assumption~\ref{ass:assumption_sampling}, we obtain
	\begin{equation}
	    \expc_{\xi_k^{\all}}[\Vert d(x_k,\tilde{w}_k) - d(x_k,\tilde{w}_k,\xi_k) \Vert]
	    \; \le \; L_{BSG} \, (C_D + 2 C_e) \alpha_k. \label{eq1001_1}
	\end{equation}
The proof can be concluded from \eqref{eq:3332}--\eqref{eq:3333}, \eqref{eq1001}, and \eqref{eq1001_1}.
\end{proof}

\subsection{Smoothness of the true objective function}
    
    Our convergence theory requires smoothness of the true function~$f$, which is given in Proposition~\ref{prop:smoothness_f_adj_grad} below. In both the~LL unconstrained and constrained cases, the Lipschitz continuity of~~$\nabla f$ can be inferred from the Lipschitz continuity of~$y(x)$, $w(x)$, and the gradients, Hessians, and Jacobians involved (along with the boundedness away from singularity of Hessian or~KKT matrices). The proof is given in Appendix~\ref{appendix:smoothness_f_adj_grad}, and again part of the proof has been reported in~\cite{NCouellan_WWang_2016} and relies on elementary arguments.

    \begin{proposition}[Smoothness of $f$]\label{prop:smoothness_f_adj_grad}
    Under Assumptions~\ref{ass:smoothness}--\ref{ass:cont_diff-2-constr},
    there exists a \linebreak constant~$L_{\nabla f} > 0$ such that the gradient~$\nabla f$ is Lipschitz continuous in $x$, i.e.,
    \begin{equation}\label{eq:eq_lipschitz_cont_adj_grad}
        \Vert \nabla f (x_1) - \nabla f (x_2)\Vert \; \le \; L_{\nabla {f}} \Vert x_1 - x_2\Vert \quad \text{ for all } (x_1,x_2) \in \mathbb{R}^n \times \mathbb{R}^n.
    \end{equation}
    \end{proposition}
    
An important and well-known consequence that follows from~\eqref{eq:eq_lipschitz_cont_adj_grad} is
	\begin{equation}\label{ass_37_result}
	    f(x) \; \leq \; f(\Bar{x}) + \nabla f(\Bar{x})^\top (x - \Bar{x}) + \frac{1}{2}L_{\nabla f} \|x - \Bar{x}\|^2 \; \text{ for all } \; (x,\Bar{x})\in\mathbb{R}^n\times\mathbb{R}^n.
	\end{equation}

	\section{Convergence rate of the BSG method}\label{sec:convergence_rate_BSGmethod}

In this section, we extend the convergence theory of the~SG method to the bilevel case when the stepsize is assumed to be decaying. The~BLP objective function~$f$ is assumed to be nonconvex, strongly convex (leading to a $1/k$ sublinear convergence rate), or simply convex ($1/\sqrt{k}$ rate).

Using the notation introduced in~\eqref{eq:adjoint_general}--\eqref{eq:notation}, the convergence theory developed in this section covers both the~LL unconstrained and constrained cases. The~BSG method under consideration takes an inexact solution of the~LL problem, for which the stochasticity is rigorously included in the analysis for the first time. Moreover, such a theory also applies when the~BSG direction~\eqref{eq:adjoint_general} is computed inexactly regardless of the approach used (see Subsection~\ref{sec:inexact_adjoint}), thus leading to an analysis that is considerably more general than the ones proposed in the literature~\cite{SGhadimi_MWang_2018,FPedregosa_2016}. 

\subsection{Rate in the nonconvex case}
	
In this section, the true objective function~$f$ is assumed to be possibly nonconvex. We now present two lemmas that will allow us to prove the convergence result. Such lemmas and the resulting theorem are based on the theory provided in~\cite{LBottou_FECurtis_JNocedal_2018}, where the main differences lie in the use of the projection operator to handle the upper-level constraints~$x \in X$ and, importantly, in how inexactly~$d(x_k,\tilde{w}_k,\xi_k)$ approximates~$-\nabla f(x_k)$ (see Lemmas~\ref{lemma:Gd(r)} and~\ref{lemma:Cd(r)}). The first lemma is just a Taylor bound derived as a result of~\eqref{eq:eq_lipschitz_cont_adj_grad}.

\begin{lemma}\label{lemma:41}
Under Assumptions~\ref{ass:smoothness}--\ref{ass:cont_diff-2-constr}, the iterates of Algorithm~\ref{alg:BSG_DS} satisfy the following inequality for all~$k\in\mathbb{N}$
\begin{equation}\label{eq:3325}
\begin{split}
    \mathbb{E}_{\xi_k^{\all}}\left[ f(x_{k+1}) \right] - f(x_k) & \leq \alpha_k \left(P_X \nabla f(x_k)\right)^\top \mathbb{E}_{\xi_k^{\all}}\left[d(x_k, \tilde{w}_k, \xi_k)\right] \\ 
    &\quad + \frac{1}{2}\alpha_k^2 L_{\nabla f} \mathbb{E}_{\xi_k^{\all}}\left[\|d(x_k, \tilde{w}_k, \xi_k)\|^2\right].
\end{split}
\end{equation}
\end{lemma}
\begin{proof} From equation \eqref{ass_37_result}, the iterates generated by Algorithm~\ref{alg:BSG_DS} satisfy
\[f(x_{k+1}) - f(x_k) \leq \nabla f(x_k)^\top (x_{k+1} - x_k) + \frac{1}{2}L_{\nabla f} \|x_{k+1} - x_k\|^2.\]
Recalling that Algorithm 1 uses the update $x_{k+1} = P_X\left( x_k + \alpha_k d(x_k, \tilde{w}_k, \xi_k) \right)$, and since $x_k$ is in the feasible region $X$, we know that $x_k = P_X(x_k)$, which yields
\begin{alignat}{2}
    & f(x_{k+1}) - f(x_k) && \; \leq \; 
    \alpha_k \nabla f(x_k)^\top P_Xd(x_k, \tilde{w}_k, \xi_k) + \frac{1}{2}\alpha_k^2 L_{\nabla f} \|P_X d(x_k, \tilde{w}_k, \xi_k)\|^2. \nonumber
\end{alignat}
Since~$P_X$ is an orthogonal projection, we know~$P_X = P_X^\top = P_X^2$ and~$\|P_X\| \leq 1$. Using this fact and 
taking expectations with respect to the distribution of $\xi_k^{\all}$, we obtain~\eqref{eq:3325}.
\end{proof}

The following lemma further extends the result of Lemma~\ref{lemma:41} by using the inexactness of~$d(x_k, \tilde{w}_k, \xi_k)$ and the bound on its second-order moment.

\begin{lemma}\label{lemma:42}
        Under Assumptions~\ref{ass:smoothness}--\ref{ass:assumption_residual}, the iterates generated by Algorithm~\ref{alg:BSG_DS} satisfy the following inequality for all $k\in\mathbb{N}$:
	\begin{equation}\label{lemma_2}
	    \mathbb{E}_{\xi_k^{\all}}\left[ f(x_{k+1}) \right] - f(x_k) \leq - \alpha_k \|P_X\nabla f(x_k)\|^2 + \alpha_k^2 C_{\nabla f} C_d + \frac{1}{2}\alpha_k^2 L_{\nabla f} G_d,
	\end{equation}
 where $C_{\nabla f} > 0$ is a bound on the norm of~$\nabla f$.
\end{lemma}
    
\begin{proof}
From inequality~\eqref{eq:3325} and Lemma~\ref{lemma:Gd(r)}, we have
\[\mathbb{E}_{\xi_k^{\all}}\left[ f(x_{k+1}) \right] - f(x_k) \leq \alpha_k \left(P_X \nabla f(x_k)\right)^\top \mathbb{E}_{\xi_k^{\all}}\left[d(x_k, \tilde{w}_k, \xi_k)\right] + \frac{1}{2}\alpha_k^2 L_{\nabla f} G_d.\]
Adding and subtracting $\alpha_k \left(P_X \nabla f(x_k)\right)^\top \mathbb{E}_{\xi_k^{\all}}\left[ \nabla f(x_k)\right]$ to the right-hand side and simplifying, we obtain
\begin{alignat}{2}
    \mathbb{E}_{\xi_k^{\all}}\left[ f(x_{k+1}) \right] - f(x_k) & \; \leq \; 
    \alpha_k \left(P_X \nabla f(x_k)\right)^\top \mathbb{E}_{\xi_k^{\all}}\left[ d(x_k, \tilde{w}_k, \xi_k) + \nabla f(x_k) \right] \nonumber\\
    &\quad  \; - \; \alpha_k \left(P_X \nabla f(x_k)\right)^\top \nabla f(x_k) + \frac{1}{2}\alpha_k^2 L_{\nabla f} G_d. \nonumber
\end{alignat}
Applying the properties of orthogonal projections along with the Cauchy-Schwarz and Jensen's inequalities, we obtain
\begin{alignat}{2}
    \mathbb{E}_{\xi_k^{\all}}\left[ f(x_{k+1}) \right] - f(x_k) & \; \leq \; \alpha_k \|P_X \nabla f(x_k)\| \; \mathbb{E}_{\xi_k^{\all}}\left[ \|d(x_k, \tilde{w}_k, \xi_k) + \nabla f(x_k)\|\right] \nonumber\\
    &\quad \; - \; \alpha_k \|P_X\nabla f(x_k)\|^2 + \frac{1}{2}\alpha_k^2 L_{\nabla f} G_d. \nonumber
\end{alignat}
Assumptions~\ref{ass:smoothness}, \ref{ass:cont_diff-2}, \ref{ass:smoothness-constrained}, and~\ref{ass:cont_diff-2-constr} imply that there exists a constant~$C_{\nabla f} > 0$ such that the gradients generated by the sequence of iterates $\{x_k\}_{k \ge 0}$ are bounded, i.e.,~$\|\nabla f(x_k)\| \; \leq \; C_{\nabla f}$.
Finally, from the boundedness of~$\nabla f$ and inequality~\eqref{eq:5555} along with the properties of orthogonal projections, we obtain the desired result.
\end{proof}

We will now introduce the final assumptions that are needed for the convergence result of the nonconvex case. The first of these states that, for Algorithm~\ref{alg:BSG_DS} to converge, the sequence of function values must be bounded below by some minimum value.
    
    \begin{assumption}\label{ass:bounded f}
	The sequence~$\{f(x_k)\}_{k \ge 0}$ is bounded below by~$f_{\text{inf}}$.
	\end{assumption}

Lastly, we require the stepsize to be of decaying type.

    \begin{assumption}\label{ass:diminishing stepsize}
	The sequence of decaying stepsizes~$\{\alpha_k\}_{k \ge 0}$ satisfies
	\[
	\sum_{k=0}^\infty \alpha_k = \infty \; \text{ and } \; \sum_{k = 0}^\infty \alpha_k^2 < \infty.
	\]
	\end{assumption}

We can now establish the convergence result for the nonconvex case. We use $\mathbb{E}[\cdot]$ to refer to the {\it total expectation} of $f$, namely, the expected value with respect to the joint distribution of all the random vectors $\xi_k^{\all}$. 

\begin{theorem}\label{theorem:nonconvex}
    Under Assumptions~\ref{ass:smoothness}--\ref{ass:assumption_residual} and~\ref{ass:bounded f}, suppose that Algorithm~\ref{alg:BSG_DS} is run with a decaying stepsize sequence that satisfies Assumption~\ref{ass:diminishing stepsize}. Then, with $A_K := \sum_{k = 0}^K \alpha_k$,
    \begin{equation}\label{nonconv_eq1}
        \lim_{K\rightarrow\infty} \mathbb{E}\left[ \sum_{k = 0}^K \alpha_k \|P_X\nabla f(x_k)\|^2 \right] < \infty,
    \end{equation}
    and therefore
    \begin{equation}\label{nonconv_eq2}
        \lim_{K\rightarrow\infty} \mathbb{E}\left[ \frac{1}{A_K} \sum_{k = 0}^K \alpha_k \|P_X\nabla f(x_k)\|^2 \right] = 0.
    \end{equation}
\end{theorem}
    
\begin{proof}
The proof follows~\cite[Theorem 4.10]{LBottou_FECurtis_JNocedal_2018} closely. Taking the total expectation of \eqref{lemma_2}, we have
\[\mathbb{E}\left[ f(x_{k+1}) \right] - \mathbb{E}\left[f(x_k)\right] \; \leq \; - \alpha_k \mathbb{E}\left[\|P_X\nabla f(x_k)\|^2\right] + \alpha_k^2 C_{\nabla f} C_d + \frac{1}{2}\alpha_k^2 L_{\nabla f} G_d.\]
Summing both sides of this inequality for $k\in\{0,1,...,K\}$ and by Assumption \ref{ass:bounded f}, we have
\begin{alignat}{2}
    & f_{\text{inf}} - \mathbb{E}\left[f(x_{0})\right] && \; \leq \; \mathbb{E}\left[ f(x_{K+1}) \right] - \mathbb{E}\left[f(x_{0})\right] \nonumber\\
    & && \; \leq \; - \sum_{k = 0}^K \alpha_k \mathbb{E}\left[\|P_X\nabla f(x_k)\|^2\right] + C_{\nabla f} C_d \sum_{k = 0}^K \alpha_k^2 + \frac{1}{2} L_{\nabla f} G_d \sum_{k = 0}^K \alpha_k^2. \nonumber
\end{alignat}
Rearranging, we obtain
\[\sum_{k = 0}^K \alpha_k \mathbb{E}\left[\|P_X\nabla f(x_k)\|^2\right] \; \leq \; \mathbb{E}\left[f(x_{0})\right] - f_{\text{inf}} + C_{\nabla f} C_d \sum_{k = 0}^K \alpha_k^2 + \frac{1}{2} L_{\nabla f} G_d \sum_{k = 0}^K \alpha_k^2.\]
Assumption~\ref{ass:diminishing stepsize} implies that the right-hand side of this inequality converges to a finite limit when~$K$ increases, which proves~\eqref{nonconv_eq1}. To obtain~\eqref{nonconv_eq2}, we can divide by~$A_K$ as follows:
\[\frac{1}{A_K}\sum_{k = 0}^K \alpha_k \mathbb{E}\left[\|P_X\nabla f(x_k)\|^2\right] \; \leq \; \frac{\mathbb{E}\left[f(x_{0})\right] - f_{\text{inf}}}{A_K} + \frac{C_{\nabla f} C_d}{A_K} \sum_{k = 0}^K \alpha_k^2 + \frac{L_{\nabla f} G_d}{2A_K}  \sum_{k = 0}^K \alpha_k^2.\]
Taking the limit as $K\rightarrow\infty$, and noting Assumption~\ref{ass:diminishing stepsize}, we obtain the desired result.
\end{proof}

\subsection{Rate in the strongly convex case}\label{subsec:strongly_convex}

In this subsection, we present the convergence rate of the BSG method when $f$ is assumed to be strongly convex. In practice, such a case occurs when the~UL objective function~$f_u$ is strongly convex and~$y(x)$ is an affine function in~$x$. Hence, imposing strong convexity of~$f$ is a strong assumption in the sense of assuming in practice that the~LL problem is a~QP problem. Still, we cover this case for completeness of our convergence theory. 

Along with~\eqref{eq:eq_lipschitz_cont_adj_grad}, we also need the iterates to lie in a bounded set, which could be ensured by the boundedness of~$X$ in the BLP formulation.
    
    \begin{assumption}[Boundedness of the iterates]\label{ass:boundedness}
	The sequence of iterates $\{x_k\}_{k \ge 0}$ yielded by Algorithm~\ref{alg:BSG_DS} is contained in a bounded set.
	\end{assumption}
	
	Assumption~\ref{ass:boundedness} implies that there exists a positive constant~$\Theta$ such that, for any $(k_1,k_2)$, we have
	\begin{equation*}\label{eq:bounded_constant}
	  \Vert x_{k_1} - x_{k_2} \Vert \; \le \; \Theta \; < \; \infty.
	\end{equation*}
      
      Finally, we assume that the true function~$f$ is strongly convex. 
    
    \begin{assumption}[Strong convexity of~$f$]\label{ass:strong_convexity}
		The function $f$ is strongly convex, namely, there exists a constant $c > 0$ such that
		\begin{equation}\label{eq:strong_convexity}
	    f(\bar x) \; \ge \; f(x) + \nabla f(x)^\top (\bar x - x) + \frac{c}{2}\Vert \bar x - x \Vert^2 \text{ for all } (\bar x, x) \in \mathbb{R}^n \times \mathbb{R}^n.
	\end{equation}	
	\end{assumption}
    
    A well-known equivalent condition to \eqref{eq:strong_convexity} (see, e.g.,~\cite{YNesterov_2018}) is given by
    \begin{equation}\label{eq:strong_convexity_2}
	    (\nabla f(x) - \nabla f(\bar x))^\top (x - \bar x) \ge c\Vert x - \bar  x \Vert^2 \text{ for all } (x, \bar x) \in \mathbb{R}^n \times \mathbb{R}^n.
	\end{equation}	
    Let $x_*$ be the unique minimizer of~$f$ on $X$, which implies that $\nabla f(x_*)^\top (x - x_*) \ge 0 \text{ for all } x \in X$. Therefore, if in \eqref{eq:strong_convexity_2} we choose $x = x_k$ and $\bar x = x_*$, we obtain
    \begin{equation}\label{eq:strong_convexity_3}
	    \nabla f(x_k)^\top (x_k - x_*) \ge c\Vert x_k - x_* \Vert^2.
	\end{equation}
	The next theorem proves that under the assumption of strong convexity and decaying stepsize ($\sum_{k = 0}^{\infty} \alpha_k = \infty$ and $\sum_{k = 0}^{\infty} \alpha_k^2 < \infty$), the sequence of points yielded by Algorithm~\ref{alg:BSG_DS} generates a sequence of~$f$ values that decays sublinearly at the rate of~$1/k$. The proof of this theorem is given in Appendix~\ref{appendix:conv_2}.
	
    \begin{theorem}\label{theorem:conv_2}
    Let Assumptions~\ref{ass:smoothness}--\ref{ass:assumption_residual} and \ref{ass:boundedness}--\ref{ass:strong_convexity} hold and~$x_*$ be the unique minimizer of~$f$ on $X$. Consider the schema given by Algorithm~\ref{alg:BSG_DS} and assume a decaying step size sequence of the form $\alpha_k = \gamma/k$, where $\gamma \ge 1/(2c)$ is a positive constant. The sequence of iterates yielded by Algorithm~\ref{alg:BSG_DS} satisfies
    \begin{alignat*}{3}
        &\mathbb{E}[\Vert x_k - x_*\Vert^2] &&\; \le \; \frac{\max\{2 \, \gamma^2 M (2 c \gamma - 1)^{-1},\| x_0 - x_*\|^2\}}{k},\\
        &\mathbb{E}[f(x_k)] - f(x_*) &&\; \le \; \frac{(L_{\nabla {f}}/2) \max\{2 \, \gamma^2 M (2 c \gamma - 1)^{-1},\| x_0 - x_*\|^2\}}{k},
    \end{alignat*}
    where~$M = G_d + 2 C_d \Theta.$
    \end{theorem}
    
\subsection{Rate in the convex case}
    
In this subsection, we state the convergence rate of the BSG method assuming that~$f$ is convex and attains a minimizer~$x_*$. The same comment about the lack of practicality applies to the convex case, i.e., that~$y(x)$ would need to be affine, which considerably restricts the choice of~LL problems.

\begin{assumption}[Convexity of~$f$]\label{ass:convexity}
	Given $(\bar y,y) \in \mathbb{R}^m \times \mathbb{R}^m$, the (continuously differentiable) function $f$ is convex in $x$, namely,
	\begin{equation}\label{eq:convexity}
    f(\bar x) \; \ge \; f(x) + \nabla f(x)^\top (\bar x - x) \text{ for all } (\bar x, x) \in \mathbb{R}^n \times \mathbb{R}^n.
\end{equation}
Moreover, $f$ attains a minimizer.
\end{assumption}
The next theorem states that the~BSG method exhibits a sublinear convergence rate of~$1/\sqrt k$, which implies that the convergence is slower than in the strongly convex case~(Theorem~\ref{theorem:conv_2}). The proof of this theorem is given in Appendix~\ref{appendix:conv_convex}.

 \begin{theorem}\label{theorem:conv_convex}
    Let Assumptions~\ref{ass:smoothness}--\ref{ass:assumption_residual}, \ref{ass:boundedness}, and~\ref{ass:convexity} hold. Consider the schema given by Algorithm~\ref{alg:BSG_DS} and assume a decaying step size of the form $\alpha_k = \bar{\alpha}/\sqrt{k}$, with $\bar{\alpha} > 0$. Given a minimizer~$x_*$ of~$f$, the sequence of iterates yielded by Algorithm~\ref{alg:BSG_DS} satisfies
    \[ 
    \min_{s = 0, \ldots, k} \mathbb{E}[f(x_s)] - f(x_*) \; \le \; \frac{\frac{\Theta^2}{2 \bar \alpha} + \bar \alpha (G_dM + 2 C_dM \Theta)}{\sqrt k}.
    \]
\end{theorem}

    \subsection{Imposing a bound on the distance from the LL optimal solution}\label{sec:bound_SGD}

    In this subsection, we want to discuss a way to enforce Assumption~\ref{ass:assumption_SGD} when using the stochastic gradient~(SG) method to solve the LL~problem at~$x_k$. We focus on the~LL unconstrained case. 
Given an initial point~$\tilde{y}_k^0$ and a sequence of stepsizes $\{ \beta_i \}$, such a SG method can be described as 
\begin{equation} \label{eq:SG}
\tilde{y}_k^{i+1} \; = \; \tilde{y}_k^{i} - \beta_i g^\ell_y (x_k, \tilde{y}_k^{i}, \xi^{\sone}_{k,i}), \quad i=0,\ldots,i_k.
\end{equation} 
 We start by introducing the sampling assumptions that are standard in the literature related to the~SG method. First, suppose that the stochastic gradient ~$g^\ell_y(x_k,\tilde{y}_k,\xi^{\sone}_{k,i})$ is unbiased, i.e.,~$\mathbb{E}_{\xi^{\sone}_{k,i}}[g_y^\ell(x_k,\tilde{y}_k,\xi^{\sone}_{k,i})] = \nabla_y f_\ell(x_k,\tilde{y}_k)$, and there exists a positive constant~$Q > 0$ such that~$\mathbb{E}_{\xi^{\sone}_{k,i}}[\|g^\ell_y(x_k,\tilde{y}_k,\xi^{\sone}_{k,i})\|^2] \le Q^2$. Also, suppose that~$f_\ell$ is strongly convex in the~$y$ variables with constant~$\mu$ with a unique minimizer~$y(x_k)$.

    Recalling that the convergence rate of the~SG method~\eqref{eq:SG} with decaying stepsize is~${\cal O}(1/\sqrt{i})$, and by choosing~$i_k$ equal to~$k^2$, one guarantees the existence of a positive constant~$C_y$ such that Assumption~\ref{ass:assumption_SGD} holds.
    In fact, by choosing a decaying step size sequence~$\{\beta_i\}$ given by~$\beta_i = \gamma/i$, where~$\gamma \ge 1/(2\mu)$ is a positive constant, and under the classical assumptions stated in the previous paragraph, from \cite[Equation~(2.9)]{ANemirovski_etal_2009} it follows that the choice~$i_k = k^2$ implies (with~$\tilde{y}_k = \tilde{y}_k^{i_k+1}$)
    \[
    \mathbb{E}_{\xi^{\sone}_{k}}[\Vert \tilde{y}_k - y(x_k)\Vert^2] \; \le \; \frac{\max\{\gamma^2 Q (2 \mu \gamma - 1)^{-1},\|  \tilde{y}_k^0 - y(x_k) \|^2\}}{k^2}.
    \]
    Such a result also holds in the~LL constrained case when the scheme~\eqref{eq:SG} incorporates a projection onto~$Y(x_k)$, as long as~$Y(x_k)$ is a closed convex set. We refer the reader to the discussion in~\cite[Appendix~B]{PKhanduri_etal_2023}, which applies to bilevel problems with an~LL strongly convex objective function and~LL linear constraints. It is also possible to obtain an inequality like~\eqref{ass:assumption_SGD} for the~LL constrained case when using a primal-dual stochastic method~\cite{YXu_2020}.

    \subsection{Imposing a bound on dynamic sampling}\label{sec:bound_biasedness}

In this subsection, we want to mention a dynamic sampling strategy to enforce the inequality in Assumption~\ref{ass:assumption_sampling} in both the~LL unconstrained and constrained cases.
For the sake of simplicity, we will omit the subscript~$k$ in this subsection.
Such a dynamic sampling strategy allows reducing the level of noise by increasing the size of the batch. 
Recalling the definition of~$D(x,w)$ given in~\eqref{eq:D(x,y)} and~\eqref{eq:D(x,v)}, the definition of~$D(x, w,\xi)$ given in~\eqref{eq:D(x,y,epsilon)} and~\eqref{eq:D(x,v,epsilon)}, and the unified notation introduced in Subsection~\ref{subsec:unified_notation}, let us denote
\begin{alignat*}{2}
&D(x,w) &&= (a(x,y), b(x,y), A(x,w), B(x,w)),\\ 
&D(x,w,\xi) &&= (a(x,y,\vartheta^u), b(x,y,\vartheta^u), A(x,w,\gamma), B(x,w,\gamma)),
\end{alignat*}
where~$\gamma$ is either $\vartheta^{\ell}$ (in the~LL unconstrained case) or $\varsigma^{\ell}$ (in the~LL constrained case).
Let us assume that the stochastic estimates~$a(x,y,\vartheta^u)$, $b(x,y,\vartheta^u)$, $A(x,w,\gamma)$, and~$B(x,w,\gamma)$ are normally distributed with means~$a(x,y)$, $b(x,y)$, $A(x,w)$, and $B(x,w)$, respectively, and variances~$\sigma^2_a$, $\sigma^2_b$, $\sigma^2_A$, and $\sigma^2_B$, respectively. Such an assumption implies that the stochastic estimates in~$D(x, w,\xi)$ are unbiased estimates of the corresponding true gradients, Hessians, and Jacobians in~$D(x,w)$.

To increase the accuracy of the stochastic estimates in~$D(x,w,\xi)$, we can choose larger batch sizes, which we denote by~$n_{a}$, $n_{b}$, $n_{A}$, and $n_{B}$. Let~$\bar a(x,y,\vartheta^u) = (1/n_{a}) \sum_{r=1}^{n_{a}} a(x,y,(\vartheta^u)_r)$ be the mini-batch stochastic estimate for~$a(x,y)$, where~$\{(\vartheta^u)_r\}_{r=1}^{n_{a}}$ are values sampled from the distribution of~$\vartheta^u$. It is known that (for details, see, for instance, \cite[Section~5.3]{SLiu_LNVicente_2019})
\begin{equation*}
    \mathbb{E}_{\xi_k^{\all}}[\Vert a(x,y) - \bar a(x,y,\vartheta^u) \Vert] \; \le \; \frac{\sigma_a \sqrt{n}}{\sqrt{n_{a}}}.\label{eq:iiiii}
\end{equation*}
One can repeat similar arguments for~$b$, $A$, and~$B$, and their corresponding mini-batch stochastic estimates~$\bar{b}$, $\bar{A}$, and $\bar{B}$, respectively. Let us denote
\[\bar{D}(x,w,\xi) = (\bar{a}(x,y,\vartheta^u), \bar{b}(x,y,\vartheta^u), \bar{A}(x,w,\gamma), \bar{B}(x,w,\gamma)).\]
From the equivalence of norms, there exists a positive constant~$\hat{C}$ such that    
\begin{alignat*}{2}
        \Vert D(x,w) - \bar D(x,w,\xi) \Vert \; &\le \; \hat{C} (\Vert a(x,y) - \bar{a}(x,y,\vartheta^u) \Vert + \Vert b(x,y) - \bar{b}(x,y,\vartheta^u) \Vert) \\
        & \; + \hat{C} (\Vert A(x,w) - \bar{A}(x,w,\gamma) \Vert + \Vert B(x,w) - \bar{B}(x,w,\gamma) \Vert).
\end{alignat*}
Taking expectations with respect to~$\xi_k^{\all}$ on both sides, one obtains
 \begin{alignat}{2}
 \mathbb{E}_{\xi_k^{\all}}[\Vert D(x,w) - \bar D(x,w,\xi) \Vert] \; &\le \; \hat{C} (\frac{\sigma_a}{\sqrt{n_{a}}} + \frac{\sigma_b}{\sqrt{n_{b}}} + \frac{\sigma_A}{\sqrt{n_{A}}} + \frac{\sigma_B}{\sqrt{n_{B}}})\sqrt{n}.\label{eq:bound_sampling_gaussian}
\end{alignat}
To guarantee that Assumption~\ref{ass:assumption_sampling} holds, we need to choose mini-batch sizes~$n_{a}$, $n_{b}$, $n_{A}$, and $n_{B}$ and sample standard deviations~$\sigma_a$, $\sigma_b$, $\sigma_A$, and $\sigma_B$ such that the right-hand side in~\eqref{eq:bound_sampling_gaussian} is less than or equal to~$C_D \, \alpha_k$. Therefore, when $\alpha_k$ decreases, the dynamic sampling strategy would increase the mini-batch sizes.
	
\section{Numerical experiments}\label{sec:numerical_experiments}

All code was written in Python and the experimental results were obtained on a desktop computer (32GB of RAM, Intel(R) Core(TM) i9-9900K processor running at 3.60GHz).\footnote{All the code for our implementation is available at \url{https://github.com/GdKent/BSG_Methods_Con_Unc}.}  
We averaged all the results over 10 trials by using different random seeds. 

\subsection{Our practical BSG methods}\label{sec:stochasticBSG1}

A major difficulty in the adjoint formulas~\eqref{adjoint} and~\eqref{constr_direction} is the use of second-order derivatives of~$f_\ell$ and~$\mathcal{L}_\ell$, respectively, and the need to solve the adjoint equation or use a truncated Neumann series, which prevents its application to large-scale ML application problems. We propose two approaches to get around this problem, which lead to two different practical versions of the~BSG method, referred to as~BSG-N-FD and~BSG-1. In the numerical experiments considered for both~LL unconstrained and constrained~BLPs, we are mainly interested in testing these two practical implementations (Algorithms~\ref{alg:BSG-N-FD} and~\ref{alg:BSG1} in Subsections~\ref{subsubsec:BSG-N-FD} and~\ref{subsubsec:BSG-1} below, respectively) as opposed to the standard BSG method (Algorithm~\ref{alg:BSG_DS}). We will consider two types of test problems: synthetic quadratic bilevel problems and continual learning problems. Since these test problems do not have upper-level constraints (i.e.,~$X = \mathbb{R}^n$), we have omitted the use of the orthogonal projection operator~$P_X$ in all of the algorithms.

In the numerical experiments for the synthetic problems, we will also test the~BSG method with stochastic Hessians, where the (negative)~BSG direction~$d(x_k, \tilde{y}_k, \xi_k)$ is calculated from~\eqref{stochastic-adjoint} or~\eqref{stoch_constr_eq1}. This version is referred to as BSG-H. In the~LL unconstrained case, BSG-H applies the linear conjugate gradient~(CG) method~\cite{JNocedal_SJWright_2006} to solve the adjoint system $H^{\ell}_{yy}(x_k,\tilde{y}_k,\vartheta^\ell_k) \lambda = g^u_y(x_k,\tilde{y}_k,\vartheta^u_k)$ until non-positive curvature is detected. In the~LL constrained case, BSG-H solves the adjoint system~$\mathcal{G}_v(x_k, \tilde{v}_k, \varsigma_k^\ell) \lambda =  L \, g^u_y(x_k,\tilde{y}_k,\vartheta^u_k)$ by applying the~GMRES method~\cite{YSaad_MHSchultz_1986}. Note that we only include~BSG-H in the experiments for the sake of comparison. For very large problems, one must use~BSG-N-FD or~BSG-1. 

In the~BSG-N-FD, BSG-1, and~BSG-H versions of the BSG method for the~LL unconstrained case, we will apply the SG method~(\ref{eq:SG}) to the~LL problem for a certain budget~$i_k$ of iterations, obtaining an approximation~$\tilde{y}_k$ to the LL optimal solution~$y(x_k)$. To obtain an approximation~$\tilde{w}_k$ to~$w(x_k)$, given~$x_k$, in the~LL constrained case, BSG-N-FD, BSG-1, and~BSG-H will first determine an approximation~$\tilde{y}_k$ to~$y(x_k)$ by minimizing the following exact penalty function over~$y$
\begin{equation}\label{eq:pen_func}
\Phi(x_k, y; \mu) = f_{\ell}(x_k, y) + \frac{1}{\mu} \sum_{i \in I} \max\{0, c_i(x_k,y)\} + \frac{1}{\mu} \sum_{i \in E} \left| c_i(x_k,y)\right|,
\end{equation}
where~$\mu$ is a penalty parameter and the functions~$c_i$, with~$i \in I \cup E$, are the~LL constraints defined in Subsection~\ref{sec:constrained}. We recall that for sufficiently small and positive values
of~$\mu$, the minimization of such an unconstrained problem will yield the optimal solution~$y(x_k)$ of the constrained~LL problem~\cite{JNocedal_SJWright_2006}. To minimize~\eqref{eq:pen_func}, which is a nonsmooth function, the stochastic subgradient method will be applied. Then, given~$x_k$ and~$\tilde{y}_k$, to determine approximations~$({\tilde{z}_I},{\tilde{z}_E})$ to the optimal multipliers~$
(z_I(x_k), z_E(x_k))$, the linear~CG method will be applied to solve the~KKT system~$G(x_k, (\tilde{y}_k, z_I, z_E)) = 0$ for the variables~$z_I$ and~$z_E$, where~$G$ is the vector function introduced in Subsection~\ref{sec:constrained}.

    For the solution of the LL~problem, we consider an inexact scheme, denoted as {\sf inc.~acc.}, which consists of obtaining~$\tilde{y}_k$ by taking multiple steps of the stochastic gradient method applied to~$f_{\ell}$ ($i_k \ge 1$, $\forall k$, in~\eqref{eq:SG}, for the~LL unconstrained case) or stochastic subgradient method applied to~\eqref{eq:pen_func} (for the~LL constrained case). In particular, the number of steps of the stochastic gradient/subgradient method increases by~1 every time the difference of the UL~objective function between two consecutive iterations is less than a given threshold, thus leading to an increasing accuracy strategy. In such an inexact scheme, $\tilde{y}_k$ is determined by using the approximation~$\tilde{y}_{k-1}$ obtained at the previous iteration as a starting point. In the~ML community, iterative schemes with multiple~LL steps, like the {\sf inc.~acc.} strategy above, are referred to as double-loop schemes~\cite{KJi_JYang_YLiang_2020,TChen_YSun_WYin_2021_closingGap}.

\subsubsection{BSG-N-FD}\label{subsubsec:BSG-N-FD}
Our first proposed method, BSG-N-FD, solves the adjoint system by using an iterative method where each Hessian vector product is approximated with a finite-difference~(FD) scheme. In particular, in the~LL unconstrained case, the adjoint equation~$\nabla_{yy}^2 f_\ell \, \lambda = \nabla_y f_u$ is solved for the adjoint variables~$\lambda$ by using the linear~CG method, with~$\nabla_{yy}^2 f_\ell \, \lambda$ being approximated as follows:
\begin{equation}\label{eq:BSGNFD_FD1}
\nabla^2_{y y} f_\ell (x_k, y_k) \lambda \; \approx \; \frac{\nabla_{y} f_\ell(x_k, y^+_k) - \nabla_{y} f_\ell(x_k,y^-_k)}{2 \varepsilon},
\end{equation}
where $y^\pm_k \;=\; y_k \pm \varepsilon \lambda, \mbox{ with } \varepsilon > 0$.
Then, the adjoint gradient is calculated from
\begin{equation}\label{eq:BSGNFD_FD2}
\nabla f \; \approx \; \nabla_x f_u - \nabla_{xy}^2 f_\ell \, \lambda,
\end{equation}
where~$\nabla_{xy}^2 f_\ell \, \lambda$ is approximated using an additional~FD scheme. In practice, we use a scaling parameter value of~$\varepsilon = 0.1$.

In the~LL constrained case, we can use~FD schemes similar to the ones in the~LL unconstrained case to approximate the Jacobian vector products when solving the adjoint equation and computing the adjoint gradient. In particular, the adjoint equation~$\nabla_v G \lambda = L \nabla_y f_u$ is solved for the adjoint variables~$\lambda = (\lambda_y, \lambda_I, \lambda_E)^{\top}$ by using the~GMRES method, with $\nabla_v G \lambda$ being approximated as follows:
\begin{equation}\label{eq:BSGNFD_FD1constr}
\nabla_v G \lambda \; \approx \; \begin{pmatrix} \overline{\nabla_{yy}^2 \mathcal{L}_\ell \lambda_y} + (z_I^{\top} \circ \nabla_y c_I) \lambda_I + \nabla_y c_E \lambda_E \\
\nabla_y c_I^{\top} \lambda_y + C_{I} \lambda_I \\
\nabla_y c_E^{\top} \lambda_y \end{pmatrix},
\end{equation}
where 
\begin{equation*}
\overline{\nabla^2_{y y} \mathcal{L}_\ell \lambda_y} \; = \; \frac{\nabla_{y} \mathcal{L}_\ell(x_k, y^+_k, (z_I)_k, (z_E)_k) - \nabla_{y} \mathcal{L}_\ell(x_k, y^-_k, (z_I)_k, (z_E)_k)}{2 \varepsilon}
\end{equation*}

and $y^\pm_k = y_k \pm \varepsilon \lambda_y\mbox{, with } \varepsilon > 0$. Then, the adjoint gradient is calculated from
\begin{equation}\label{eq:BSGNFD_FD2constr}
\nabla f \; \approx \; \nabla_x f_u - \nabla_x G \lambda,
\end{equation}
where~$\nabla_x G \lambda$ is approximated using an additional~FD scheme. In practice, what worked better for us was again~$\varepsilon = 0.1$.


We denote the algorithm corresponding to this approach as~BSG-N-FD, where the~``N'' stands for the Newton-type system given by the adjoint equation, and the ``FD'' for the finite-difference approximations we employ. Its schema is described in Algorithm~\ref{alg:BSG-N-FD}. In the practical~BSG-N-FD method, we will use the stochastic data defined by~\eqref{eq:D(x,y,epsilon)} and~\eqref{eq:D(x,v,epsilon)}.

     \begin{algorithm}[H]
	\caption{BSG-N-FD Method}\label{alg:BSG-N-FD}
	\begin{algorithmic}[1]
		\medskip
		\item[] {\bf Input:} $(x_0,\tilde{w}_0)$, $\{\alpha_k\}_{k \geq 0} > 0$.
		\medskip
		\item[] {\bf For $k = 0, 1, 2, \ldots$ \bf do}
		\item[] \quad\quad {\bf Step 1.}
		 Obtain an approximation $\tilde{w}_k$ to the LL optimal solution $w(x_k)$. 
 \nonumber
        \item[] \quad\quad {\bf Step 2.} Compute~$d(x_k, \tilde{w}_k, \xi_k)$ by drawing the stochastic gradients and/or Jacobians from~\eqref{eq:BSGNFD_FD1}--\eqref{eq:BSGNFD_FD2} for the~LL unconstrained case or~\eqref{eq:BSGNFD_FD1constr}--\eqref{eq:BSGNFD_FD2constr} for the~LL constrained case.
		\item[] \quad\quad {\bf Step 3.} Compute $x_{k+1} = x_{k} + \alpha_k \, d(x_k, \tilde{w}_k, \xi_k)$.
		\item[] {\bf End do}
		\par\bigskip\noindent
    	\end{algorithmic}
    \end{algorithm}

\subsubsection{BSG-1}\label{subsubsec:BSG-1}
Our second proposed method, BSG-1, approximates the second-order derivatives in the adjoint formulas~\eqref{adjoint} and~\eqref{constr_direction} with outer products of the corresponding gradients, i.e.,
\begin{alignat}{6}
    & \nabla_{xy}^2 f_\ell \; &&\simeq \;  \nabla_x f_\ell \nabla_y f_\ell^\top  \; &\text{ and } \; &\nabla_{yy}^2 f_\ell \; &&\simeq \; \nabla_y f_\ell \nabla_y f_\ell^\top, \label{unconstrained_rank1_approx}\\
    & \nabla_{yx}^2 \mathcal{L}_\ell \; &&\simeq \; \nabla_y \mathcal{L}_\ell \nabla_x \mathcal{L}_\ell^\top \; &\text{ and } \; & \nabla_{yy}^2 \mathcal{L}_\ell \; &&\simeq \; \nabla_{y} \mathcal{L}_\ell \nabla_y \mathcal{L}_\ell^\top.\label{constrained_rank1_approx}
\end{alignat}
As mentioned in Subsection~\ref{sec:cont}, such approximations are inspired by Gauss-Newton methods for nonlinear least-squares problems, where the Hessian matrix of the objective function $\sum_{i=1}^{p} (r_i-a_i)^2$ (in which each $r_i$ is a scalar function and $a_i$ a scalar) is approximated by $\sum_{i=1}^{p} \nabla r_i \nabla r_i^\top$, and also from the fact that the empirical risk of misclassification in ML is often a sum of non-negative terms matching a function to a scalar which can then be considered in a least-squares fashion~\cite{ABotev_2017,MGargiani_etal_2020}. In the numerical experiments, the rank-1 approximations are observed to perform well when the LL function $f_{\ell}$ has a Gauss-Newton structure, such as the binary cross-entropy loss function used for the continual learning instances (see Subsection~\ref{sec:continual_learning}).

In the~LL unconstrained case, the resulting approximate adjoint equation $(\nabla_y f_\ell \nabla_y f_\ell^\top)\, \lambda = \nabla_y f_u$ is most likely infeasible, and we suggest solving it in the least-squares sense. One solution is $\lambda = \nabla_y f_u/(\nabla_y f_\ell^\top \nabla_y f_\ell)$.
Plugging this and $\nabla_{xy}^2 f_\ell \simeq \nabla_x f_\ell \nabla_y f_\ell^\top$ in the adjoint formula~(\ref{adjoint}) gives rise to our practical~BSG-1 calculation 
\begin{equation} \label{approx-BSG}
\nabla f \; \simeq \; \nabla_x f_u - \frac{\nabla_y f_\ell^\top \nabla_y f_u}{\nabla_y f_\ell^\top \nabla_y f_\ell} \nabla_x f_\ell.
\end{equation} 
This approximate BSG allows us to use the adjoint formula without computing Hessians or even Hessian-vector products, which is prohibitively expensive for the large bilevel problems arising in ML applications.

In the~LL constrained case, we can use the outer products~\eqref{constrained_rank1_approx} to approximate the second-order derivatives of~$\mathcal{L}_\ell$ in the Jacobian matrices~$\nabla_x G^{\top}$ and~$\nabla_v G^{\top}$, introduced in~\eqref{eq:jacobians}, and obtain corresponding approximate Jacobians~$\tilde{G}_x^{\top}$ and~$\tilde{G}_v^{\top}$, respectively.
The resulting approximate adjoint equation is given by~$\tilde{G}_v \tilde{\lambda} = L \nabla_y f_u$, where~$L$ is the matrix used in~\eqref{constr_direction}, and can be solved by using an iterative method for non-symmetric linear systems. Plugging a solution~$\tilde{\lambda}$ into~$\nabla_x f_u - \tilde{G}_x \, \tilde{\lambda}$, we obtain the practical~BSG-1 calculation
\begin{equation}\label{constr_direction_approx}
    \nabla f \; \simeq \; \nabla_x f_u - \tilde{G}_x \, \tilde{\lambda}, \quad \text{ where } \quad \tilde{G}_v \tilde{\lambda} =  L \nabla_y f_u.
\end{equation}

Both of these rank-1 approaches for the~LL unconstrained and constrained cases will be referred to as~BSG-1, the~``1'' standing for first-order rank-1 approximations of the Hessian and Jacobian matrices. In the practical~BSG-1 method, we will use the stochastic data defined by~\eqref{eq:D(x,y,epsilon)} and~\eqref{eq:D(x,v,epsilon)}. 

     \begin{algorithm}[H]
	\caption{BSG-1 Method}\label{alg:BSG1}
	\begin{algorithmic}[1]
		\medskip
		\item[] {\bf Input:} $(x_0,\tilde{w}_0)$, $\{\alpha_k\}_{k \geq 0} > 0$.
		\medskip
		\item[] {\bf For $k = 0, 1, 2, \ldots$ \bf do}
		\item[] \quad\quad {\bf Step 1.}
		 Obtain an approximation $\tilde{w}_k$ to the LL optimal solution $w(x_k)$. 
 \nonumber
        \item[] \quad\quad {\bf Step 2.} Compute~$d(x_k, \tilde{w}_k, \xi_k)$ by drawing the stochastic gradients and/or Jacobians from~\eqref{approx-BSG} for the~LL unconstrained case or~\eqref{constr_direction_approx} for the~LL constrained case.
		\item[] \quad\quad {\bf Step 3.} Compute $x_{k+1} = x_{k} + \alpha_k \, d(x_k, \tilde{w}_k, \xi_k)$.
		\item[] {\bf End do}
		\par\bigskip\noindent
    	\end{algorithmic}
    \end{algorithm}
    
\subsection{DARTS}\label{sec:stochasticBLP_DARTS}

DARTS was proposed in~\cite{HLiu_KSimonyan_YYang_2019} for the solution of stochastic BLPs arising from NAS, and was briefly introduced in Section~\ref{sec:BSD}.
Only the LL unconstrained case ($Y(x) = \mathbb{R}^m$) has been considered.
To avoid the computation of the second-order derivatives in~(\ref{eq:DARTS}), 
DARTS approximates the matrix-vector product $\nabla^2_{xy} f_\ell (x_k, y_k) \nabla_{y} f_u (x_k, \tilde y_k)$ by a finite-difference scheme~\cite{HLiu_KSimonyan_YYang_2019}:
	\[
	\nabla^2_{x y} f_\ell (x_k, y_k) \nabla_{y} f_u (x_k,\tilde y_k) \; \approx \; \frac{\nabla_{x} f_\ell(x_k, y^+_k) - \nabla_{x} f_\ell(x_k,y^-_k)}{2 \varepsilon},
	\]
	where 
	\begin{equation} \label{eq:ypm}
	y^\pm_k \;=\; y_k \pm \varepsilon \nabla_{y} f_u(x_k, \tilde y_k) \quad \mbox{with} \quad \varepsilon \;=\; 0.01/\| \nabla_{y} f_u (x_k,\tilde y_k) \|.
	\end{equation}
	Algorithm~\ref{alg:DARTS} reports the schema of DARTS for the stochastic setting. In Step~1, a single step of SG (with fixed stepsize $\eta$) is applied to the LL problem to obtain an approximation $\tilde{y}_k$ to the LL optimal solution. Then, in Step~2, the UL variables are updated by moving along the ``approximated'' descent direction using a stepsize~$\alpha_k$.
	
	 \begin{algorithm}[H]
	\caption{Differentiable Architecture Search (DARTS)}\label{alg:DARTS}
	\begin{algorithmic}[1]
		\medskip
		\item[] {\bf Input:} $(x_0,y_0) \in \mathbb{R}^n \times \mathbb{R}^m$, $\{\alpha_k\}_{k \geq 0} > 0$, $\eta > 0$.
		\medskip
		\item[] {\bf For $k = 0, 1, 2, \ldots$ \bf do}
		\item[] \quad\quad {\bf Step 1.} Compute $\tilde y_k = y_k - \eta \, g_y^\ell(x_k, y_k, \vartheta_k^\ell)$.
		\item[] \quad\quad {\bf Step 2.} Compute $x_{k+1} = x_k - \alpha_k \left( g_x^u(x_k, \tilde y_k, \vartheta_k^u) - \frac{\eta}{2 \varepsilon} (g_x^\ell(x_k, y_k^+, \vartheta_k^\ell) - g_x^\ell(x_k, y_k^-, \vartheta_k^\ell) ) \right)$, with $y_k^\pm$ and $\varepsilon$ as in~(\ref{eq:ypm}), with~$g_y^u(x_k, \tilde y_k, \vartheta_k^u)$ instead of~$\nabla_{y} f_u(x_k, \tilde y_k)$, and set $y_{k+1}=\tilde{y}_{k}$.
		\item[] {\bf End do}
		\par\bigskip\noindent
    	\end{algorithmic}
    \end{algorithm}

\subsection{Numerical results for synthetic quadratic bilevel problems}
	
    We first report results for a ``synthetic'' bilevel problem, where both levels are defined by quadratic objective functions. Given~$h_1 \in \mathbb{R}^n$, $h_2 \in \mathbb{R}^m$, symmetric positive definite matrices~$H_2 \in \mathbb{R}^{n \times n}$ and~$H_3 \in \mathbb{R}^{m \times m}$, and matrices~$H_1 \in \mathbb{R}^{n \times m}$ and~$H_4 \in \mathbb{R}^{m \times n}$, we consider the following problem 
    \begin{equation}\label{prob:synthetic_prob}
        \begin{gathered}
            \min_{x \in \mathbb{R}^n} ~  f_u(x, y) \; = \; h_1^{\top} x + h_2^{\top} y + \frac{1}{2} x^{\top} H_1 y + \frac{1}{2} x^{\top} H_2 x \\
            \text{s.t.} \ \ \ \ y \in \argmin_{y \in Y(x)} ~ f_{\ell}(x, y) \; = \; \frac{1}{2} y^{\top} H_3 y - y^{\top} H_4 x,
        \end{gathered}
    \end{equation}
    where the set~$Y(x)$ used for the numerical experiments will be specified in Subsections~\ref{sec:unc_results}--\ref{sec:quad_constr_results} below. In particular, Subsection~\ref{sec:unc_results} focuses on the~LL unconstrained case and Subsections~\ref{sec:lin_con_results}--\ref{sec:quad_constr_results} address the~LL constrained case.

    For all of the algorithms, we used the best~UL and~LL fixed stepsizes (i.e.,~$\alpha^{u}$ and~$\alpha^{\ell}$, respectively) found by performing a grid search over the set~$\{10^{-s_u} ~|~ s_u \in \{\underline{s}_u, \ldots, \bar{s}_u\}\}$ for the~UL stepsize and~$\{10^{-s_{\ell}} ~|~ s_{\ell} \in \{\underline{s}_{\ell}, \ldots, \bar{s}_{\ell}\}\}$ for the~LL stepsize. We used independent bounds~$\underline{s}_{u}, \bar{s}_{u}, \underline{s}_{\ell}, \bar{s}_{\ell}$ for each algorithm with the goal of selecting stepsize values that capture their best performances on the different types of problems. Using a single arbitrary grid to test all of the algorithms would have led to stepsizes biased toward one set of algorithms over the others.    
    The domain of possible values for the bounds~$\underline{s}_{u}, \bar{s}_{u}, \underline{s}_{\ell}, \bar{s}_{\ell}$ was restricted to the set~$\{1,\ldots,8\}$. We chose the values of such bounds to include two to three consecutive values in the grid searches for both~$s_u$ and~$s_{\ell}$.

    With the exception of~DARTS in the~LL unconstrained case, we utilized the~LL {\sf inc.~acc. strategy (introduced in Subsection~\ref{sec:stochasticBSG1}) for the rest of the algorithms, with an~$f_u$ difference threshold for increasing the number of~LL iterations equal to~$10^{-1}$ (and a maximum limit of~30~LL iterations). In all the figures included in this subsection, we plot the true function~$f$ of the BLP ($f(x_k)=f_u(x_k,y(x_k))$). The number of~UL} iterations and running time were both used as metrics for the comparison of the algorithms.
    In the~LL constrained case, in accordance with the procedure described in Subsection~\ref{sec:stochasticBSG1}, we obtain approximate Lagrange multipliers at each iteration by solving the corresponding~KKT system with the linear conjugate gradient method (with a maximum number of iterations equal to~3 and tolerance equal to~$10^{-4}$). The system in~\eqref{constr_direction_approx} is solved by using the~GMRES method (with a maximum number of~$3$ iterations and tolerance equal to~$10^{-4}$).

\subsubsection{Results for the LL unconstrained case}\label{sec:unc_results}
	In the numerical experiments for the LL unconstrained version ($Y(x) = \mathbb{R}^m$) of problem~\eqref{prob:synthetic_prob}, we considered a dimension of~300 at both the upper and lower levels (i.e.,~$n = m = 300$), with~$H_2$ and~$H_3$ randomly generated and~$H_1$ and~$H_4$ set equal to the identity matrix. Note that problem~\eqref{prob:synthetic_prob} is deterministic. To investigate the numerical performance of the stochastic methods considered in the experiments, we computed stochastic gradient and Hessian estimates by adding Gaussian noise with mean~0 to each corresponding deterministic gradient (i.e.,~$\nabla_x f_u$, $\nabla_y f_u$, $\nabla_x f_{\ell}$, $\nabla_y f_{\ell}$) and Hessian (i.e.,~$\nabla^2_{xy} f_{\ell}$, $\nabla^2_{yy} f_{\ell}$). In the stochastic case, the standard deviations for the stochastic estimates of the gradients and Hessians were set to~$5$ and~$0.05$, respectively. 
 
    In the experiments for this subsection, we compared~BSG-N-FD, BSG-1, and~BSG-H \linebreak against~DARTS and~StocBiO, introduced in Subsection~\ref{sec:BSD}. Regarding the grid searches for the stepsizes~$\alpha_u$ and~$\alpha_{\ell}$, in the deterministic case, we used~$\underline{s}_u = \underline{s}_{\ell} = 2$ and~$\bar{s}_u = \bar{s}_{\ell} = 4$ for~BSG-N-FD, BSG-H, BSG-1, and~StocBiO, and~$\underline{s}_u = 3$, $\bar{s}_u = 5$, $\underline{s}_{\ell} = 2$, and~$\bar{s}_{\ell} = 4$ for~DARTS. In the stochastic case, we used~$\underline{s}_u = \underline{s}_{\ell} = 2$ and~$\bar{s}_u = \bar{s}_{\ell} = 4$ for~BSG-N-FD, BSG-1, and~StocBiO, and~$\underline{s}_u = 3$, $\bar{s}_u = 5$, $\underline{s}_{\ell} = 2$, and~$\bar{s}_{\ell} = 4$ for~BSG-H and~DARTS. For StocBiO, we set the constant~$C_0$ introduced in Subsection~\ref{sec:BSD} to~0.05 and the parameter~$q$ introduced in Subsection~\ref{sec:inexact_adjoint} to~2, which led to the best results.

    \begin{figure}
    \centering
          \includegraphics[scale=0.44]{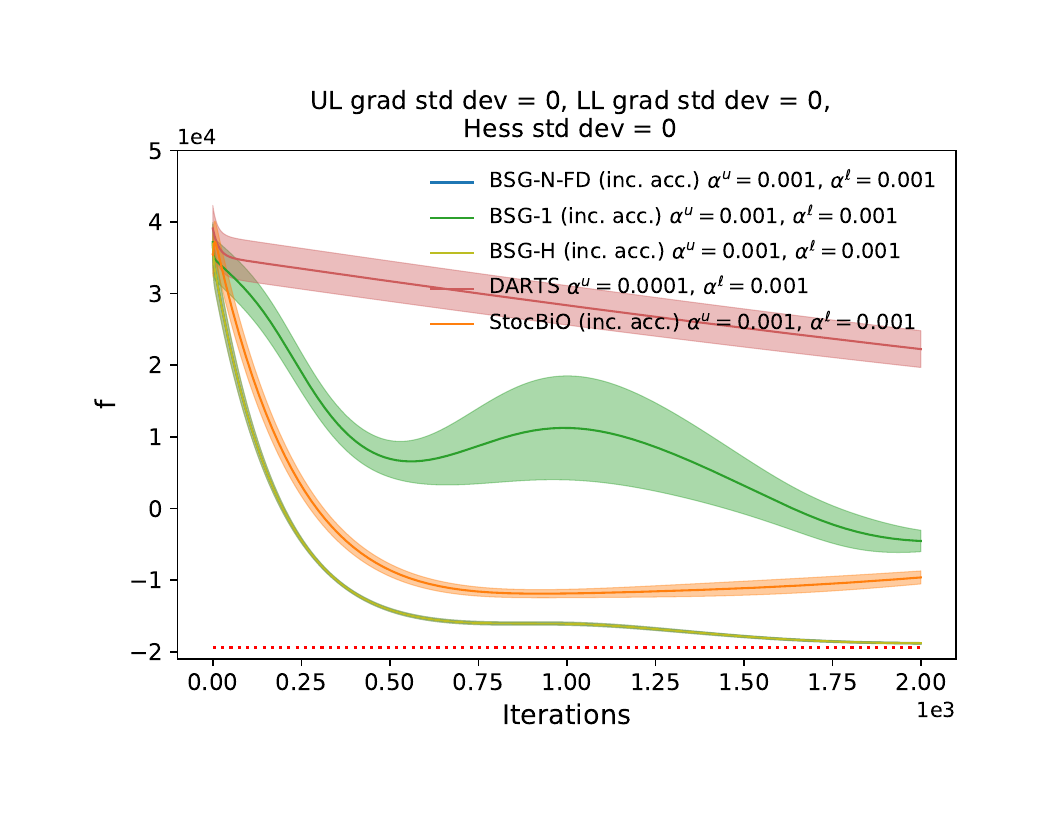}
          \includegraphics[scale=0.44]{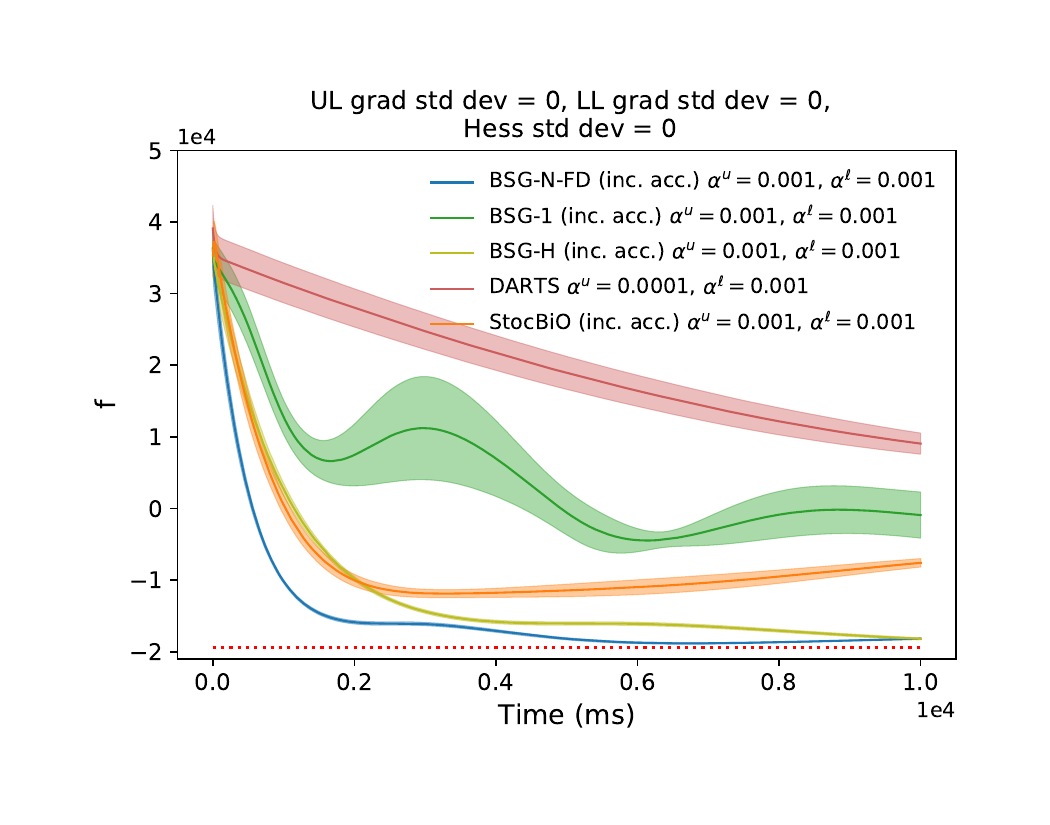}
          \includegraphics[scale=0.44]{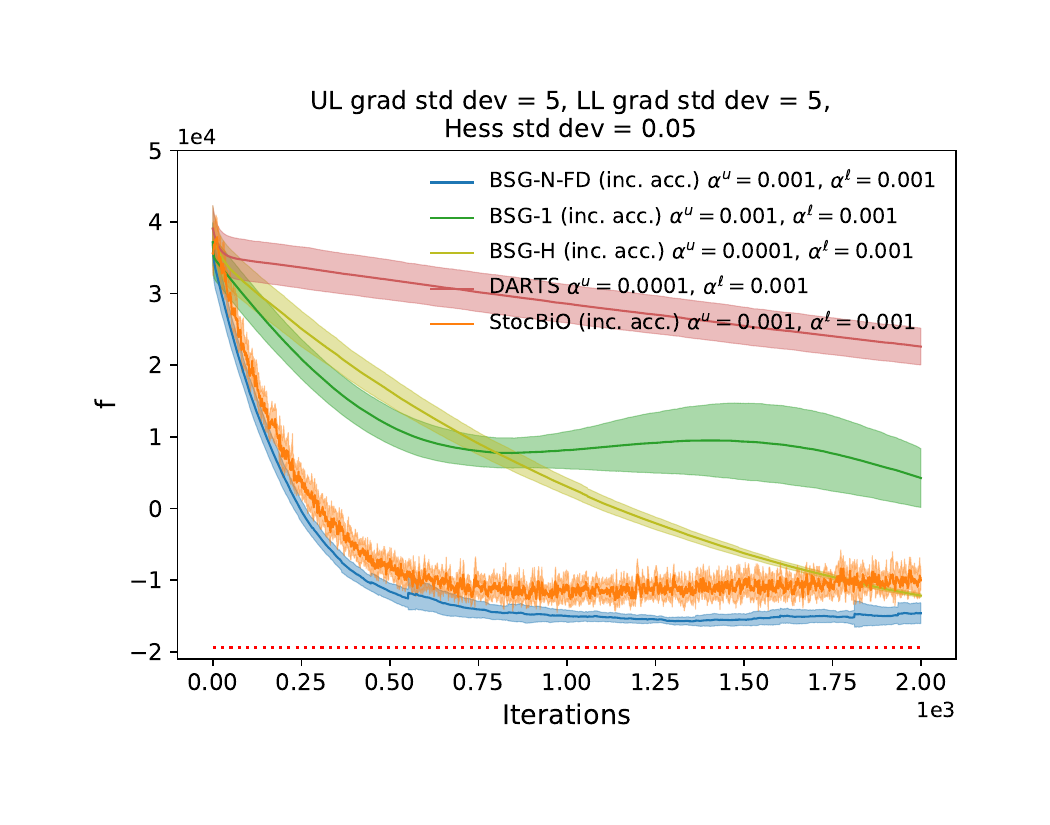}
          \includegraphics[scale=0.44]{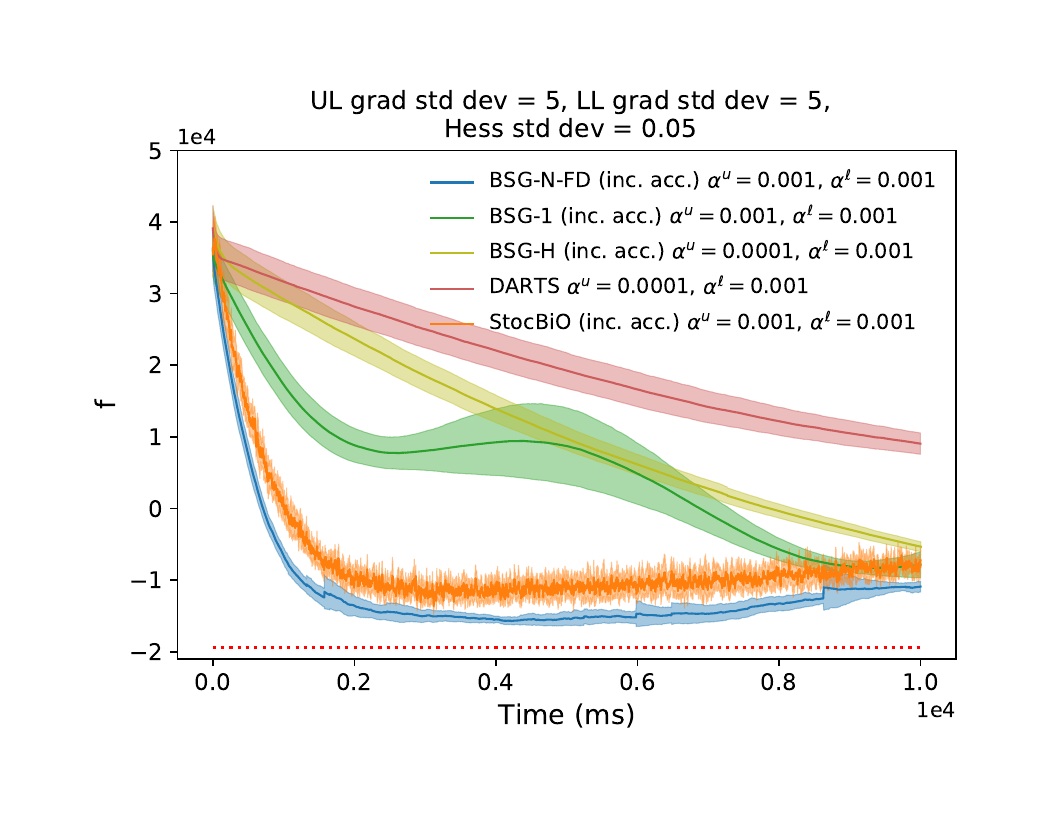}
        \caption{Numerical results of the~BSG-N-FD, BSG-H, BSG-1, DARTS, and~StocBiO algorithms on problem~\eqref{prob:synthetic_prob} for the~LL unconstrained case in terms of both iterations and time (in milliseconds).}\label{fig:synth_unc_plots}
\end{figure}

    Starting with the deterministic results displayed in the top two plots of Figure~\ref{fig:synth_unc_plots}, we can see that both BSG-N-FD and BSG-H (these methods completely overlap in the iterations plot) clearly outperform all the other methods in terms of both iterations and time, with StocBiO performing slightly worse. Note that compared to~BSG-N-FD, BSG-H is less efficient in terms of time due to the high computational cost of computing Hessian matrices. As somewhat expected, due to the lack of Gauss-Newton structure in problem~\eqref{prob:synthetic_prob}, BSG-1 performs worse than the other versions of the~BSG method, but it is still able to yield a decrease in the true function and outperform DARTS, which has the worst performance out of the five methods. It bears mentioning that only BSG-N-FD and BSG-H are able to achieve the optimal value of the true function~$f$ (represented by the red dotted horizontal lines in all the plots of Figure~\ref{fig:synth_unc_plots}) within the iteration and time limit used.

    We will now focus on the results for the stochastic case displayed in the bottom two plots of Figure~\ref{fig:synth_unc_plots}. Firstly, BSG-N-FD still performs the best in terms of both iterations and time, closely followed by StocBiO. BSG-H is able to achieve a similar function value to those achieved by~BSG-N-FD and~StocBiO, but after several more iterations and much more time. The decrease in the performance of~BSG-H can be explained by sample sizing of the stochastic Hessian matrices, and it is well known that Hessians require more samples than gradients in the presence of noise~\cite[Section 6.1.1]{LBottou_FECurtis_JNocedal_2018}. In particular, BSG-H is very sensitive to the value of the Hessian standard deviation, and its performance significantly degrades for values much larger than~$0.05$. Finally, BSG-1 and~DARTS seem to be relatively robust to the noise and exhibit similar performance to the deterministic case in terms of iterations. However, in terms of time, BSG-1 outperforms~BSG-H and is even able to achieve a similar level of accuracy to BSG-N-FD and StocBiO.
 

\subsubsection{Results for the LL linearly constrained case}\label{sec:lin_con_results}
    In the experiments for the~LL linearly constrained version of problem~\eqref{prob:synthetic_prob}, the~LL constraint set~$Y(x) = Y$ was defined by the following~$\vert I \vert$ linear inequality constraints in~$y$
    \begin{equation}\label{linear_constr_y}
        Y \; = \; \{y ~|~ Wy \;\leq\; s\},
    \end{equation}
    where~$W\in\mathbb{R}^{\vert I \vert\times m}$ and~$s\in\mathbb{R}^{\vert I \vert}$ were both randomly generated according to a uniform distribution, respectively.
    In this section, we focus on~LL linear constraints in~$y$ because this allows us to compare the performance of the~BSG algorithms developed in this paper against~SIGD, which is only designed for handling these types of constraints (as mentioned in Subsection~\ref{sec:BSD}).
    We again considered a dimension of~300 at both the upper and lower levels (i.e., $n=m=300$) along with~$\vert I \vert=50$  constraints, with~$H_1$, $H_2$, $H_3$, and~$H_4$ chosen in the same manner as in Subsection~\ref{sec:unc_results}.

    In the stochastic case, we added Gaussian noise with mean~0 to the gradients and Hessians as in Subsection~\ref{sec:unc_results} and also to the Jacobians~$\nabla_{x} c_I$, $\nabla_{x} c_E$, $\nabla_{y} c_I$, and $\nabla_{y} c_E$. Note that we did not add noise to~$\nabla_{yx}^2 c_i$ and~$\nabla_{yy}^2 c_i$, with~$i \in I \cup E$, because they are null matrices in the~LL linearly constrained case. The value of the standard deviation was chosen from~$\{0.005,0.05\}$ (referred to as the ``low'' and ``high'' values, respectively) for the stochastic Hessian matrices, $\nabla_{xy}^2 f_{\ell}$ and~$\nabla_{yy}^2 f_{\ell}$, and was set to~$0.5$ for all the other stochastic estimates. Using two different values for the standard deviation of the stochastic Hessian allows us to demonstrate the impact of different noisy Hessian estimates on the performance of~BSG-H and~SIGD, which use second-order derivatives.
    Regarding the stepsizes~$\alpha_u$ and~$\alpha_{\ell}$, in the deterministic case, we used~$\underline{s}_u = 2$, $\bar{s}_u = 4$, $\underline{s}_{\ell} = 3$, and~$\bar{s}_{\ell} = 5$ for~BSG-N-FD and~BSG-H and~$\underline{s}_u = \underline{s}_{\ell} = 2$ and~$\bar{s}_u = \bar{s}_{\ell} = 4$ for~SIGD. For the stochastic case with Hessian standard deviation equal to~0.005, we used the same~UL and~LL stepsizes for SIGD as in the deterministic case. However, we changed to values of~$\underline{s}_{\ell} = 7$ and~$\bar{s}_{\ell} = 8$ for BSG-N-FD along with~$\underline{s}_{\ell} = 3$ and~$\bar{s}_{\ell} = 5$ for~BSG-H. When the Hessian standard deviation was equal to~0.05, we used the same~UL and~LL stepsizes as in the deterministic case, with the exception of~$\underline{s}_{u} = 5$ and~$\bar{s}_{u} = 7$ for~SIGD. We do not include any results for BSG-1 and for BSG-H when the Hessian standard deviation is~0.05 because we were unable to find stepsizes that allowed the algorithms to converge.

    \begin{figure}
    \centering
          \includegraphics[scale=0.44]{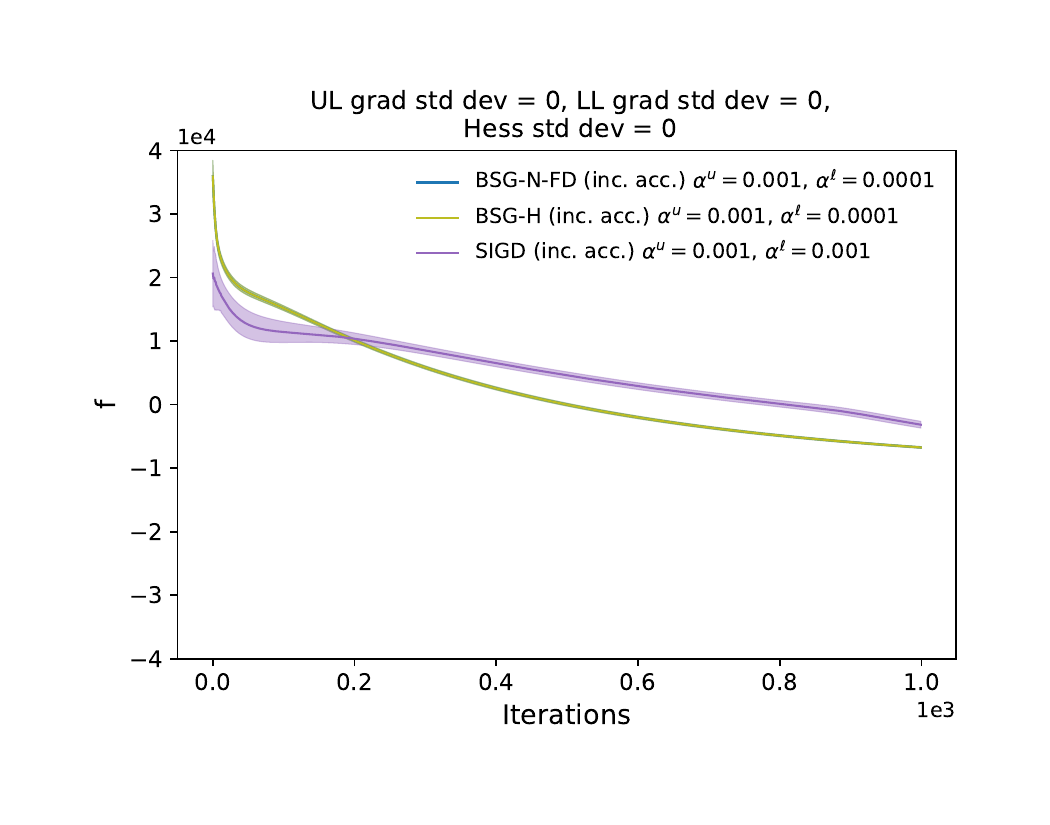}
          \includegraphics[scale=0.44]{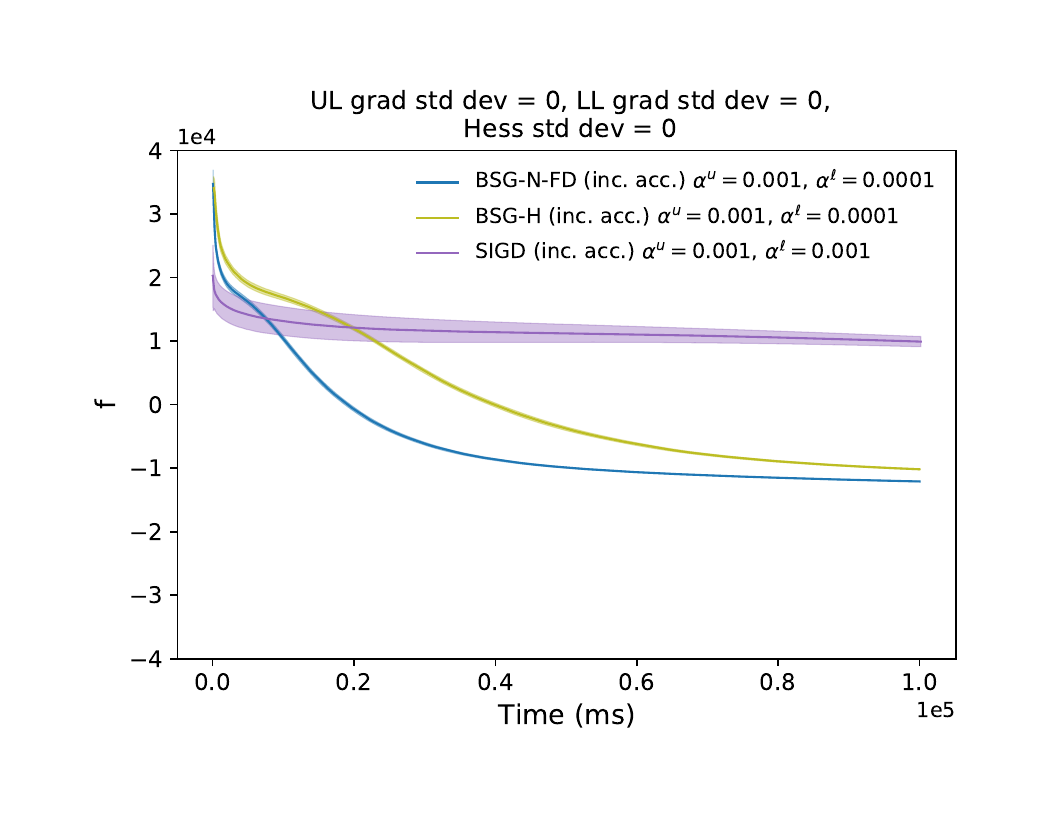}
          \includegraphics[scale=0.44]{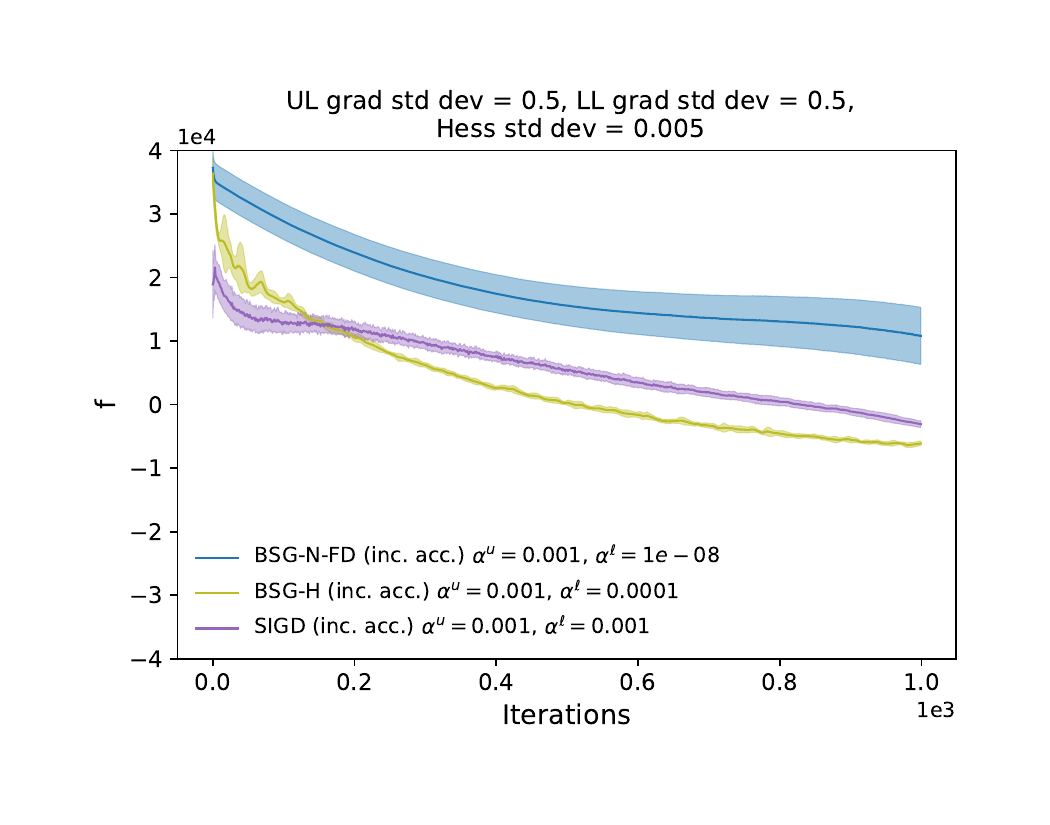}
          \includegraphics[scale=0.44]{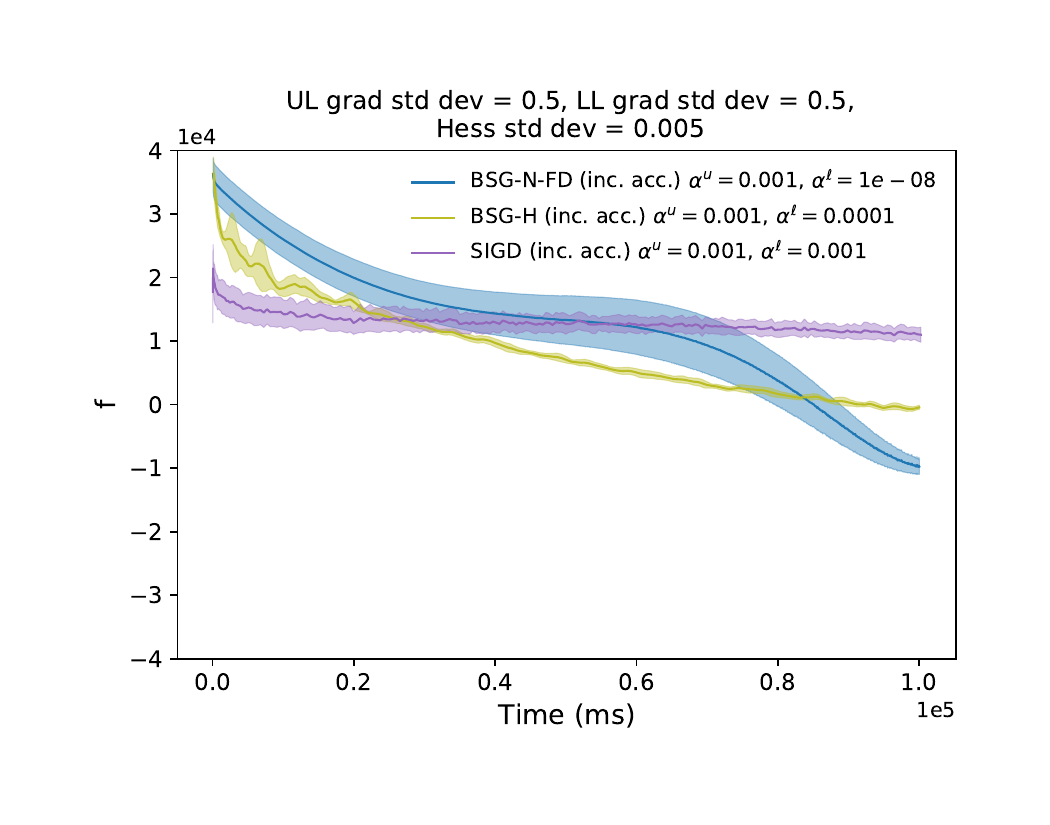}
          \includegraphics[scale=0.44]{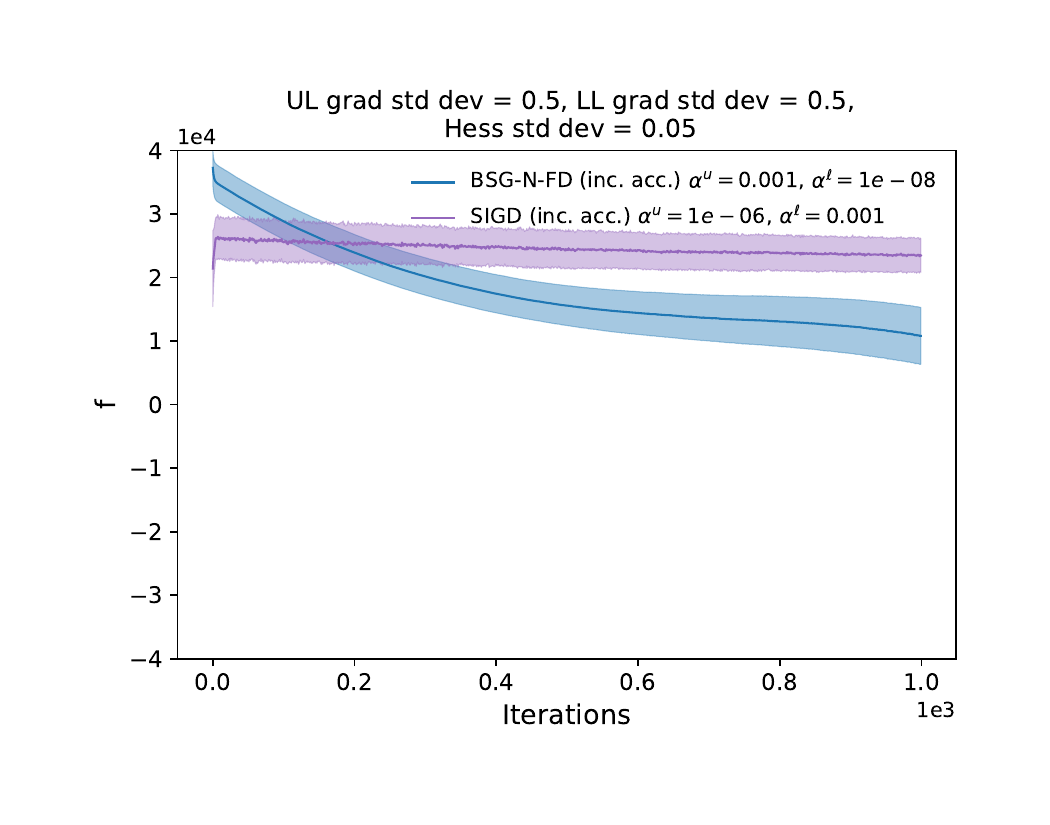}
          \includegraphics[scale=0.44]{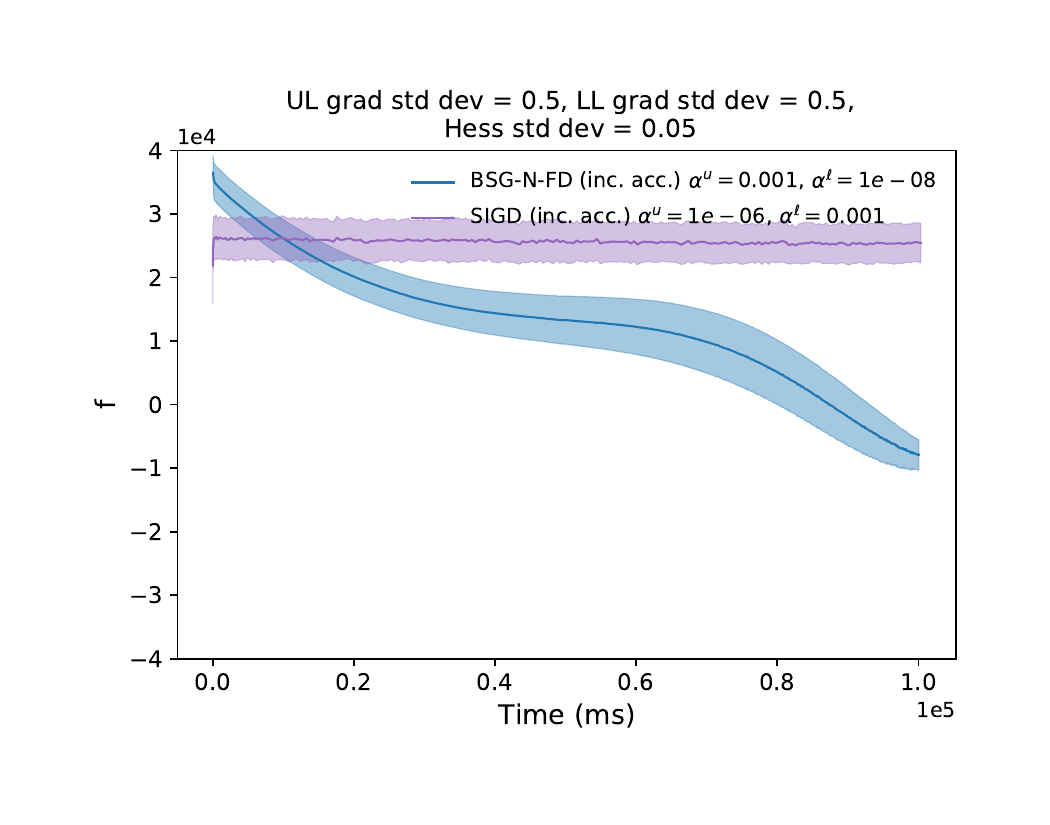}
        \caption{Numerical results of the~BSG-N-FD, BSG-H, and~SIGD algorithms on problem~\eqref{prob:synthetic_prob} with linear constraints in $y$ defined by~\eqref{linear_constr_y} in terms of both iterations and time (in milliseconds).}
         \label{fig:lin_y_constr_results}
\end{figure}

    Starting with the deterministic results displayed in the top two plots of Figure~\ref{fig:lin_y_constr_results}, we can clearly see that both~BSG-N-FD and~BSG-H outperform~SIGD in terms of iterations and time. In fact, BSG-N-FD and~BSG-H yield the exact same performance in terms of iterations, resulting in overlapping lines. 
    
    In the stochastic setting with low Hessian standard deviation displayed in the two middle plots of Figure~\ref{fig:lin_y_constr_results}, we can see that~BSG-H is still able to outperform~SIGD both in terms of iterations and time, despite its inferior performance compared to the deterministic case. Although BSG-N-FD has the worst performance here in terms of iterations, it is able to outperform both~BSG-H and~SIGD in the long run due to its superior efficiency in terms of time. Looking now at the setting with high Hessian standard deviation displayed in the two bottom plots of Figure~\ref{fig:lin_y_constr_results}, we can see that~BSG-N-FD is clearly the superior method.
    Similar to the behavior of BSG-H, SIGD is no longer able to converge. By contrast, BSG-N-FD is less affected by the noise, yielding the best performance.
    

\subsubsection{Results for the LL quadratically constrained case}\label{sec:quad_constr_results}
    In the experiments for the~LL quadratically constrained version of problem~\eqref{prob:synthetic_prob}, the~LL constraint set~$Y(x)$ was defined by the following~$\vert I \vert$ quadratic inequality constraints
    \begin{equation}\label{quad_constr}
        Y(x) \;=\; 
            \begin{cases}
                y^\top Q_1^{(1)} y \;+\; x^\top Q_2^{(1)} y \; \leq \; s^{(1)},\\
                y^\top Q_1^{(2)} y \;+\; x^\top Q_2^{(2)} y \; \leq \; s^{(2)},\\
                \;\;\;\;\;\;\;\;\;\;\;\;\;\;\;\;\;\;\;\;\;\vdots\\
                y^\top Q_1^{(\vert I \vert)} y \;+\; x^\top Q_2^{(\vert I \vert)} y \; \leq \; s^{(\vert I \vert)},\\
            \end{cases}
    \end{equation}
    where~$Q_1^{(i)}\in\mathbb{R}^{m\times m}$, $Q_2^{(i)}\in\mathbb{R}^{n\times m}$, and~$s^{(i)} \in \mathbb{R}$, for all $i\in\{1,2,...,\vert I \vert\}$, were all randomly generated according to a uniform distribution. We again considered a dimension of 300 at both the upper and lower levels (i.e., $n=m=300$) along with $\vert I \vert=5$ constraints, with $H_1$, $H_2$, $H_3$, and $H_4$ chosen in the same manner as in subsection~\ref{sec:unc_results}.

    We now present numerical results for the~BSG methods developed in this paper on a deterministic and two stochastic versions of problem~\eqref{prob:synthetic_prob} with constraints defined by~\eqref{quad_constr}, again testing two different levels of Hessian noise (with standard deviation values chosen from $\{0.005,0.5\}$) as in Subsection~\ref{sec:lin_con_results}. To the best of our knowledge, there do not exist any other bilevel stochastic algorithms that can handle general nonlinear constraints in the~LL problem, specifically quadratic constraints in this case, and as a result, the following numerical experiments are the first for this type of problem. We added noise to each gradient, Jacobian, and Hessian (including~$\nabla_{xy}^2 c_i$ and~$\nabla_{yy}^2 c_i$, for all~$i \in I \cup E$), as described in Subsection~\ref{sec:unc_results}.
    Regarding the stepsizes~$\alpha_u$ and~$\alpha_{\ell}$, we used~$\underline{s}_u = 2$, $\bar{s}_u = 4$, $\underline{s}_{\ell} = 6$, and~$\bar{s}_{\ell} = 8$ for both~BSG-N-FD and~BSG-H in both the deterministic and stochastic cases.

    \begin{figure}
        \centering
              \includegraphics[scale=0.44]{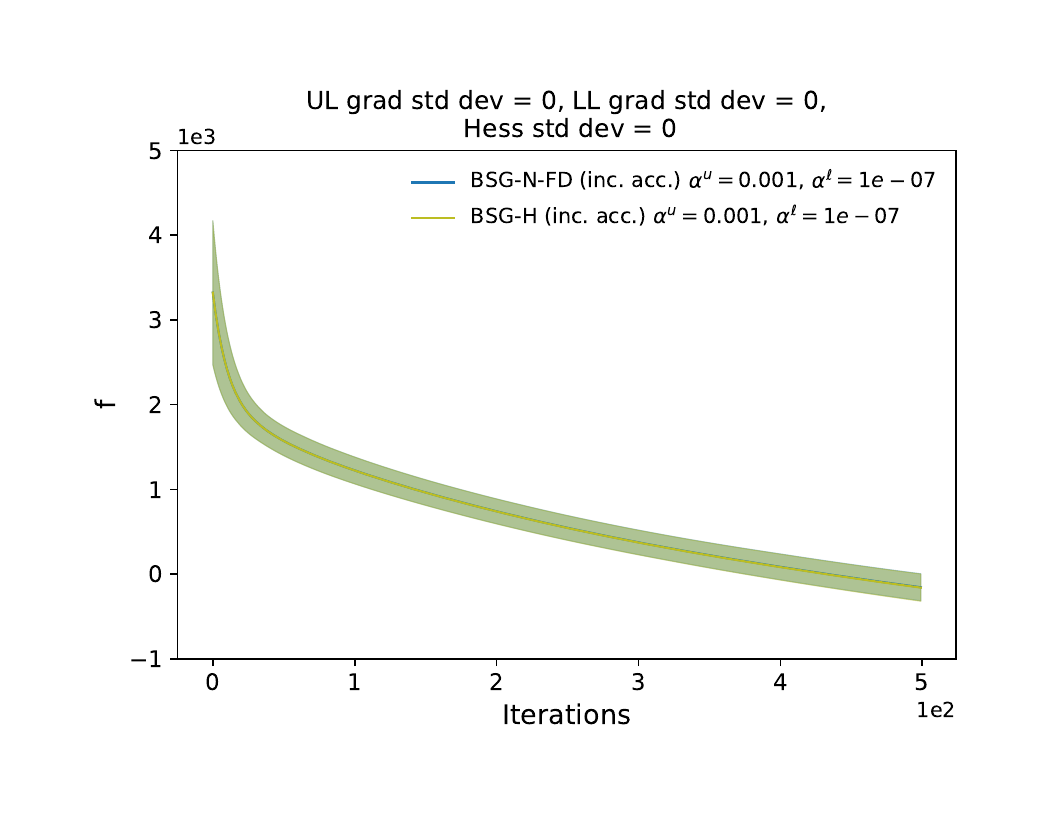}
              \includegraphics[scale=0.44]{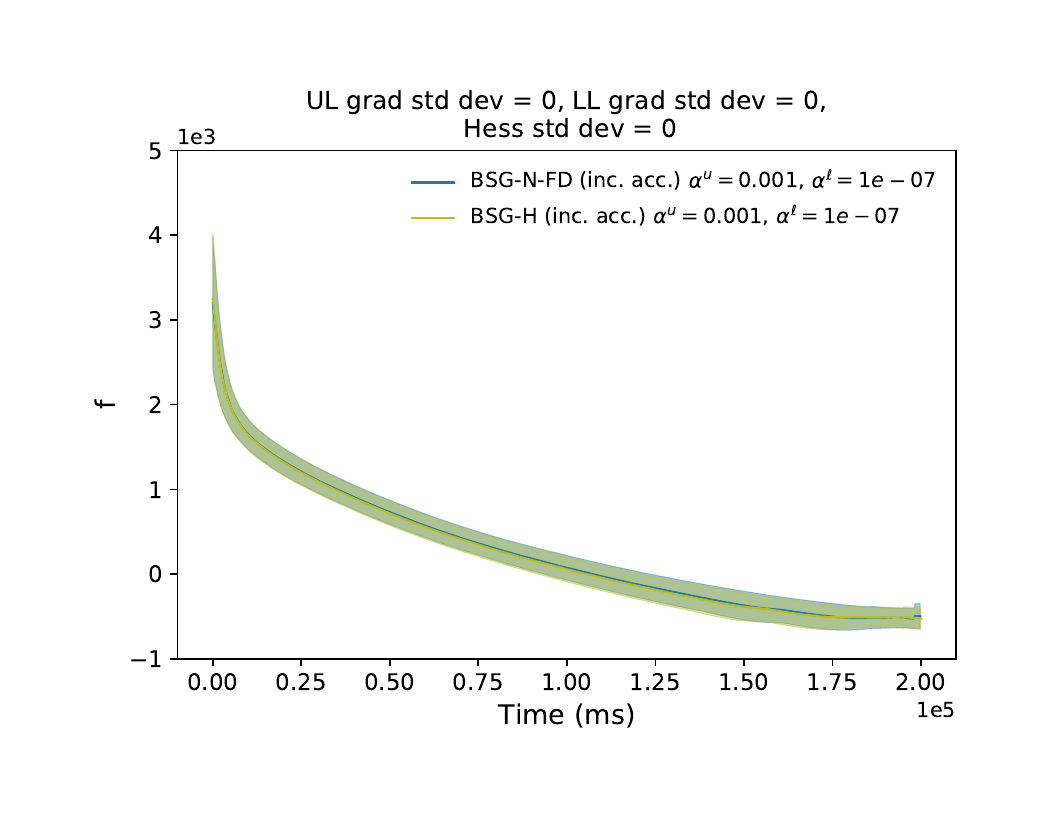}
              \includegraphics[scale=0.44]{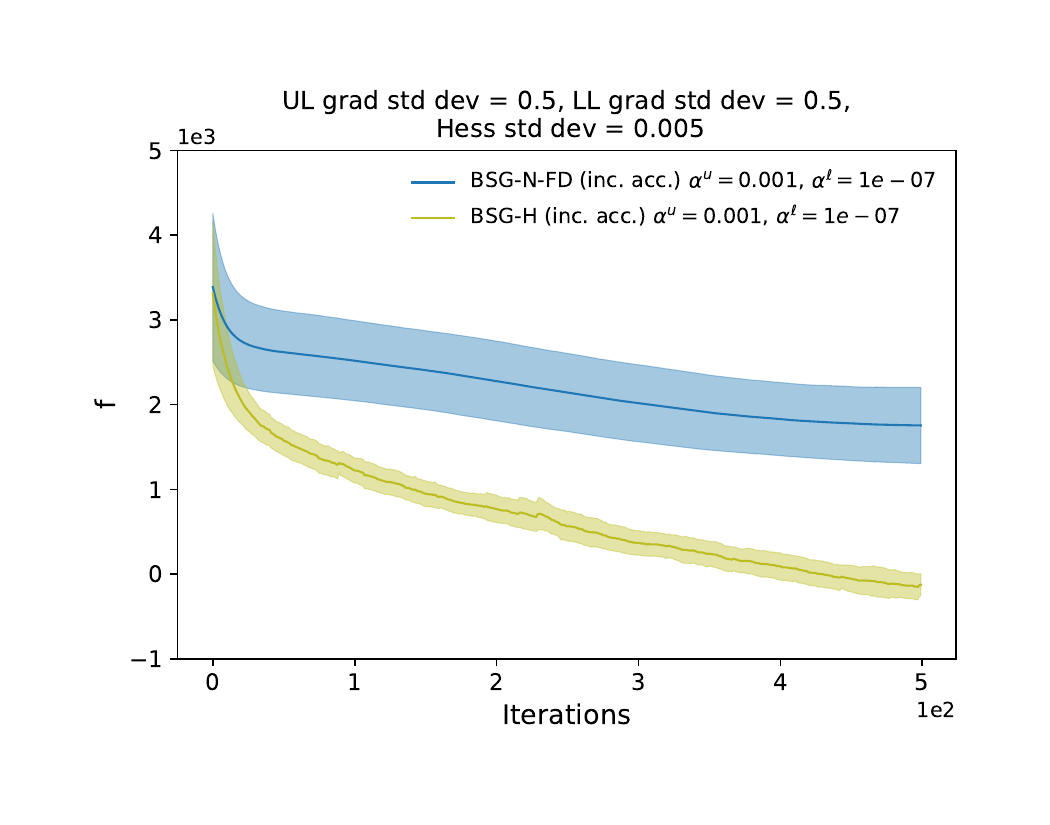}
              \includegraphics[scale=0.44]{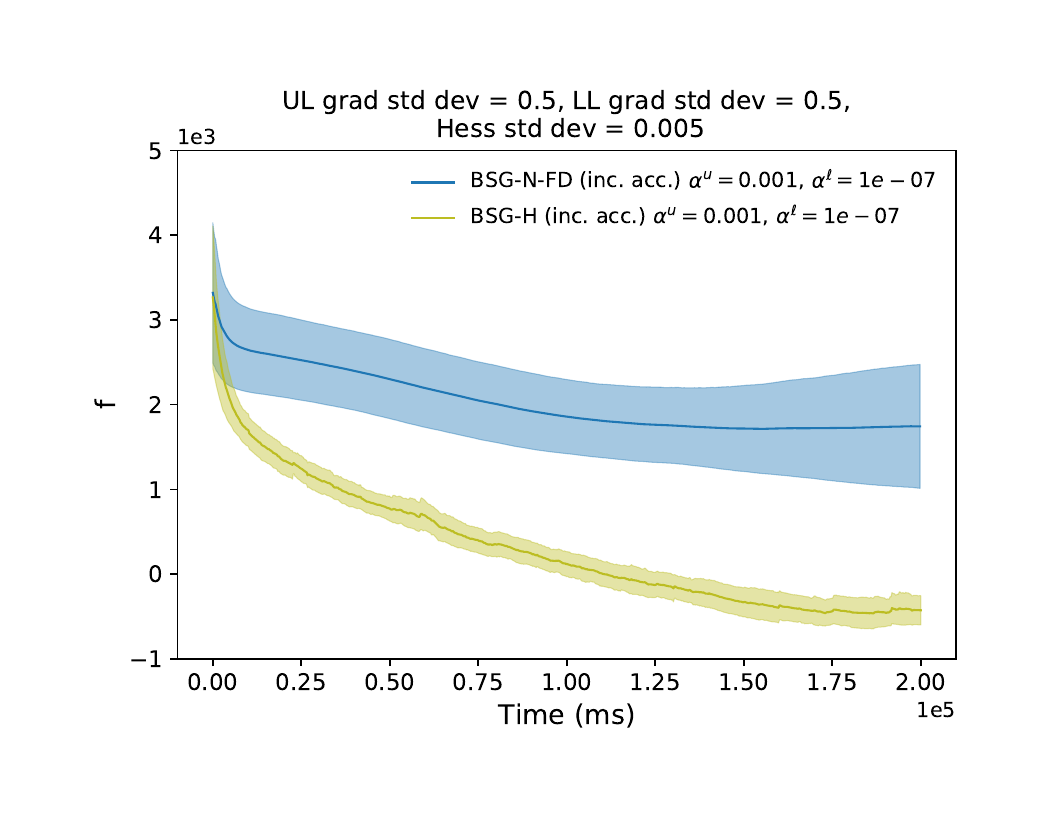}
              \includegraphics[scale=0.44]{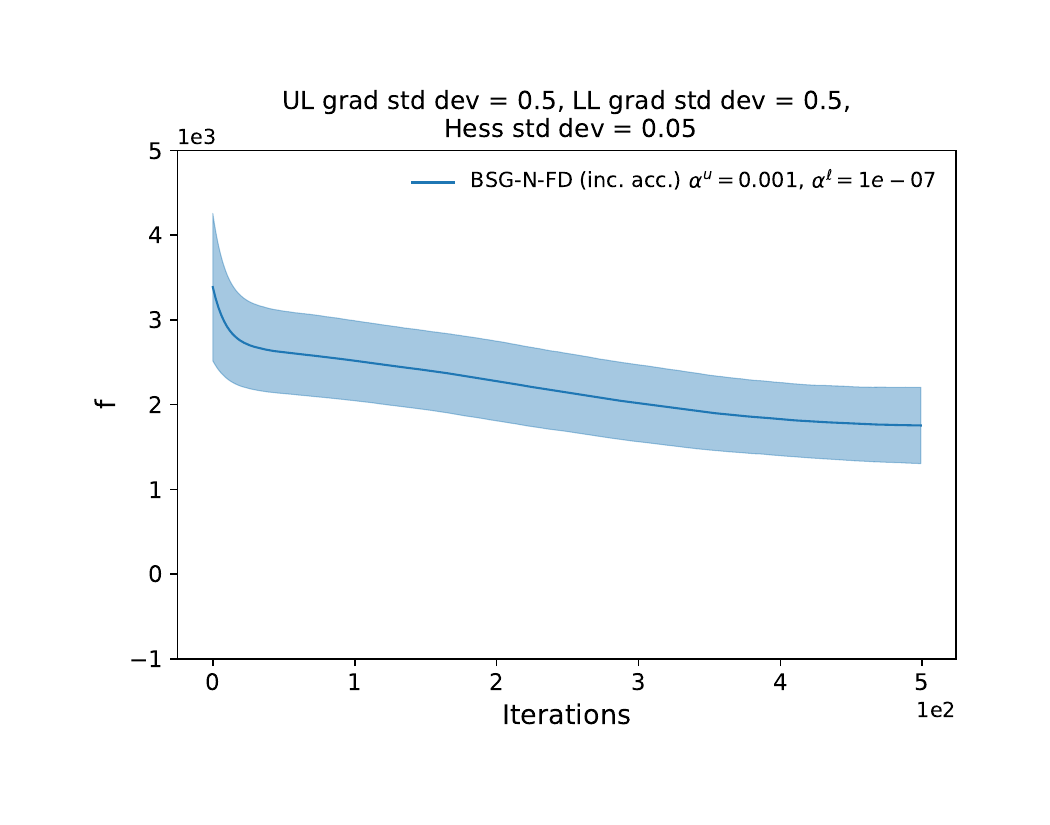}
              \includegraphics[scale=0.44]{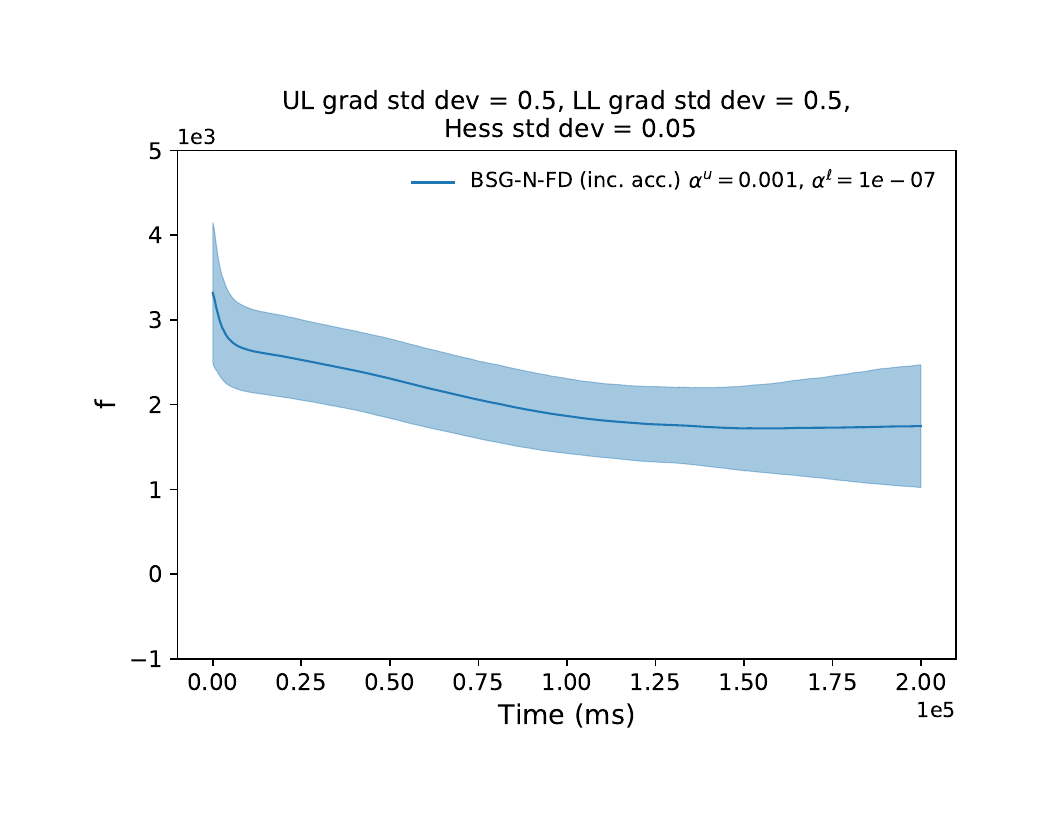}
            \caption{
            Numerical results of the BSG-N-FD and BSG-H algorithms on problem~\eqref{prob:synthetic_prob} with quadratic constraints defined by~\eqref{quad_constr} in terms of both iterations and time (in milliseconds).
             }\label{fig:quad_constr_results}
    \end{figure}

    Starting with the deterministic results displayed in the top two plots of Figure~\ref{fig:quad_constr_results}, we can clearly see that~BSG-N-FD and BSG-H have the same exact performance in terms of iterations and time. It would be expected that~BSG-N-FD yields better results in terms of time due to the efficiency we noted in Subsection~\ref{sec:lin_con_results}. In fact, this still holds true, except that this time the differences are not as substantial, and they would only become more visually apparent when allowing the algorithms to run for much longer.
    
    In the stochastic setting with a low level of Hessian noise (the middle two plots), we notice that~BSG-N-FD is impacted by the noise, while~BSG-H is still able to retain almost the same behavior as in the deterministic case. Although~BSG-H looks very favorable here, the stochastic setting with a high level of Hessian noise (the bottom two plots) shows similar results to the linearly constrained case in Subsection~\ref{sec:lin_con_results} (bottom two plots of Figure~\ref{fig:lin_y_constr_results}). Specifically, the performance of~BSG-N-FD remains unchanged, while we were not able to find stepsize values that allowed BSG-H to converge.


    \subsection{Continual learning}\label{sec:continual_learning}
    
    We are going to use instances of Continual Learning (CL) as practical stochastic bilevel problems to test the performance of BSG-N-FD, BSG-1, StocBiO, and DARTS. CL was briefly described in Section~\ref{sec:introduction}, and is now introduced in more detail. Let us denote a whole features/labels dataset by ${\cal D} = \{ (\textbf{u}_j, \textbf{v}_j), \, j \in \{ 1, \ldots , N\} \}$, consisting of $N$~pairs of a feature vector $\textbf{u}_j$ and the corresponding true label~$\textbf{v}_j$. For any data point~$j$, the classification is deemed correct if the right label is predicted. To evaluate the loss incurred when using the prediction function $\phi (\textbf{u}; \theta)$, which in this section is supposed to be a DNN, we use a loss function $\ell(\phi (\textbf{u}; \theta), \textbf{v})$.  
    
    The goal of CL is to minimize the prediction error over a sequence of tasks that become available one at a time. Among the many different formulations proposed for CL, hierarchical objectives have been used in~\cite{QPham_etal_2020,AShaker_etal_2020}. In this section, we present an incremental setting where each task is available as a subset of samples, similar to the formulation in~\cite{QPham_etal_2020}. Given $t \in \{1, \ldots, T\}$, let ${\cal D}_t$ be the set of samples for a new task $T_t$, which can be split into a training set $D^t_{\tr}$ and a validation set $D^t_{\va}$. Moreover, let us split the parameters $\theta_t$ of the model $\phi (\textbf{u}; \theta_t)$ into two subvectors, $\lambda_t$ and $\delta_t$, whose roles are to give us flexibility in minimizing the classification error on the training and validation data along the sequence of tasks. According to \cite{LFranceschi_etal_2018}, when using a neural network as a prediction function, a reasonable strategy is to choose $\lambda_t$ and $\delta_t$ as the vectors of weights in the hidden and output layers, respectively.  
    
    To solve the overall CL problem, one starts from the first task $T_1$ and, after an arbitrary number of iterations or amount of time, we include in the problem the second task $T_2$. One reiterates this procedure until all the tasks have been added to the problem. Let us now suppose that one has already added $t$ tasks. At this stage, the goal of the UL and LL problems is to determine the values of $\lambda_t$ and $\delta_t$ that ensure a small classification error on $T_t$ and on all the previous tasks~$T_{i}$, with $i < t$. To this end, the UL problem determines~$(\lambda_t,\delta_t)$ by minimizing the prediction error on~$D^t_{\val} = \cup_{i \le t} D^i_{\va}$, which is composed of the data sampled from the validation sets associated with the current and previous tasks. Similarly, the LL problem determines~$\delta_t$ by minimizing the error on $D^t_{\train} = \cup_{i \le t} D^i_{\tr}$. 
    Note that at each stage one solves a different problem since the objective functions of the UL and LL change as new tasks are included in the problem. The formulation of the problem solved at stage $t$, with~$t \in \{1, \ldots, T\}$, can be written as follows:
    \begin{equation}\label{prob:bilevel_continual_learning}
    	\begin{split}
    		\min_{(\lambda_t, \, \delta_t)} ~~ & f_u(\lambda_t, \, \delta_t) \; = \; \frac{1}{|D_{\val}^t|} \sum_{(\textbf{u}, \textbf{v}) \in \, D_{\val}^t} \ell(\phi (\textbf{u}; \lambda_t, \delta_t), \textbf{v}) \\
    		\mbox{s.t.}~~ & \delta_t \in \argmin_{\delta_t} ~~ f_{\ell}(\lambda_t, \, \delta_t) \; = \; \frac{1}{|D_{\train}^t|} \sum_{(\textbf{u}, \textbf{v}) \in \, D_{\train}^t} \ell(\phi (\textbf{u};  \lambda_t,\delta_t), \textbf{v}).
    	\end{split}
    \end{equation}
    
    We point out that the large dimension of the datasets usually considered in ML may prevent the use of the whole sets $D^i_{\tr}$ and $D^i_{\va}$ from previous tasks $i$'s, where $i < t$ and $t$ is the current task. In such cases, it may be necessary to resort to subsets $\bar D^i_{\tr} \subset D^i_{\tr}$ and $\bar D^i_{\va} \subset D^i_{\va}$, which we will not do in this paper given that our interest focuses on the solution of stochastic BLPs.
    
    
    Once a new task is included in the problem, the classification accuracy of the~DNN on the previous tasks tends to deteriorate, thus resulting in the well-studied phenomenon of \textit{catastrophic forgetting} \cite{IJGoodfellow_MMirza_DXiao_ACourville_YBengio_2013}, which can be alleviated by adding~LL inequality constraints to the lower level of problem~\eqref{prob:bilevel_continual_learning}.
    Such inequality constraints are inspired by~\cite{DLopezPaz_MARanzato_2017} and ensure that, at each stage, the current model outperforms the old model on all the previous tasks, thus preventing the deterioration of the classification accuracy when learning new tasks. In particular, for all $i < t$ and $t \ge 2$, we have
    \begin{equation}\label{eq:constraints_CL_catastrophic}
    \sum_{(\textbf{u}, \textbf{v}) \in \, D_{\tr}^i} \ell(\phi (\textbf{u};  \lambda_t,\delta_t), \textbf{v}) - \sum_{(\textbf{u}, \textbf{v}) \in \, D_{\tr}^i} \ell(\phi (\textbf{u};  \lambda_{t-1},\delta_{t-1}), \textbf{v}) \; \le \; 0.
    \end{equation}
    We point out that similar constraints were also used in the bilevel formulation proposed in~\cite{AShaker_etal_2020}, where the violation of the constraints is penalized. Note that, in general, the constraints in~\eqref{eq:constraints_CL_catastrophic} may be nonconvex when using a neural network as the prediction function~$\phi (\textbf{u};  \lambda_t,\delta_t)$.
    However, our theory does encompass nonconvexity of the LL constraints as long as the~LL SOSC is satisfied (see Assumption~\ref{ass:LL_assumptions_constr}).

\subsection{Results for continual learning instances}
	
    We now present numerical results comparing~BSG-N-FD and~BSG-1 against both~DARTS and~StocBiO on the CL problem~(\ref{prob:bilevel_continual_learning}) that was posed in Section~\ref{sec:continual_learning}. We also include numerical results on the~CL problem with LL constraints defined by~\eqref{eq:constraints_CL_catastrophic} for~BSG-N-FD and~BSG-1, as no other method applies in this case.
	In our implementation, we determine $\lambda_t$ (the UL variables) on the current problem by starting from the parameter values found from the previous problem. However, since each consecutive task increases the output space of the DNN, we entirely re-initialize $\delta_t$ (the LL variables) at the start of each new task (when first applying an~LL step) so that the model outputs are not biased from previous tasks.

    \begin{figure}
    \centering
        \includegraphics[scale=0.55]{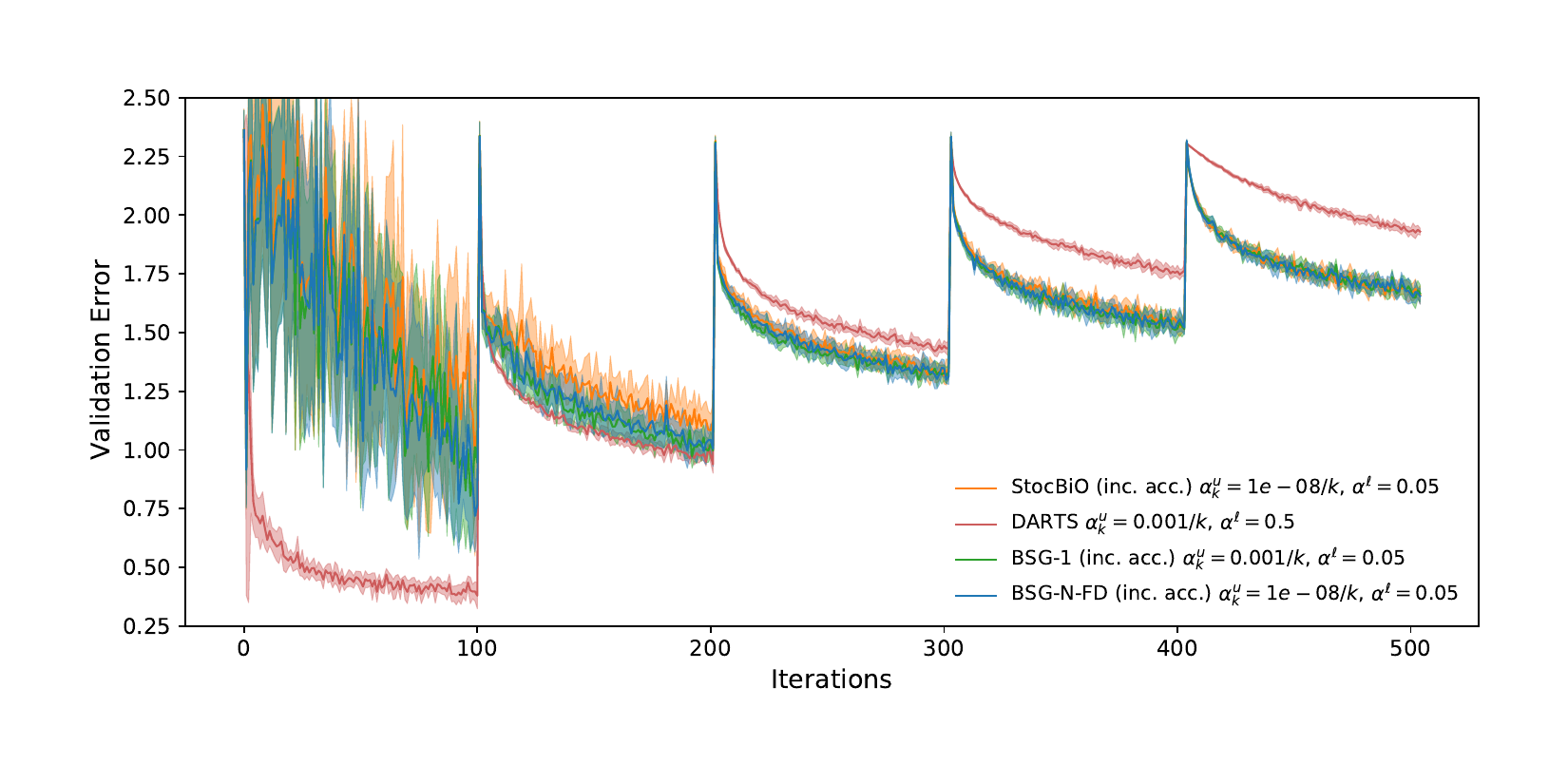}
        \includegraphics[scale=0.55]{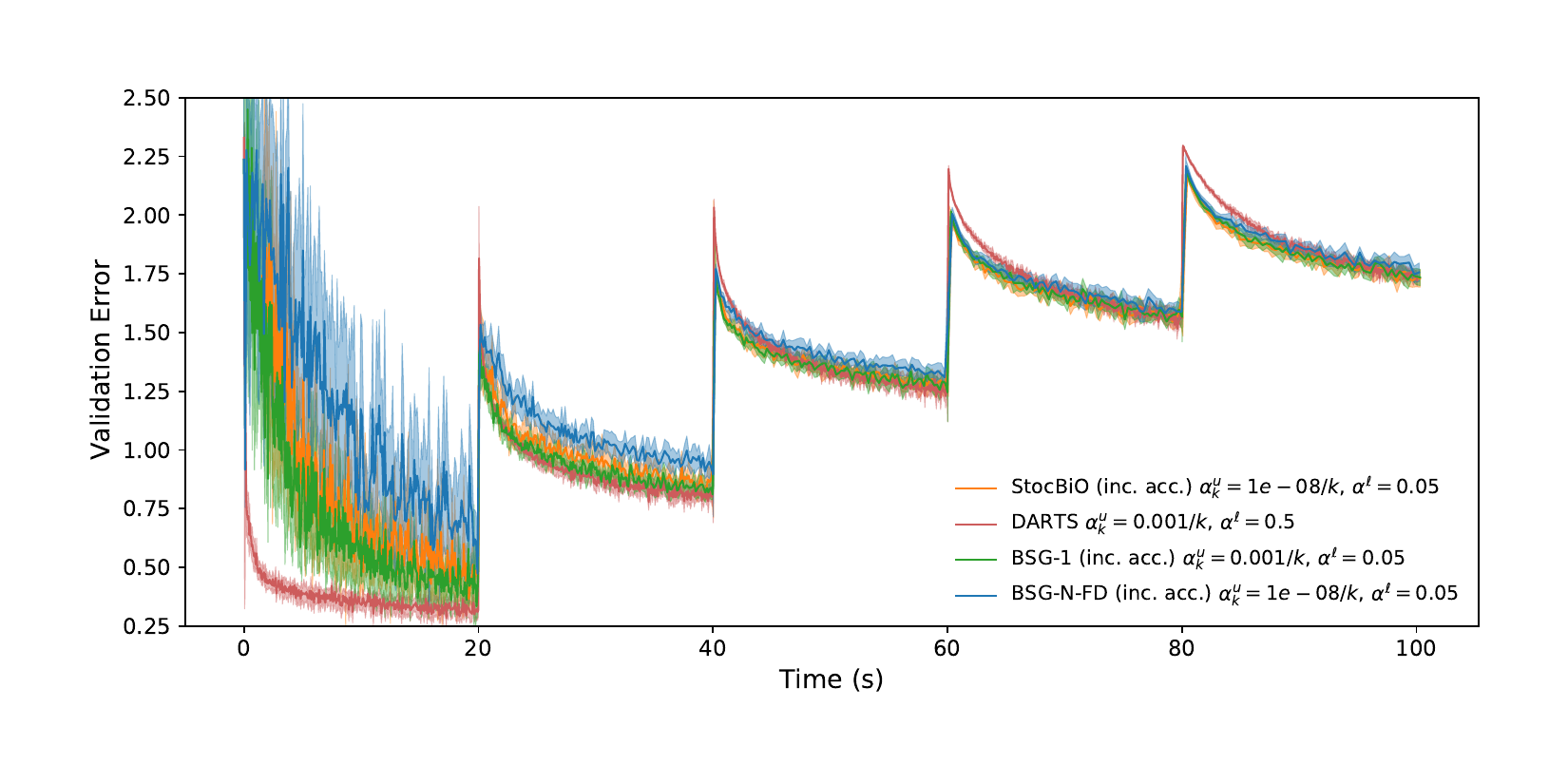}
        \caption{Comparison of BSG-N-FD, BSG-1, DARTS, and StocBiO on the~CL problem~(\ref{prob:bilevel_continual_learning}) in terms of both iterations (top plot) and time (bottom plot, in seconds).}\label{fig:CL_MNIST}
    \end{figure}
    

    \begin{figure}
    \centering
        \includegraphics[scale=0.55]{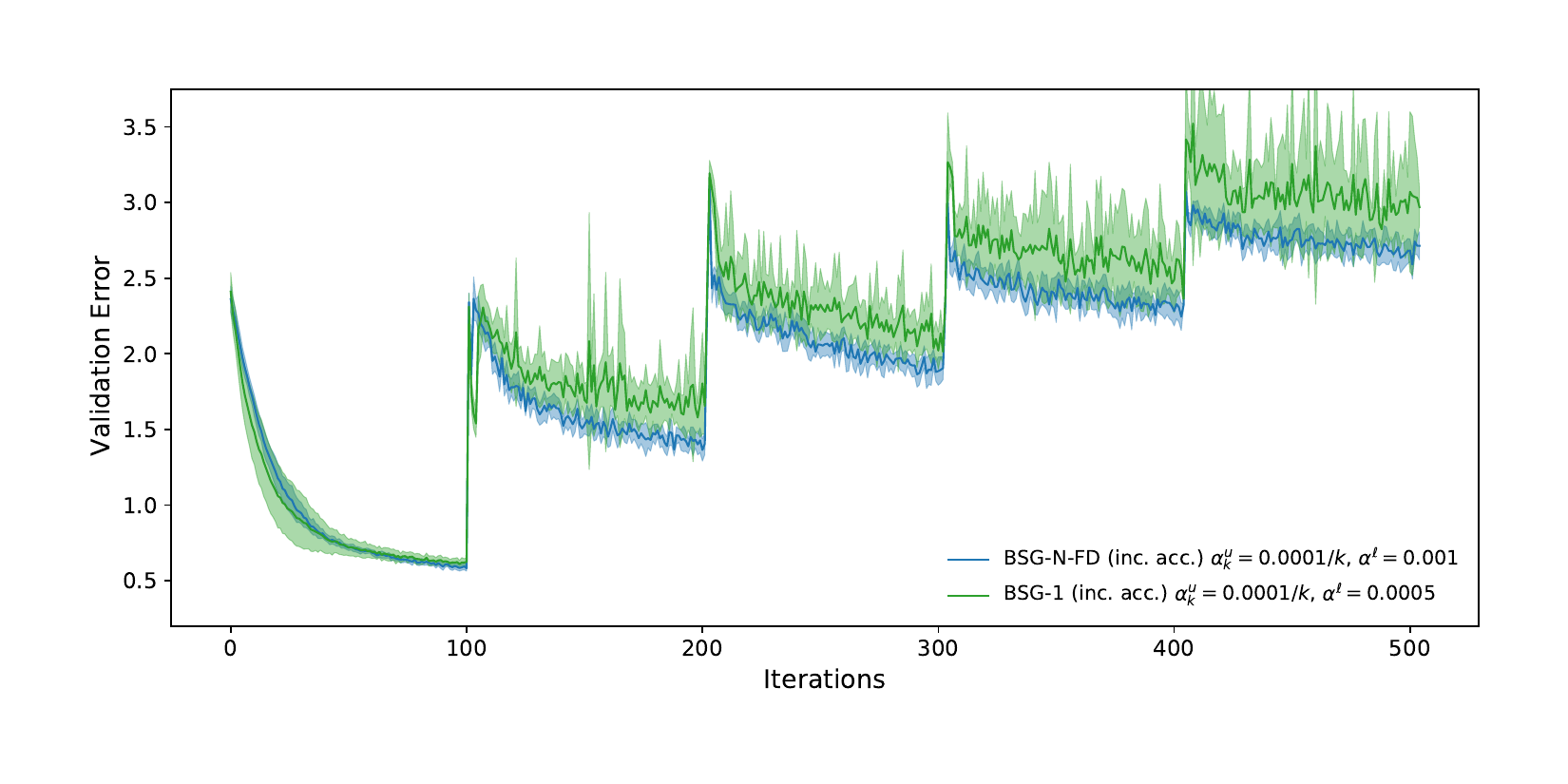}
        \includegraphics[scale=0.55]{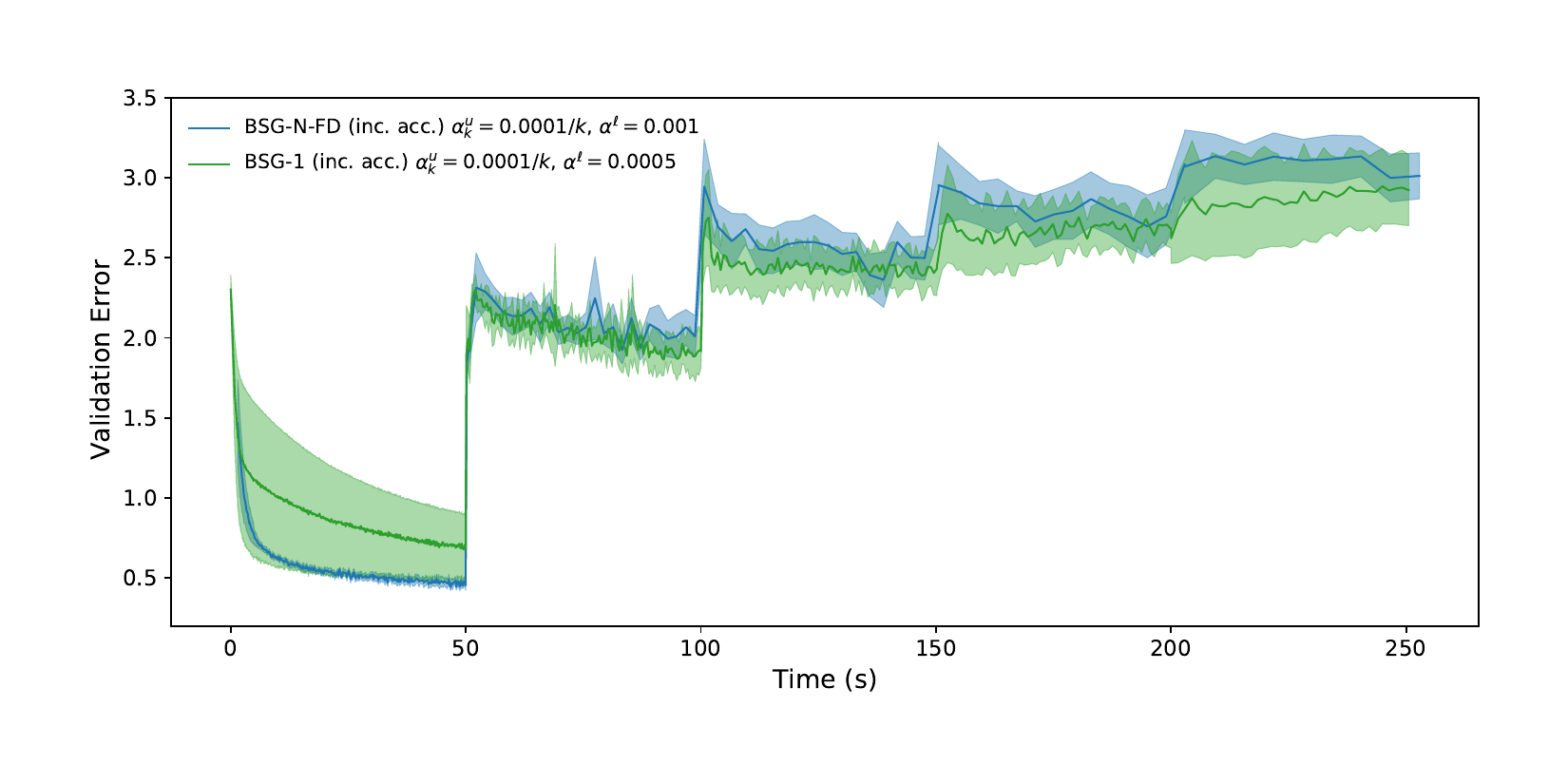}
        \caption{Comparison of BSG-N-FD and BSG-1 on the~CL problem~(\ref{prob:bilevel_continual_learning}) with constraints~\eqref{eq:constraints_CL_catastrophic} in terms of both iterations (top plot) and time (bottom plot, in seconds).}\label{fig:CL_CIFAR_constrained}
    \end{figure}
    
    In order to test our algorithm on a large-scale ML scenario, we chose the well-studied CIFAR-10 dataset~\cite{AKrizhevsky_2009}, which consists of a total of 60,000 colored images (32 $\times$ 32) of 10 different classes (i.e., airplanes, automobiles, birds, cats, deer, dogs, frogs, horses, ships, and trucks). The dataset is split into a training set that consists of 50,000 images and a testing set that consists of 10,000 images; however, for our experiments we only used the first set of images. We used a subset of 40,000 images for training and the remaining 9,999 images for validation (since one of the images had an issue, we removed it from the dataset). We solved five problems~(\ref{prob:bilevel_continual_learning}) with an increasing number of tasks from~$1$ to~$5$, where the first task datasets ($D_{\text{val}}^1$ and $D_{\text{train}}^1$) consist of only the images with class labels in $\{0,1\}$ (these correspond to airplanes and automobiles), the second task datasets ($D_{\text{val}}^2$ and $D_{\text{train}}^2$) consist of the images with class labels in $\{0,1,2,3\}$ (airplanes, automobiles, birds, and cats), etc., until the final task datasets ($D_{\text{val}}^5$ and $D_{\text{train}}^5$), which are the original training and validation sets and consist of all the class labels $\{0,1,...,9\}$. Further, we implemented a DNN with two convolutional layers, a max-pooling layer, and one linear fully-connected layer as our model. The network consisted of~19,392 and~163,840 weights in the hidden and output layers, respectively. For the UL and LL~problems, we have used batch sizes equal to~0.05\% and~0.01\% of the sizes of the current task's validation and training datasets, respectively. 
    
    All of the algorithms, with the exception of~DARTS in the~LL unconstrained case, were run while using an increasing accuracy strategy in the LL problem with an~$f_u$ difference threshold for increasing the number of~LL iterations equal to $10^{-2}$ (and a maximum limit of~30~LL iterations).
    The number of~UL iterations and running time (in seconds) were both used as metrics for the comparison. 
    As a loss function, we used the well-known binary cross-entropy loss.
    In all of the figures in this section, we plot the approximation~$f_u$ of the true function~$f$, as it is typically done in bilevel ML~\cite{MHong_etal_2020,JYang_KJi_YLiang_2021,DSow_KJi_YLiang_2021,FPedregosa_2016,JLorraine_PVicol_DDuvenaud_2019}.
    
    The results for the unconstrained LL case are reported in Figure~\ref{fig:CL_MNIST}. We are not reporting BSG-H because of the extremely high computational cost of dealing with second-order derivatives given the choices of DNN and dataset.
    We compare four algorithms, BSG-N-FD and BSG-1 against DARTS and StocBiO, using the best~UL decaying stepsize sequence~$\{\alpha_k^{u}\}_{k \in \mathbb{N}}$ and the best~LL fixed stepsize~$\alpha^{\ell}$ found for each algorithm. The step sizes for each algorithm were obtained by performing grid searches over the following sets: $\alpha^u \in \{5 \cdot 10^{-3}/k, 10^{-3}/k, 5 \cdot 10^{-4}/k\}$ for~BSG-1, and~DARTS and~$\alpha^u \in \{10^{-7}/k, 10^{-8}/k\}$ for~BSG-N-FD and StocBiO; $\alpha^\ell \in \{10^{-1}, 5 \cdot 10^{-2}, 10^{-2}\}$ for~BSG-N-FD, BSG-1, and~StocBiO, $\alpha^\ell \in \{1, 5 \cdot 10^{-1}, 10^{-1}\}$ for~DARTS.
    Again, for StocBiO, we set the constant~$C_0$ introduced in Subsection~\ref{sec:BSD} to~0.05 and the parameter~$q$ introduced in Subsection~\ref{sec:inexact_adjoint} to~2, which led to the best results.
    By the nature of the~CL problem, we expect to see five separate ``jumps'' in the validation error (the UL objective function) indicating the start of a new task. Among the four algorithms, BSG-N-FD, BSG-1, and~StocBiO have similar performance in terms of iterations. In particular, they perform the best on all tasks excluding the first two, with StocBiO performing the worst on task two. On the first two tasks, DARTS performs the best in terms of both iterations and time, while~BSG-N-FD, BSG-1, and StocBiO experience some initial noise. All methods seem to perform similarly on the last three tasks, except for DARTS, which seems to fall behind on the last two tasks.
    
     
     Lastly, in Figure~\ref{fig:CL_CIFAR_constrained}, we provide numerical results for BSG-N-FD and BSG-1 on the~CL problem~\eqref{prob:bilevel_continual_learning} when considering~LL constraints~\eqref{eq:constraints_CL_catastrophic}. 
     In accordance with the procedure described in Subsection~\ref{sec:stochasticBSG1} for the~LL constrained case, approximate Lagrange multipliers are obtained at each iteration by solving the corresponding~KKT system with the linear conjugate gradient method (with maximum number of iterations equal to~3 and tolerance equal to~$10^{-4}$). 
     The system in~\eqref{constr_direction_approx} is solved by using the~GMRES method with maximum number of iterations equal to~$3$ when running~BSG-N-FD and~50 when running~BSG-1, with a tolerance equal to~$10^{-4}$.
     In a similar manner for the CL unconstrained case, we chose the stepsizes for each algorithm by performing the following grid searches: $\alpha^u \in \{5\cdot 10^{-4}/k, 10^{-4}/k, 5 \cdot 10^{-5}/k\}$ for both algorithms, $\alpha^\ell \in \{5 \cdot 10^{-3},10^{-3},5 \cdot 10^{-4}\}$ for~BSG-N-FD, and~$\alpha^\ell \in \{10^{-3}, 5 \cdot 10^{-4}, 10^{-4}\}$ for~BSG-1. 
     Referring to Figure~\ref{fig:CL_CIFAR_constrained}, BSG-N-FD performs the best in terms of iterations on all tasks except for the first, on which both algorithms perform similarly. In terms of time, BSG-1 yields slightly superior performance compared to~BSG-N-FD besides the first task. It should be noted that both algorithms seem to plateau after the second task in terms of time as the number of iterations in each consecutive task decreases substantially. This is due to the amount of time allotted to each task.
     The results in Figure~\ref{fig:CL_CIFAR_constrained} demonstrate that our BSG methods are able to solve large-scale bilevel optimization problems with nonlinear constraints in the LL problem, and similar to Section~\ref{sec:quad_constr_results}, these results are the first of their kind for this type of problem, to the best of our knowledge.

\section{Concluding remarks and future work}\label{sec:conclusions}

In this paper, we proposed a general framework for bilevel stochastic gradient~(BSG) methods that applies to both the LL unconstrained and constrained cases, we provided a corresponding convergence theory that allows for any inexactness in the calculation of adjoint gradients and that also rigorously covers the inexact solution of the~LL problem, and we introduced practical~BSG methods for large-scale bilevel optimization problems (BSG-N-FD and BSG-1). The numerical results showed that BSG-N-FD, which is consistent with the theory, performs well on the synthetic quadratic bilevel problem. On the continual learning instances, BSG-N-FD has a similar performance to the practical algorithms~BSG-1 and~StocBiO in terms of iterations and is slightly outperformed in terms of time.

The results on the ML instances considered in this proposal suggest that our \mbox{BSG} methods have the potential to perform well on the unconstrained bilevel formulations of NAS, which in the literature are mostly still tackled by using DARTS when a continuous relaxation of the (discrete) search space is used~\cite{PRen_etal_2021}. We point out that using finite differences like in~BSG-N-FD or rank-1 Hessian approximations like in~BSG-1 is crucial to allow the application of the BSG method to NAS, which would not be possible otherwise due to the extreme dimensions of the resulting bilevel problems. Moreover, the fact that our \mbox{BSG} methods can solve bilevel optimization problems with constrained LL problems paves the way for the solution of new NAS formulations. In particular, one could think of including in the LL problem constraints that help the model avoid overfitting~\cite{SNRavi_TDinh_VSLokhande_VSingh_2019} or constraints that depend on the specific learning instances considered~\cite{SSangalli_etal_2021}. Also left for future work are variance reduction techniques, which can be incorporated into our \mbox{BSG} methods to ensure faster convergence, as already proposed in~\cite{TChen_YSun_WYin_2021,JYang_KJi_YLiang_2021}.

\section*{Acknowledgments}
This work is partially supported by the U.S. Air Force Office of Scientific Research (AFOSR) award FA9550-23-1-0217.


\appendix

\section{Proposition \ref{prop:sensitivity_BSG_dir}}\label{appendix:sensitivity_BSG_dir}

\begin{proof}
Let us first prove~\eqref{Ass:33}. There are two cases to consider: inexact adjoint system and truncated Neumann series. In both we will use the fact that when~$B_1$ and~$B_2$ are non-singular,
\begin{equation}\label{eq:_P}
\Vert B_1^{-1} - B_2^{-1} \Vert \le \Vert B_1^{-1} (B_1 - B_2) B_2^{-1} \Vert \le \Vert B_1^{-1} \Vert \Vert B_2^{-1} \Vert \Vert B_1 - B_2 \Vert.
\end{equation}

\noindent
{\it 1) Inexact adjoint system.} \\
    The approximate~BSG direction is~$d(D) \; = \; -(a - A B^{-1}\tilde{b})$, where~$\tilde{b} = b + \tilde{r}$ and~$\tilde{r}$ is the residual error due to the inexact solution of the ajoint equation (see Subsection~\ref{sec:inexact_adjoint}). Now, we have
    \[
    \|d(D_1) - d(D_2)\| \; = \; \|-a_1 + A_1B_1^{-1} \tilde b_1 + a_2 - A_2B_2^{-1} \tilde b_2\|.
    \]
    Adding and subtracting $A_1B_1^{-1} \tilde b_2$ and using the triangle inequality, we obtain
    \begin{alignat}{1}
        \|d(D_1) - d(D_2)\|  \; \leq \; \|a_1 - a_2\| + \|A_1B_1^{-1}\| \|\tilde b_1 - \tilde b_2 \| + \|\tilde b_2 \| \|A_1B_1^{-1} - A_2B_2^{-1} \|. \nonumber
    \end{alignat}
    Adding and subtracting $A_2B_1^{-1}$ in the last norm on the right, we have
    \begin{alignat}{2}
        \|d(D_1) - d(D_2)\|  \; \leq \;\;& \|a_1 - a_2\| + \|A_1\| \|B_1^{-1}\| \|\tilde b_1 - \tilde b_2 \| \nonumber \\ &+ \|\tilde b_2 \| \|A_2 \| \|B_1^{-1} - B_2^{-1} \| + \|\tilde b_2 \| \|B_1^{-1}\| \|A_1 - A_2 \|. \label{eq:eq0000_100}
    \end{alignat}
    Let us define the positive constants $C_1 = \max\{ C C_{\ell}$, $\bar{C} \bar{C}_{\ell}\}$, $C_2 = \max\{C^2 C_{\ell}^2, \bar{C}^2 \bar{C}_{\ell}^2\}$, $C_3 = \max\{C C_{\ell}^2, \bar{C} \bar{C}_{\ell}^2\}$, and $C_4 = \max\{C_{\ell}, \bar{C}_{\ell}\}$, where~$C$, $C_{\ell}$, $\bar{C}$, and~$\bar{C}_{\ell}$ are the constants introduced in Remark~\ref{remark}.
    From the assumptions of Proposition~\ref{prop:sensitivity_BSG_dir} and~\eqref{eq:_P}, and setting~$r_1 = \tilde{r}_1$ and~$r_2 = \tilde{r}_2$, there exists~$L = \max\{1,C_1,C_2 + C_3\|r_2\|,C_1 + C_4 \|r_2\|\}$, such that
    \begin{alignat}{2}
        \| d(D_1) - d(D_2) \| & \; \leq \; \|a_1 - a_2\| + C_1(\left\|b_1 - b_2\right\| + \left\|r_1 - r_2\right\|) \nonumber\\
        &\quad  + (C_2 + C_3\|r_2\|) \left\|B_1 - B_2\right\| + (C_1 + C_4 \|r_2\|) \left\|A_1 - A_2\right\| \nonumber\\
        & \; \leq \; L (\|a_1 - a_2\| + \|b_1 - b_2\| + \left\|B_1 - B_2\right\| + \left\|A_1 - A_2\right\|) + L \|r_1 - r_2\|. \label{eq:eq0000}
    \end{alignat}    
    From the equivalence of norms, there exists a positive constant~$\tilde{C}$ such that the proof of Part~1) is completed with~$L_{BSG} = L \max\{\tilde{C},1\}$.\\

\noindent{\it 2) Truncated Neumann series}. \\
    The approximate~BSG direction is~$ d(D) \; = \; -(a - A \mathscr{B} b)$, where $\mathscr{B} = B^{-1} - \tilde{R}$ and~$\tilde{R}$ is a residual matrix (see Subsection~\ref{sec:inexact_adjoint}). By repeating the reasoning used to prove Part~1) until~\eqref{eq:eq0000_100}, we arrive at
     \begin{alignat}{2}
     \| d(D_1) - d(D_2) \|
        & \; \leq \; \|a_1 - a_2\| + \left\|A_1\right\|\|\mathscr{B}_1\|\left\|b_1 - b_2\right\| \nonumber\\ 
        &\quad  \; + \; \left\|b_2\right\|\left\|A_2\right\|\|\mathscr{B}_1 - \mathscr{B}_2\| +\left\|b_2\right\| \|\mathscr{B}_1\|\left\|A_1 - A_2\right\|. \nonumber 
    \end{alignat}
    Let us define the positive constants $C_1 = \max\{ C C_{\ell}, \bar{C} \bar{C}_{\ell}\}$, $C_2 = \max\{C, \bar{C}\}$, $C_3 = \max\{C^2 C_{\ell}^2,$ \linebreak $\bar{C}^2 \bar{C}_{\ell}^2\}$, $C_4 = \max\{C^2, \bar{C}^2\}$, and $C_5 = \max\{C, \bar{C}\}$, where~$C$, $C_{\ell}$, $\bar{C}$, and~$\bar{C}_{\ell}$ are again the constants introduced in Remark~\ref{remark}. 
Therefore, by using the same arguments as in Part~1), but now with~$r_1 = \tilde{R}_1$ and~$r_2 = \tilde{R}_2$, there exists~$L = \max\{1,C_1 + C_2 \|r_1\|,C_3,C_4,C_1 + C_5 \|r_1\|\}$, such that 
\begin{alignat}{2}
    \| d(D_1) - d(D_2) \| & \leq \|a_1 - a_2\| + (C_1 + C_2 \|r_1\|) \left\|b_1 - b_2\right\| \nonumber\\
    &\quad  + C_3 \|B_1 - B_2\| + C_4 \|r_1 - r_2\| + (C_1 + C_5 \|r_1\|) \left\|A_1 - A_2\right\| \nonumber\\
    & \; \leq \; L (\|a_1 - a_2\| + \|b_1 - b_2\| + \left\|B_1 - B_2\right\| + \left\|A_1 - A_2\right\|) + L \|r_1 - r_2\|. \label{eq:eq0000_2}
\end{alignat}
Again, from the equivalence of norms, there exists a positive constant~$\tilde{C}$ such that the proof of Part~2) is completed with~$L_{BSG} = L \max\{\tilde{C},1\}$.

    We will now prove~\eqref{Ass:33-additional} for both the~LL unconstrained and constrained cases. From the equivalence of norms, there exists a positive constant~$\hat{C}$ such that 
     \begin{equation*}\label{eq:lipschitz_stoch_estimates}
        \Vert D_1 - D_2 \Vert \; \le \; \hat{C} (\Vert a_1 - a_2 \Vert + \Vert b_1 - b_2 \Vert + \Vert B_1 - B_2 \Vert + \Vert A_1 - A_2 \Vert),
    \end{equation*}
    where~$D_1 = D(x, w(x), \xi)$, $D_2 = D(x, \tilde{w}, \xi)$, and $(a_i,b_i,A_i,B_i)$ is the stochastic data in~$D_i$, for~$i \in \{1, 2\}$.    
    From the assumption on the Lipschitz continuity of the stochastic gradients, Hessians, and Jacobians in~$D_1$ and~$D_2$ for all~$\xi$, there exists a positive constant~$\bar{L}_{LL}$ such that~\eqref{Ass:33-additional} is satisfied.
\end{proof}

\section{Proposition~\ref{prop:smoothness_f_adj_grad}}\label{appendix:smoothness_f_adj_grad}

\begin{proof}
    We start by handling the LL unconstrained case. Taking the norm of equations~\eqref{eq:jacobian_formula} and \eqref{nabla_y} and from Remark~\ref{remark}, there exists a positive constant~$L_{\ell} = C C_{\ell}$ such that~$\|\nabla y(x) \| \leq L_{\ell}$ and~$\|\nabla w(x) \| \leq L_{\ell}$. It is well known that these two inequalities imply that~$y(x)$ and~$w(x)$ are Lipschitz continuous in~$x$ with constant~$L_{\ell}$ (see, e.g.,~\cite[Chapter~5]{ABeck_2017}). Therefore, one can write
    \begin{equation}\label{ineq2}
        \left\| y(x_1) - y(x_2)\right\| \; \leq \; L_{\ell} \| x_1 - x_2 \| \quad \text{ and } \quad \left\| w(x_1) - w(x_2)\right\| \; \leq \; L_{\ell} \| x_1 - x_2 \|.
    \end{equation}
    
    To prove~\eqref{eq:eq_lipschitz_cont_adj_grad} for the~LL unconstrained case, using the derivation followed for the proof of Proposition~\ref{prop:sensitivity_BSG_dir} until~\eqref{eq:eq0000}--\eqref{eq:eq0000_2} and considering~$r_1 = r_2 = 0$, we have
    \begin{equation}\label{bigeq}
        \|\nabla f(x_1) - \nabla f(x_2)\| \; \leq \; L \, (\|a_1 - a_2\| + \left\|b_1 - b_2\right\| + \left\|B_1 - B_2\right\| + \left\|A_1 - A_2\right\|),
    \end{equation}
    with~$a_i = \nabla_x f_u (x_i,y(x_i))$, $b_i = \nabla_y f_u (x_i,y(x_i))$, $B_i = \nabla_{yy}^2 f_{\ell} (x_i,y(x_i))$, and~$A_i = \nabla_{xy}^2 f_{\ell} (x_i,y(x_i))$, $i \in \{1, 2\}$.
    
    By the Lipschitz continuity of $\nabla_x f_u$ in $(x,y)$ due to Assumption~\ref{ass:smoothness}, we have
    \begin{alignat*}{2}
        & \|a_1 - a_2\| && \; = \; \|\nabla_x f_u (x_1,y(x_1)) - \nabla_x f_u (x_2,y(x_2))\| \; \leq \; L_1 \| (x_1-x_2, y(x_1) - y(x_2) )^{\top}\|,
    \end{alignat*}
    where~$L_1$ denotes the Lipschitz constant.
    Squaring both sides and using~\eqref{ineq2}, we obtain
    \begin{alignat}{2}
    & \|a_1 - a_2\|^2 && \; \leq \; 
    L_1^2 (\|x_1-x_2\|^2 + L_{\ell}^2\|x_1 - x_2\|^2) \; = \; L_1^2(1+L_{\ell}^2) \|x_1 - x_2\|^2. \nonumber
    \end{alignat}
    Taking the square root of both sides yields $\|a_1 - a_2\| \leq L_a\|x_1 - x_2\|$, where $L_a = L_1(1+L_{\ell}^2)^{\frac{1}{2}}$. Since~$\nabla_y f_{u}$ is Lipschitz continuous in~$(x,y)$, after performing the same process as above, we will obtain the following bound:~$\left\|\nabla_y f_u (x_1,y(x_1)) - \nabla_y f_u (x_2,y(x_2))\right\| \leq L_b \|x_1 - x_2\|$. Similarly, since~$\nabla_{yy}^2 f_{\ell}$ and~$\nabla_{xy}^2 f_{\ell}$ are Lipschitz continuous in~$(x,w)$, we have~$\|\nabla_{yy}^2 f_{\ell}(x_1,y(x_1)) - \nabla_{yy}^2 f_{\ell}(x_2,y(x_2))\| \leq L_B \|x_1 - x_2\|$ and~$\left\|\nabla_{xy}^2 f_{\ell}(x_1,y(x_1)) - \nabla_{xy}^2 f_{\ell}(x_2,y(x_2))\right\| \leq L_A \|x_1 - x_2\|$.
    Thus, substituting all of these into \eqref{bigeq}, we obtain
    \[
    \|\nabla f(x_1) - \nabla f(x_2)\| \; \leq \; L_{\nabla f} \|x_1 - x_2\|, \nonumber 
    \]
    where $L_{\nabla f} = L \, (L_a + L_b + L_A + L_B)$. This concludes the proof for the~LL unconstrained case. 
    
    The proof for the~LL constrained case follows very similar steps. However, given the structure of the adjoint gradient~(\ref{constr_direction}) in the constrained case, we must first establish the Lipschitz continuity in~$x$ of $\nabla_v G$ and $\nabla_x G$ given in~(\ref{eq:jacobians}). At this point of the paper, this follows from already-seen arguments, which we will briefly summarize here to avoid repetition.
    The Lipschitz continuity of the Hessian terms $\nabla_{yy}^2 \mathcal{L}_\ell$ and $\nabla_{yx}^2 \mathcal{L}_\ell$ results from the Lipschitz continuity of the Hessians defining the problem, the Lipschitz continuity of the multipliers~(\ref{ineq2}), the Lipschitz continuity of sums and products, and the boundedness of all terms by Assumption~\ref{ass:smoothness-constrained}.
    The Lipschitz continuity of the terms $z_I \circ \nabla_x c_I^{\top}$ and $z_I \circ \nabla_y c_I^{\top}$ results from the Lipschitz continuity of the Jacobians, the Lipschitz continuity of the multipliers~(\ref{ineq2}), the Lipschitz continuity of sums and products, and the boundedness of all terms.
    The remaining elements in~(\ref{eq:jacobians}) are constraint functions and their Jacobians, which are Lipschitz continuous per Assumption~\ref{ass:smoothness-constrained}.
\end{proof}

\section{Theorem~\ref{theorem:conv_2}}\label{appendix:conv_2}

	\begin{proof}
    For any $k \in \mathbb{N}$, we can write
    \begin{align}
        \hspace{-0.3cm} \mathbb{E}_{\xi_k^{\all}}[\Vert x_{k+1} - x_*\Vert^2] & = \mathbb{E}_{\xi_k^{\all}}[\Vert P_X (x_k + \alpha_k \, d(x_k,\tilde w_k,\xi_k)) - x_*\Vert^2] \nonumber\\
        & \le \mathbb{E}_{\xi_k^{\all}}[\Vert x_k + \alpha_k \, d(x_k,\tilde w_k,\xi_k) - x_*\Vert^2] \nonumber\\
        & = \Vert x_k - x_*\Vert^2 + \alpha_k^2 \mathbb{E}_{\xi_k^{\all}}[\Vert d(x_k,\tilde w_k,\xi_k) \Vert^2] \nonumber\\
        &\quad + 2\alpha_k \mathbb{E}_{\xi_k^{\all}}[d(x_k,\tilde w_k,\xi_k)]^\top (x_k - x_*). \nonumber
    \end{align}
    Adding and subtracting the term $2\alpha_k (\mathbb{E}_{\xi_k^{\all}}[d(x_k,w(x_k))])^\top(x_k - x_*)$, noting that $d(x_k,w(x_k)) = -\nabla f(x_k)$, and applying the Cauchy-Schwarz and Jensen's inequalities, we obtain
    \begin{align}
        \hspace{-0.3cm} \mathbb{E}_{\xi_k^{\all}}[\Vert x_{k+1} - x_*\Vert^2]
        &\le \Vert x_k - x_*\Vert^2 + \alpha_k^2 \mathbb{E}_{\xi_k^{\all}}[\Vert d(x_k,\tilde w_k,\xi_k)\Vert^2] - 2\alpha_k \nabla f(x_k)^\top (x_k - x_*) \nonumber\\
        & \quad + 2\alpha_k \mathbb{E}_{\xi_k^{\all}}[d(x_k,\tilde w_k,\xi_k) - d(x_k,w(x_k))]^\top (x_k - x_*) \nonumber \\
        &\le \Vert x_k - x_*\Vert^2 + \alpha_k^2 \mathbb{E}_{\xi_k^{\all}}[\Vert d(x_k,\tilde w_k,\xi_k)\Vert^2] - 2\alpha_k \nabla f(x_k)^\top (x_k - x_*) \nonumber\\
        & \quad + 2\alpha_k \mathbb{E}_{\xi_k^{\all}}[\Vert d(x_k,\tilde w_k,\xi_k) - d(x_k, w(x_k)) \Vert] \Vert x_k - x_* \Vert. \nonumber 
    \end{align}
     Then, by using Assumption~\ref{ass:boundedness} and inequalities~\eqref{eq501}, \eqref{eq:5555}, and \eqref{eq:strong_convexity_3},
    \begin{equation}\label{eq20}
       \mathbb{E}_{\xi_k^{\all}}[\Vert x_{k+1} - x_*\Vert^2]
        \; \le \; (1 - 2 c\alpha_k)\Vert x_{k} - x_*\Vert^2 + (G_d + 2 C_d \Theta) \,  \alpha_k^2.\nonumber
    \end{equation}
    Denoting $M = G_d + 2 C_d \Theta$ and taking the total expectation on both sides, one obtains
    \begin{align}
       \mathbb{E}[\Vert x_{k+1} - x_*\Vert^2] 
        \; \le \; (1 - 2c\alpha_k) \mathbb{E}[\Vert x_{k} - x_*\Vert^2] + M \alpha_k^2.\nonumber
    \end{align}
    Using $\alpha_k = \frac{\gamma}{k}$ for some constant $\gamma > \frac{1}{2 \, c}$, it follows by induction \cite[Eq. (2.9) and (2.10)]{ANemirovski_etal_2009} that
    \begin{align}\label{eq:444}
       \mathbb{E}[\Vert x_{k} - x_*\Vert^2] 
        \; \le \; \frac{\max\{2 \, \gamma^2 M (2 c \gamma - 1)^{-1},\| x_0 - x_*\|^2\}}{k},
    \end{align}
    which proves the first result. 
    
    From \eqref{ass_37_result}, one obtains (see, e.g.,~\cite[Lemma 5.7]{ABeck_2017} and~$P_X(x_k - x_*) = (x_k - x_*)$) 
    \begin{equation}\label{eq:333}
        f(x_k) \; \le \; f(x_*) + (P_X\nabla f(x_*))^\top (x_k - x_*) + \frac{1}{2} L_{\nabla{f}} \Vert x_k - x_*\Vert^2.
    \end{equation}
    From \eqref{eq:444} and \eqref{eq:333}, by taking the total expectation and recalling~$P_X\nabla f(x_*) = 0$, one can obtain the optimality gap in terms of function values, yielding
    \begin{alignat*}{2}
      & \mathbb{E}[f(x_k)] - f(x_*) 
        && \; \le \; \frac{1}{2} L_{\nabla {f}} \mathbb{E}[\Vert x_k - x_*\Vert^2] \\ & && \; \le \; \frac{(L_{\nabla {f}}/2) \max\{2 \, \gamma^2 M (2 c \gamma - 1)^{-1},\| x_0 - x_*\|^2\}}{k}.\nonumber
    \end{alignat*}
    \end{proof}

\section{Theorem~\ref{theorem:conv_convex}}\label{appendix:conv_convex}

\begin{proof}
    Assumption~\ref{ass:convexity} implies that
    \begin{equation}\label{eq200}
         \nabla f(x_k)^\top (x_k - x_*) \; \ge \; f(x_k) - f(x_*).
    \end{equation}
    Repeating the same arguments that in the proof of Theorem~\ref{theorem:conv_2} led to~\eqref{eq20}, but now using~\eqref{eq200} instead, we obtain
    \begin{align}
        \mathbb{E}_{\xi_k^{\all}}[\Vert x_{k+1} - x_*\Vert^2]
        \;\leq\; \Vert x_{k} - x_*\Vert^2 + 2\alpha_k(f(x_*) - f(x_k)) + (G_d + 2 C_d \Theta)\alpha_k^2. \nonumber
    \end{align}
    Letting $M = G_d + 2 C_d \Theta$, we have
    \[
    \mathbb{E}_{\xi_k^{\all}}\left[ \|x_{k+1} - x_* \|^2\right] \;\leq\; \|x_{k} - x_*\|^2 + 2\alpha_k (f(x_*) - f(x_k)) + M\alpha_k^2.
    \]
    Rearranging, taking total expectations, and dividing by $\alpha_k$, we obtain
    \[
    2(\mathbb{E}[f(x_k)] - f(x_*)) \;\leq\; \frac{\mathbb{E}\left[\|x_{k} - x_*\|^2\right]}{\alpha_k} - \frac{\mathbb{E}\left[ \|x_{k+1} - x_* \|^2\right]}{\alpha_k} + M\alpha_k.
    \]
    If we replace~$k$ by~$s$ and sum over $s=0,1,\ldots,k$, we obtain
    \begin{align}
        2\sum_{s=0}^k(\mathbb{E}[f(x_s)] - f(x_*))
        &\; \leq \;\sum_{s=0}^k\left(\frac{\mathbb{E}\left[\|x_{s} - x_*\|^2\right]}{\alpha_s} - \frac{\mathbb{E}\left[ \|x_{s+1} - x_* \|^2\right]}{\alpha_s}\right) + M\sum_{s=0}^k\alpha_s \nonumber\\
        &\;= \;\frac{\mathbb{E}\left[\|x_{0} - x_*\|^2\right]}{\alpha_0} + \sum_{s=1}^k\left( \frac{1}{\alpha_s} - \frac{1}{\alpha_{s-1}} \right)\mathbb{E}\left[\|x_{s} - x_*\|^2\right] + M\sum_{s=0}^k\alpha_s \nonumber\\
        &\;= \;\frac{\mathbb{E}\left[\|x_{k} - x_*\|^2\right]}{\alpha_k} + M\sum_{s=0}^k\alpha_s. \nonumber 
    \end{align}
    Using Assumption~\ref{ass:boundedness} and repeating the same steps used in~\cite[Theorem 5.3]{SLiu_LNVicente_2019}, we can obtain the desired result.
    \end{proof}
    
\end{document}